\documentclass[12pt]{article}

\usepackage{listings}
\usepackage[margin=1in]{geometry}
\usepackage{amsmath}
\usepackage{bm}
\usepackage{url}
\usepackage{xcolor}
\usepackage{graphicx}
\graphicspath{{./figures/}}
\usepackage{booktabs}
\usepackage{xspace}
\usepackage{amssymb} 
\usepackage[numbers,sort&compress]{natbib} 
\usepackage{hyperref} 
\usepackage[acronym,toc]{glossaries}
\newacronym{SUPG}{SUPG}{Streamline-Upwind Petrov--Galerkin}
\newacronym{PSPG}{PSPG}{Pressure-Stabilized Petrov--Galerkin}
\newacronym{MMS}{MMS}{Method of Manufactured Solutions}
\newacronym{LSIC}{LSIC}{Least-Squares Incompressiblity Constraint}
\newacronym{INS}{INS}{Incompressible Navier--Stokes}
\newacronym{MOOSE}{MOOSE}{Multiphysics Object Oriented Simulation Environment}
\newacronym{LBB}{LBB}{Ladyzhenskaya--Babu\v{s}ka--Brezzi}
\newacronym{VMS}{VMS}{Variational Multiscale}
\newacronym{SGS}{SGS}{Subgrid Scale}
\newacronym{GLS}{GLS}{Galerkin Least-Squares}
\newacronym{AMG}{AMG}{Algebraic Multigrid}
\newacronym{GMG}{GMG}{Geometric Multigrid}
\newacronym{ODE}{ODE}{ordinary differential equation}
\newacronym{DG}{DG}{discontinuous Galerkin}
\newacronym{PDE}{PDE}{partial differential equations}
\newacronym{FEM}{FEM}{finite element method}
\makeglossaries

\usepackage{authblk}

\usepackage{subcaption}

\usepackage{upquote}

\definecolor{codegreen}{rgb}{0.0, 0.5, 0.0} 
\definecolor{codegray}{rgb}{0.5,0.5,0.5} 
\definecolor{codepurple}{rgb}{0.58,0,0.82}
\definecolor{backcolour}{rgb}{0.9,0.9,0.9} 
\definecolor{blue(pigment)}{rgb}{0.2, 0.2, 0.6}
\definecolor{bleudefrance}{rgb}{0.19, 0.55, 0.91}
\definecolor{bluegray}{rgb}{0.4, 0.6, 0.8}
\lstdefinestyle{mystyle}{
    backgroundcolor=\color{backcolour},
    commentstyle=\color{codegreen},
    numberstyle=\tiny\color{codegray},
    stringstyle=\color{codepurple},
    basicstyle=\scriptsize\ttfamily,
    breakatwhitespace=false,
    breaklines=true,
    captionpos=b,
    keepspaces=true,
    numbers=left,
    numbersep=5pt,
    showspaces=false,
    showstringspaces=false,
    showtabs=false,
    tabsize=2
}
\newcommand{\us}{\char`_}
\newcommand{\unit}[1]{\hat{#1}}
\newcommand{\ns}{\texttt{navier{\us}stokes}\xspace}
\newcommand{\heatcond}{\texttt{heat{\us}conduction}\xspace}
\newcommand{\tm}{\texttt{tensor{\us}mechanics}\xspace}
\newcommand{\pf}{\texttt{phase{\us}field}\xspace}
\newcommand{\grad}{\text{grad}\,}
\newcommand{\curl}{\text{curl}\,}
\newcommand{\Div}{\text{div}\,}
\newcommand{\Lap}{\text{Lap}\,}
\newcommand{\meshspace}{\mathcal{G}}
\newcommand{\elem}{\Omega_e}
\newcommand{\pp}[2]{$P^{#1}P^{#2}$}
\newcommand{\qq}[2]{$Q^{#1}Q^{#2}$}
\newcommand{\hcomp}[2]{{#1}^h_{#2}}
\newcommand{\dirichletvec}{\vec{g}}
\newcommand{\neumannscalar}{s}
\newcommand{\neumannvec}{\vec{\neumannscalar}}
\newcommand{\mat}[1]{\bm{#1}}
\newcommand{\highlight}[1]{%
  \colorbox{red!50}{$\displaystyle#1$}}
\newcommand{\cpp}{C{\tiny{$^{++}$}}\xspace}
\newcommand{\freefem}{FreeFem{\tiny{$^{++}$}}\xspace}

\title{Overview of the Incompressible Navier--Stokes simulation capabilities in the MOOSE Framework}
\author{John~W.~Peterson}
\author{Alexander~D.~Lindsay}
\author{Fande Kong}
\affil{ \normalsize\texttt{\{jw.peterson,alexander.lindsay,fande.kong\}@inl.gov} }
\affil{Department of Modeling and Simulation\\Idaho National Laboratory\\2525 N. Fremont Ave.\\Idaho Falls, ID 83415}
\date{\today}

\begin{document}
\maketitle

\begin{abstract}
The \gls{MOOSE} framework is a high-performance, open source, \cpp finite
element toolkit developed at Idaho National Laboratory. \gls{MOOSE}
was created with the aim of assisting domain scientists and engineers
in creating customizable, high-quality tools for multiphysics
simulations.  While the core \gls{MOOSE} framework itself does not
contain code for simulating any particular physical application, it is
distributed with a number of physics ``modules'' which are tailored to
solving e.g.\ heat conduction, phase field, and solid/fluid mechanics
problems. In this report, we describe the basic equations, finite
element formulations, software implementation, and
regression/verification tests currently available in \gls{MOOSE}'s \ns
module for solving the \gls{INS} equations.
\end{abstract}

\section{Introduction}
The \gls{MOOSE} framework~\cite{Gaston_2015}, which has been under
development since 2008 and LGPL-2.1 licensed open source software on
GitHub\footnote{\url{www.github.com/idaholab/moose}} since March
2014~\cite{Slaughter_2015}, was originally created to facilitate the
development of sophisticated simulation tools by domain experts in
fields related to nuclear power generation (neutron transport, nuclear
fuel performance, mesoscale material modeling, thermal hydraulics,
etc.) who were not necessarily experts in computational science.
A guiding principle in the development of the \gls{MOOSE} framework is:
by lowering the typical computational science-related barriers of entry,
e.g.\ knowledge of low-level programming languages, familiarity with
parallel programming paradigms, and complexity/lack of interfaces to
existing numerical software libraries, one can substantially increase the
number and quality of scientific contributions to the aforementioned fields
in a cost-effective, sustainable, and scientifically sound manner.

The generality and applicability of \gls{MOOSE} stems directly from
that of the \gls{FEM} itself. \gls{MOOSE} relies on the
libMesh~\cite{Kirk_2006} library's mesh handling, geometric/finite
element families, and numerical linear algebra
library~\cite{PETSc_2016} interfaces to avoid code duplication, and
focuses on providing a consistent set of interfaces which allow users
to develop custom source code and input file syntax into
application-specific simulation tools.  In the research and
development environment, one frequently simulates physics involving the
combination of a well-established mathematical model (heat conduction,
solid mechanics, etc.) alongside a new or as-yet-untested collection
of models.  The \gls{MOOSE} physics modules were created to help
ensure that these types of combinations were not only possible, but
also relatively painless to create.

There are currently in the neighborhood of a dozen actively developed \gls{MOOSE} physics modules.
The \heatcond module, for example, can be used to solve the transient solid heat conduction equation,
\begin{align}
  \label{eq:heat_conduction}
  \rho c_p \frac{\partial T}{\partial t} &= \nabla \cdot k \nabla T + f
\end{align}
where $\rho$, $c_p$, and $k$ are the density, specific heat, and thermal conductivity, respectively,
of the medium, $T$ is the unknown temperature, and $f$ is a volumetric heat source term.
This module supports both single- and multi-domain solves
using Lagrange multiplier and gap heat transfer finite element formulations. The \tm module
focuses on solving linear elastic and finite strain solid mechanics equations
of the form
\begin{align}
  \nabla \cdot \left(\bm{\sigma} + \bm{\sigma}_0\right) + \vec{b} = \vec{0}
\end{align}
where $\bm{\sigma}$ is the Cauchy stress tensor, $\bm{\sigma}_0$ is an ``extra''
stress due to e.g.\ coupling with other physics, and $\vec{b}$ is a possibly
temperature-dependent body force which can be used to couple together the
\tm and \heatcond modules. Likewise, the
\pf module can be used to solve the Allen--Cahn~\cite{Allen_1972}
and Cahn--Hilliard~\cite{Cahn_1958,Zhang_2013} equations:
\begin{align}
  \label{eq:ac}
  \frac{\partial \eta_j}{\partial t} &= - L_j \frac{\delta F}{\delta \eta_j}
  \\
  \label{eq:ch}
  \frac{\partial c_i}{\partial t} &= \nabla \cdot \left(M_i \nabla \frac{\delta F}{\delta c_i}\right)
\end{align}
for the non-conserved ($\eta_j$) and conserved ($c_i$) order parameters,
where the free energy functional $F$ can depend on both the temperature $T$ and stress field $\bm{\sigma}$
computed by the \heatcond and \tm modules, respectively.

\gls{MOOSE}'s \ns module, which is the subject of the present work, is capable
of solving both the compressible and incompressible Navier--Stokes
equations using a variety of Petrov--Galerkin, \gls{DG}, and finite
volume (implemented as low-order \gls{DG})
discretizations. Here, we focus only on the incompressible version of
the Navier--Stokes equations, since the numerical techniques used for
the compressible version are significantly different and will be
discussed in a future publication. One of the more compelling reasons
for developing an \gls{INS} solver capability as a \gls{MOOSE}
module is the potential for coupling it to the \gls{MOOSE} modules just
listed, as well as others: fluid-structure interaction problems could eventually
be tackled via coupling to the \texttt{tensor{\us}mechanics} module,
immiscible fluid interface topology tracking with realistic divergence-free
flow fields could be implemented via coupling to the \texttt{level{\us}set}
module, and chemically-reacting flows could be investigated by coupling
to the \texttt{chemical{\us}reactions} module.

The \gls{INS} equations are an important tool for investigating many
types of physical applications where viscous effects cannot readily be
neglected, including aerodynamic drag calculations, turbulence
modeling, industrial lubrication, boundary layers, and pipe flows, to
name just a few. A number of classic
textbooks~\cite{White_1991,White_2010,Panton_2012,Munson_2012} have
been written on the historical development and present-day
understanding of incompressible flow models.
The widespread use of the \gls{FEM} for simulating incompressible
flow problems postdates its use as the tool of choice for structural analysis
problems, but relatively recent theoretical, numerical, and computational
advances from the 1980s to the present day have led to a dramatic increase in its
popularity. The numerical methods employed in the \ns module generally follow
well-established trends in the field at large and are based largely on
the many reference texts~\cite{Carey_1986,Gallagher_1988,Gresho_1998,Gresho_2002,Lewis_2004,Zienkiewicz_2005}
available in this area. The intent of this article is therefore not to
describe new developments in the field of finite element modeling for
incompressible fluids, but rather to elucidate an open, extensible,
community-developed approach for applying established and proven techniques.

The rest of this paper is arranged as follows: in \S\ref{sec:goveqn} the relevant
forms of the governing equations are briefly described, and in \S\ref{sec:numerical_method}
details of the numerical method employed in the \ns module are discussed, including
weak formulations, stabilized \gls{FEM}, time discretizations, and
solvers. In \S\ref{sec:software} we describe in more detail the overarching design
of the \gls{MOOSE} framework itself and the implications of that design in the
development of the \ns module. In \S\ref{sec:verification_tests},
several of the verification tests which form the backbone of our regression
testing and continuous integration based software development process are characterized, and finally,
in \S\ref{sec:results} we discuss some representative applications
which demonstrate the effectiveness and usefulness of the software.

\section{Governing equations\label{sec:goveqn}}
The \ns module in \gls{MOOSE} can be used to solve the \gls{INS}
equations of fluid flow in 2D and 3D Cartesian coordinates as well as
2D-axisymmetric (referred to as ``RZ'' in the code) coordinates. Two variants
of the governing equations, the so-called ``traction'' and
``Laplace'' forms, are supported. The two forms differ
primarily in the constant viscosity assumption, and the resulting form of
the natural boundary conditions~\cite{Taylor_1984}. Detailed descriptions are provided
in \S\ref{sec:traction_form} and~\S\ref{sec:laplace_form}, respectively.

\subsection{Traction form\label{sec:traction_form}}
The ``traction'' form of the \gls{INS} equations, so-called because the natural boundary condition
involves the surface traction vector, can be expressed in the following compact form:
\begin{align}
  \label{eqn:ns}
  \vec{R}(\vec{U}) = \vec{0}
\end{align}
where $\vec{U} \equiv (\vec{u}, p)$ is a vector composed of the velocity and pressure unknowns, and
\begin{align}
  \label{eqn:ns_mom_and_mass}
  \vec{R}(\vec{U}) \equiv
  \begin{bmatrix}
  \rho\left(\displaystyle \frac{\partial \vec{u}}{\partial t} + (\vec{u} \cdot \grad\!)\vec{u}\right) - \Div \bm{\sigma} - \vec{f}
  \\[12pt]
  \Div \vec{u}
  \end{bmatrix}
\end{align}
where $\rho$ is the (constant) fluid density,
$\vec{f}$ is a body force per unit volume, and $\bm{\sigma}$ is the total stress
tensor, defined by:
\begin{align}
  \nonumber
  \bm{\sigma} &\equiv -p\bm{I} + \bm{\tau}
  \\
  \label{eqn:stress_tensor}
              &= -p\bm{I} + \mu \left( \grad \vec{u} + (\grad \vec{u})^T \right)
\end{align}
In~\eqref{eqn:stress_tensor}, $\mu$ is the dynamic viscosity,
$\bm{I}$ is the identity matrix, and $\bm{\tau}$ is the viscous stress tensor.
In the present work, we consider only Newtonian fluids, but the extension of
our formulation to other types of fluids is straightforward.
In~\eqref{eqn:stress_tensor}, we have \emph{a priori} neglected the
dilatational viscosity contribution, $\lambda
\left(\Div\vec{u}\right) \bm{I}$, where the second coefficient of viscosity $\lambda = -\frac{2\mu}{3}$~\cite{Gad-el-Hak_1995}.
This term is dropped based on the divergence-free
constraint (second row of~\eqref{eqn:ns_mom_and_mass}), but it could be retained, and
in fact it is connected to a particular form of numerical stabilization
that will be discussed in \S\ref{sec:lsic_stabilization}. Finally, we note that while
variants of the convective operator in~\eqref{eqn:ns_mom_and_mass},
including $\curl \vec{u} \times \vec{u} + \frac{1}{2}\grad |\vec{u}|^2$ (rotation form)
and $\frac{1}{2} \left((\vec{u} \cdot \grad\!)\vec{u} + \Div \vec{u} \vec{u} \right)$ (skew-symmetric form),
can be useful under certain circumstances~\cite{Zang_1991},
only the form shown in~\eqref{eqn:ns_mom_and_mass} is implemented in the \ns module.

\subsection{Laplace form\label{sec:laplace_form}}
As mentioned previously, the density $\rho$ of an incompressible fluid
is, by definition, a constant, and therefore can be scaled out of the
governing equations without loss of generality. In the \ns module, however, $\rho$
is maintained as an independent, user-defined constant and the
governing equations are solved in dimensional form. Although there
are many benefits to using non-dimensional formulations, it is
frequently more convenient in practice to couple with other simulation codes
using a dimensional formulation, and this motivates the present approach.

The remaining material property, the dynamic viscosity $\mu$, is also
generally constant for isothermal, incompressible fluids, but may
depend on temperature, especially in the case of gases. In the \ns
module, we generally allow $\mu$ to be non-constant
in~\eqref{eqn:stress_tensor}, but if $\mu$
is constant, we can make the following simplification:
\begin{align}
  \nonumber
  \Div \bm{\sigma} &= \Div \left( -p\bm{I} + \mu \left( \grad \vec{u} + (\grad \vec{u})^T \right) \right)
  \\
  \nonumber
                   &= - \grad p + \mu \left( \Div \grad \vec{u} + \Div (\grad \vec{u})^T \right)
  \\
  \nonumber
                   &= - \grad p + \mu \left( \grad (\Div \vec{u}) + \Div (\grad \vec{u})^T \right)
  \\
  \nonumber
                   &= - \grad p + \mu \,  \Div (\grad \vec{u})^T
  \\
  \label{eqn:laplace_derivation}
                   &= - \grad p + \mu \, \Lap \vec{u}
\end{align}
where $\Lap \vec{u} \equiv \Div (\grad \vec{u})^T$ is the
\emph{vector} Laplacian operator (which differs from component-wise application of
the scalar Laplacian operator in e.g.\ cylindrical coordinates). Note that
in the derivation of~\eqref{eqn:laplace_derivation} we also assumed
sufficient smoothness in $\vec{u}$ to interchange the order of
differentiation in the third line, and we used the incompressibility
constraint to drop the first term in parentheses
on line 4.  The ``Laplace'' form of~\eqref{eqn:ns} is therefore given by:
\begin{align}
  \label{eqn:ns_laplace}
  \vec{R}_L(\vec{U}) = \vec{0}
\end{align}
where
\begin{align}
  \label{eqn:ns_mom_and_mass_laplace}
  \vec{R}_L(\vec{U}) \equiv
  \begin{bmatrix}
  \rho\left(\displaystyle \frac{\partial \vec{u}}{\partial t} + (\vec{u} \cdot \grad\!)\vec{u}\right) + \grad p - \mu \,  \Lap \vec{u}  - \vec{f}
  \\[12pt]
  \Div \vec{u}
  \end{bmatrix}
\end{align}
In the continuous setting, Eqn.~\eqref{eqn:ns_laplace} is equivalent
to~\eqref{eqn:ns} provided that the constant viscosity assumption is satisfied.
In the discrete setting, when the velocity field is
not necessarily pointwise divergence-free, the approximate solutions resulting from finite element
formulations of the two equation sets will not necessarily be the same. In practice, we have
observed significant differences in the two formulations near outflow boundaries where
the natural boundary condition is weakly imposed. The \ns module therefore implements
both~\eqref{eqn:ns} and~\eqref{eqn:ns_laplace},
and the practitioner must determine which is appropriate in their particular application.

Equations~\eqref{eqn:ns} and~\eqref{eqn:ns_laplace} are valid in all of the
coordinate systems mentioned at the beginning of this section, provided that the proper definitions of the
differential operators ``grad'' and ``div'' are employed.
The governing equations are solved in the
spatiotemporal domain $(\vec{x},t) \in \Omega \times (0, T)$ where
$\Omega$ has boundary $\Gamma$ and outward unit normal $\unit{n}$, and
the full specification of the problem is of course incomplete without the addition of
suitable initial and boundary conditions. In the \ns module, we
support the usual Dirichlet ($\vec{u} = \dirichletvec \in
\Gamma_D\subset\Gamma$) and surface traction ($\unit{n} \cdot
\bm{\sigma} = \neumannvec \in \Gamma_N\subset\Gamma$) boundary
conditions, as well as several others that
are applicable in certain specialized situations. Initial conditions on the
velocity field, which must satisfy the divergence-free constraint,
are also required, and may be specified as constants, data, or in functional form.

\section{Numerical method\label{sec:numerical_method}}
The \ns module solves the traction~\eqref{eqn:ns} and
Laplace~\eqref{eqn:ns_laplace} form of the \gls{INS} equations, plus initial
and boundary conditions, using a stabilized Petrov--Galerkin \gls{FEM}.
The development of this \gls{FEM} follows
the standard recipe: first, weak formulations are developed for
prototypical boundary conditions in~\S\ref{sec:weak_formulation}, then
the semi-discrete stabilized finite element formulation is introduced
in~\S\ref{sec:finite_element_method}, followed by the time
discretization in~\S\ref{sec:time_discretization}. Finally,
in~\S\ref{sec:preconditioning} some details on the different
types of preconditioners and iterative solvers available in the
\ns module are provided.

\subsection{Weak formulation\label{sec:weak_formulation}}
To develop a weak formulation for~\eqref{eqn:ns}, we
consider the combined Dirichlet/flux boundary and initial conditions
\begin{alignat}{2}
  \label{eqn:dirichlet_bc}
  \vec{u} &= \dirichletvec && \; \in \Gamma_D
  \\
  \label{eqn:flux_bc}
  \unit{n} \cdot \bm{\sigma} &= \neumannvec && \; \in \Gamma_N
  \\
  \label{eqn:ic}
  \vec{u} &= \vec{u}_0  && \; \in \Omega \text{ at } t=0
\end{alignat}
where $\dirichletvec$ and $\neumannvec$ are given (smooth) boundary data, and $\Gamma
= \Gamma_D \cup \Gamma_N$ is the entire boundary of $\Omega$.
The weak formulation of~\eqref{eqn:ns} with
conditions~\eqref{eqn:dirichlet_bc}--\eqref{eqn:ic} is then:
find $\vec{U} \in \mathcal{S}$ such that
\begin{align}
  \label{eqn:a}
  a(\vec{U},\vec{V}) = 0
\end{align}
holds for every $\vec{V} \in \mathcal{V}$, where
\begin{align}
  \nonumber
  a(\vec{U},\vec{V})
  &\equiv
  \int_{\Omega} \left[
  \rho\left(\frac{\partial \vec{u}}{\partial t} + (\vec{u} \cdot \grad\!)\vec{u}\right)\cdot \vec{v}
  -p\, \Div \vec{v} + \bm{\tau}:\grad \vec{v} - \vec{f}\cdot\vec{v} \right] \text{d} \Omega
  \\
  \label{eqn:a_definition}
  & \qquad -\int_{\Gamma_N} \neumannvec \cdot \vec{v}  \;\text{d}\Gamma -
    \int_{\Omega} q\, \Div \vec{u} \; \text{d} \Omega
\end{align}
and the relevant function spaces are
\begin{align}
  \mathcal{S} &= \{ \vec{U} =(\vec{u}, p) : u_i,p \in H^1(\Omega) \times C^0(0,T),\, \vec{u}=\dirichletvec \in \Gamma_D\}
  \\
  \mathcal{V} &= \{ \vec{V} =(\vec{v}, q) : v_i,q \in H^1(\Omega) \times C^0(0,T),\, \vec{v}=\vec{0} \in \Gamma_D\}
\end{align}
where $\Omega \subset \mathbb{R}^{n_s}$, $n_s$ is the spatial dimension of the problem, $H^1(\Omega)$ is the Hilbert space of functions
with square-integrable generalized first derivatives defined on $\Omega$, and the colon
operator denotes tensor contraction, i.e.\ $\bm{A}\colon\!\!\bm{B} \equiv A_{ij} B_{ij}$. We remark that this weak
formulation can be derived by dotting~\eqref{eqn:ns}
with a test function $\vec{V}\equiv(\vec{v},q)$, integrating over $\Omega$, and applying the following identity/divergence theorem:
  \begin{align}
    \nonumber
    \int_{\Omega} \vec{v} \cdot \Div \bm{S} \;\text{d}{\Omega} &= \int_{\Omega} \Div (\bm{S} \vec{v}) \;\text{d}{\Omega} - \int_{\Omega} \bm{S}^T\colon \grad\vec{v} \;\text{d}{\Omega}
    \\
    \nonumber
                                                               &= \int_{\Gamma} (\bm{S} \vec{v}) \cdot\unit{n} \;\text{d}{\Gamma} - \int_{\Omega} \bm{S}^T\colon \grad\vec{v} \;\text{d}{\Omega}
    \\
                                                               &= \int_{\Gamma} (\unit{n} \cdot \bm{S} ) \cdot\vec{v} \;\text{d}{\Gamma} - \int_{\Omega} \bm{S}^T\colon \grad\vec{v} \;\text{d}{\Omega}
  \end{align}
which holds for the general (not necessarily symmetric) tensor $\bm{S}$.
In the development of the \gls{FEM}, we will extract the individual component equations of~\eqref{eqn:a_definition} by considering
linearly independent test functions $\vec{V}=(v\unit{e}_i,0)$, $i=1,\ldots,n_s$, where $\unit{e}_i$ is a unit
vector in $\mathbb{R}^{n_s}$, and $\vec{V}=(\vec{0},q)$.

For the ``Laplace'' version~\eqref{eqn:ns_laplace} of the governing equations,
we consider the same Dirichlet and initial conditions as before, and a slightly modified version of the
flux condition~\eqref{eqn:flux_bc} given by
\begin{align}
  \label{eqn:flux_bc_laplace}
  \unit{n} \cdot \left( \mu \left(\grad\vec{u}\right)^T - p\bm{I}\right) &= \neumannvec\; \in \Gamma_N
\end{align}
The weak formulation for the Laplace form of the governing equations is then:
find $\vec{U} \in \mathcal{S}$ such that
\begin{align}
  \label{eqn:aL}
  a_L(\vec{U}, \vec{V}) = 0
\end{align}
holds for every $\vec{V} \in \mathcal{V}$, where
\begin{align}
  \nonumber
  a_L(\vec{U}, \vec{V})
  &\equiv
  \int_{\Omega} \left[
  \rho\left(\frac{\partial \vec{u}}{\partial t} + (\vec{u} \cdot \grad\!)\vec{u}\right)\cdot \vec{v}
  -p\, \Div \vec{v} + \mu \left(\grad \vec{u}:\grad \vec{v}\right) - \vec{f}\cdot\vec{v} \right] \text{d} \Omega
  \\
  \label{eqn:aL_definition}
  &\qquad -\int_{\Gamma_N} \neumannvec \cdot \vec{v}  \;\text{d}\Gamma
  -\int_{\Omega} q\, \Div \vec{u} \; \text{d} \Omega
\end{align}
The main differences between~\eqref{eqn:aL_definition} and~\eqref{eqn:a_definition} are thus
the form of the viscous term in the volume integral (although both terms are symmetric) and the definition of the
Neumann data $\neumannvec$. Both \eqref{eqn:a_definition} and~\eqref{eqn:aL_definition}
remain valid in any of the standard coordinate systems, given suitable definitions for the
``grad,''  ``div,'' and integration operations.
As a matter of preference and historical convention, we have chosen the
sign of the $q\, \Div \vec{u}$ term in both~\eqref{eqn:a_definition} and~\eqref{eqn:aL_definition}
to be negative. This choice makes no difference in the weak solutions $(\vec{u},p)$ which
satisfy~\eqref{eqn:a_definition} and~\eqref{eqn:aL_definition}, but in the special case of
Stokes' flow (in which the convective term $(\vec{u}\cdot \grad\!) \vec{u}$ can be neglected),
it ensures that the stiffness matrices resulting from these weak formulations are symmetric.

\subsubsection{Boundary conditions\label{sec:bcs}}
A few notes on boundary conditions are also in order. In the case of a
pure Dirichlet problem, i.e.\ when $\Gamma_D = \Gamma$ and $\Gamma_N =
\emptyset$, the pressure $p$ in both the strong and weak formulations
is only determined up to an arbitrary constant. This indeterminacy
can be handled in a variety of different ways, for instance by
requiring that a global ``mean-zero pressure'' constraint~\cite{Bochev_2005},
$\int_{\Omega} p \;\text{d}\Omega = 0$, hold in addition to the
momentum balance and mass conservation equations discussed already, or
by using specialized preconditioned Krylov solvers~\cite{Gill_1992,Murphy_2000} for which the
non-trivial nullspace of the operator can be pre-specified.  A more
practical/simplistic approach is to specify a single value of the
pressure at some point in the domain, typically on the boundary. This
approach, which we generally refer to as ``pinning'' the pressure, is
frequently employed in the \ns module to avoid the difficulties
associated with the non-trivial nullspace of the operator.

In both~\eqref{eqn:a_definition} and~\eqref{eqn:aL_definition}, we
have applied the divergence theorem to (integrated by parts) the pressure term as well as
the viscous term. While this is the standard approach for weak
formulations of the \gls{INS} equations, it is not
strictly necessary to integrate the pressure term by parts. Not integrating the pressure term
by parts, however, both breaks the symmetry of the Stokes problem
discussed previously, and changes the form and meaning of the Neumann
boundary conditions~\eqref{eqn:flux_bc}
and~\eqref{eqn:flux_bc_laplace}. In particular, the ``combined''
boundary condition formulation of the problem is no longer sufficient
to constrain away the nullspace of constant pressure fields,
and one must instead resort to either pressure pinning or one of
the other approaches discussed previously. Nevertheless, there are
certain situations where the natural boundary condition imposed by
not integrating the pressure term by parts is useful in simulation contexts,
and therefore the \ns module provides the user-selectable boolean flag
\texttt{integrate{\us}p{\us}by{\us}parts} to control this behavior.

It is common for applications with an ``outflow'' boundary to employ
the so-called ``natural'' boundary condition, which corresponds to
setting $\neumannvec = \vec{0}$, where $\neumannvec$ is defined
in~\eqref{eqn:flux_bc} and \eqref{eqn:flux_bc_laplace} for the
traction and Laplace forms, respectively.  A drawback to imposing the
natural boundary condition is that, unless the flow is fully
developed, this condition will influence the upstream behavior of the
solution in a non-physical way.  In order to avoid/alleviate this
issue, one can artificially extend the computational domain until the
flow either becomes fully developed or the upstream influence of the
outlet boundary is judged to be acceptably small.

Another alternative is to employ the so-called ``no-BC'' boundary
condition described by Griffiths~\cite{Griffiths_1997}.  In the no-BC
boundary condition, the residual contribution associated with the
boundary term is simply computed and accumulated along with the other
terms in the weak form.  The no-BC boundary condition is therefore
distinct from the natural boundary condition, and can be thought of as
not specifying any ``data'' at the outlet, which is not strictly valid
for uniformly elliptic \gls{PDE}. Nevertheless, Griffiths showed, for the 1D scalar
advection-diffusion equation discretized with finite elements, that
the no-BC boundary condition appears to influence the upstream
behavior of the solution less than the natural boundary condition.
The \ns module includes both traction and Laplace forms of the no-BC
boundary condition via the
\texttt{INSMomentumNoBCBCTractionForm} and
\texttt{INSMomentumNoBCBCLaplaceForm} classes, respectively.

\subsection{Stabilized finite element method\label{sec:finite_element_method}}
To construct the \gls{FEM}, we discretize $\Omega$
into a non-overlapping set of elements $\mathcal{T}^h$ and
introduce the approximation spaces $\mathcal{S}^h \subset \mathcal{S}$ and
$\mathcal{V}^h \subset \mathcal{V}$, where
\begin{align}
  \label{eqn:trial_space_h}
  \mathcal{S}^h &= \{ \vec{U}_h =(\vec{u}_h, p_h) : \hcomp{u}{i} \in \meshspace^k, p_h \in \meshspace^r,\, \vec{u}_h=\dirichletvec \in \Gamma_D\}
  \\
  \label{eqn:test_space_h}
  \mathcal{V}^h &= \{ \vec{V}_h =(\vec{v}_h, q_h) : \hcomp{v}{i} \in \meshspace^k, q_h \in \meshspace^r,\, \vec{v}_h=\vec{0} \in \Gamma_D\}
\end{align}
and
\begin{align}
  \meshspace^k \equiv \{ w \in C^0(\Omega) \times C^0(0,T),\, \left. w \right|_{\elem} \in \mathcal{P}^k(\elem), \;\forall\, \elem \in \mathcal{T}^h \}
\end{align}
where $\hcomp{u}{i}$ is the $i$th component of $\vec{u}_h$,
and $\mathcal{P}^k(\elem)$ is either the space of complete polynomials of
degree $k$ (on triangular and tetrahedral elements) or
the space of tensor products of polynomials of degree $k$ (on
quadrilateral and hexahedral elements) restricted to element $\elem$.
For brevity we shall subsequently follow the standard practice and refer to
discretizations of the former type as ``\pp{k}{r} elements'' and the
latter as ``\qq{k}{r} elements.'' Although it is generally possible to employ
non-conforming velocity elements and constant or discontinuous ($\subset L^2(\Omega)$) pressure elements~\cite{Olshanskii_2002},
we have not yet implemented such discretizations in the \ns module, and do not
describe them in detail here.

The spatially-discretized (continuous in time) version of the weak formulation~\eqref{eqn:a} is
thus: find $\vec{U}_h \in \mathcal{S}^h$ such that
\begin{align}
  \label{eqn:ah}
  a(\vec{U}_h,\vec{V}_h) = 0
\end{align}
holds for every $\vec{V}_h \in \mathcal{V}^h$. The semi-discrete weak
formulation for the Laplace form~\eqref{eqn:aL} is analogous.  Not all
combinations of polynomial degrees $k$ and $r$ generate viable (i.e.\
\gls{LBB} stable~\cite{Ladyzhenskaya_1969,Babuska_1973,Brezzi_1974}) mixed
finite element formulations. A general rule of thumb is that the
velocity variable's polynomial degree must be at least one order higher than
that of the pressure, i.e.\ $k \geq r+1$; however, we will discuss
equal-order discretizations in more detail momentarily.

In addition to instabilities arising due to incompatibilities between
the velocity and pressure approximation spaces, numerical
approximations of the \gls{INS} equations can also exhibit
instabilities in the so-called advection-dominated regime. The
\gls{INS} equations are said to be advection-dominated (in the
continuous setting) when $\text{Re}\gg 1$, where
$\text{Re}\equiv\frac{\rho |\vec{u}| L}{\mu}$ is the Reynolds number
based on a characteristic length scale $L$ of the
application. Galerkin finite element approximations such
as~\eqref{eqn:ah} are known to produce highly-oscillatory solutions in
the advection-dominated limit, and are therefore of severely limited
utility in many real-world applications. To combat each of these types
of instabilities, we introduce the following stabilized weak statement
associated with~\eqref{eqn:ah}: find $\vec{U}_h \in \mathcal{S}^h$ such that
\begin{align}
  \label{eqn:stabilized_formulation}
  a(\vec{U}_h, \vec{V}_h) + \int_{\Omega'} \vec{S}(\vec{U}_h, \vec{V}_h) \cdot \vec{R}(\vec{U}_h) \;\text{d}\Omega' = 0
\end{align}
holds for every $\vec{V}_h \in \mathcal{V}^h$, where
\begin{align}
  \label{eqn:stabilization_operator}
  \vec{S}(\vec{U}_h, \vec{V}_h) \equiv
  \begin{bmatrix}
    \tau_{\text{SUPG}}\, (\vec{u}_h \cdot \grad\!) \vec{v}_h + \rho^{-1}\, \tau_{\text{PSPG}}\, \grad q_h
    \\[12pt]
    \rho\, \tau_{\text{LSIC}}\, ( \Div \vec{v}_h )
  \end{bmatrix}
\end{align}
is the stabilization operator, $\Omega'$ is the union of element interiors,
\begin{align}
  \label{eqn:element_interiors}
  \int_{\Omega'}  \;\text{d}\Omega' \equiv \sum_{\elem \in \mathcal{T}^h} \int_{\elem} \;\text{d}\elem
\end{align}
and $\tau_{\text{SUPG}}$, $\tau_{\text{PSPG}}$, and
$\tau_{\text{LSIC}}$ are mesh- and solution-dependent
coefficients corresponding to \gls{SUPG}, \gls{PSPG}, and \gls{LSIC} techniques, respectively.
The stabilized weak formulation for the Laplace
form~\eqref{eqn:ns_laplace} is analogous: $\vec{S}$ remains unchanged,
$a_L$ replaces $a$, and $\vec{R}_L$ replaces $\vec{R}$.

The development of~\eqref{eqn:stabilized_formulation} and \eqref{eqn:stabilization_operator} has
spanned several decades, and is summarized
succinctly by Donea~\cite{Donea_2003} and the references therein.
An important feature of this approach is that it is ``consistent,'' which
in this context means that, if the true solution $\vec{U}$
is substituted into~\eqref{eqn:stabilized_formulation}, then
the $\vec{R}(\vec{U})$ term vanishes, and we
recover~\eqref{eqn:ah}, the original weak statement of the problem.
A great deal of research~\cite{Brezzi_1992,Hughes_1995,Brezzi_1997,Codina_2000,Franca_2005}
has been conducted over the years
to both generalize and unify the concepts of residual-based
stabilization methods such as~\eqref{eqn:stabilized_formulation} and \eqref{eqn:stabilization_operator},
and current efforts are focused around so-called \gls{VMS} and \gls{SGS} methods.
The three distinct contributions to the stabilization operator used in
the present work are briefly summarized below.

\subsubsection{SUPG stabilization\label{sec:supg_stabilization}}
The \gls{SUPG} term introduces residual-dependent artificial
viscosity to stabilize the node-to-node oscillations present in
Galerkin discretizations of advection-dominated
flow~\cite{Brooks_1982,Johnson_1986,Shakib_1991,Tezduyar_1991,Franca_1992,Franca_1993,Tezduyar_2000}. This type of
stabilization can be employed independently of the \gls{LBB} stability
of the chosen velocity and pressure finite element spaces.
The coefficient $\tau_{\text{SUPG}}$ must have physical units of
time so that the resulting term is dimensionally consistent with the momentum
conservation equation to which it is added.

Most definitions of $\tau_{\text{SUPG}}$ stem from analysis of
one-dimensional steady advection-diffusion problems. For example, the original form
of $\tau_{\text{SUPG}}$ for the 1D linear advection-diffusion equation~\cite{Brooks_1982} is given by:
\begin{align}
  \label{eq:nodally_exact}
  \tau_{\text{SUPG,opt}} &= \frac{h}{2|\vec{u}|} \left(\coth(\text{Pe}) - \frac{1}{\text{Pe}}\right)
\end{align}
where $|\vec{u}|$ is the advective velocity magnitude, $h$ is the
element size, $\text{Pe}\equiv\frac{|\vec{u}|h}{2k}$ is the element
Peclet number, and $k$ is the
diffusion coefficient. For piecewise linear elements, this formulation for
$\tau_{\text{SUPG}}$ yields a nodally exact
solution. An alternative form of $\tau_{\text{SUPG}}$ which yields
fourth-order accuracy for the 1D linear advection-diffusion equation~\cite{Shakib_1991} is given by:
\begin{align}
  \label{eq:tau_mod}
  \tau_{\text{SUPG,mod}} = \left[
  \left(\frac{2|\vec{u}|}{h}\right)^2 +
  9\left(\frac{4k}{h^2}\right)^2
  \right]^{-\frac{1}{2}}
\end{align}
Note that neither the nodally exact~\eqref{eq:nodally_exact} nor the
fourth-order accurate~\eqref{eq:tau_mod} forms of the stabilization parameter are directly applicable to
2D/3D advection-diffusion problems or the \gls{INS} equations. Nevertheless, Eqn.~\eqref{eq:tau_mod}
does motivate the $\tau_{\text{SUPG}}$ implementation used in the \ns module, which is given by
\begin{align}
  \label{eq:tau_impl}
  \tau_{\text{SUPG}} = \alpha \left[
  \left(\frac{2}{\Delta t}\right)^2 +
  \left(\frac{2 |\vec{u}| }{h}\right)^2 +
  9\left(\frac{4\nu}{h^2}\right)^2
  \right]^{-\frac{1}{2}}
\end{align}
where $\Delta t$ is the discrete timestep size
(see~\S\ref{sec:time_discretization}), $\nu=\frac{\mu}{\rho}$ is the kinematic viscosity, and $0 \leq \alpha \leq 1$ is a
user selectable parameter which, among other things, allows the value
of $\tau_{\text{SUPG}}$ to be tuned for higher-order elements~\cite{Codina_1992}.


\subsubsection{PSPG stabilization\label{sec:pspg_stabilization}}
The \gls{PSPG} method was originally
introduced~\cite{Hughes_1986_PSPG,Hughes_1987,Douglas_1989,Hansbo_1990,Tezduyar_1992}
in the context of Stokes flow with the aim of circumventing
instabilities arising from the use of \gls{LBB}-unstable element
pairs. The motivation for its introduction stems from the fact that
\gls{LBB}-stable finite element pairs, such as the well-known ``Taylor--Hood'' (\qq{2}{1}) element,
are often more expensive or simply less convenient to work with than equal-order discretizations.
\gls{PSPG} stabilization is required for \gls{LBB}-unstable element pairs regardless of
whether the flow is advection-dominated.

In the case of Stokes flow, the stabilization parameter
$\tau_{\text{PSPG}}$ is given by:
\begin{align}
  \tau_{\text{PSPG}} = \beta\frac{h^2}{4\nu}
  \label{eq:pspg_stokes}
\end{align}
where $0 \leq \beta \leq 1$ is a user-tunable parameter. Decreasing
$\beta$ improves the convergence rate of the pressure error in the
$L^2(\Omega)$ norm, however, using too small a value for
$\beta$ will eventually cause equal-order discretizations to lose
stability. When this happens, spurious modes (such as the
``checkerboard'' mode) pollute the discrete pressure solution and
produce higher overall error.

Numerical experiments suggest
that an optimal value is $\beta\approx 1/3$. Substituting $\beta=1/3$ into
\eqref{eq:pspg_stokes} yields an expression for $\tau_{\text{PSPG}}$
equivalent to~\eqref{eq:tau_impl} in the case of steady Stokes flow. Consequently, to
generalize \gls{PSPG} stabilization to problems involving advection,
we use the same expression for both $\tau_{\text{SUPG}}$ and $\tau_{\text{PSPG}}$,
namely~\eqref{eq:tau_impl}, in the \ns module.
The coefficient $\tau_{\text{PSPG}}$ therefore also has physical units of time, and the $\rho^{-1}$
factor in~\eqref{eqn:stabilization_operator} ensures that the resulting term is dimensionally consistent
with the mass conservation equation to which it is added.


\subsubsection{LSIC stabilization\label{sec:lsic_stabilization}}
The \gls{LSIC} stabilization term~\cite{Codina_1997,Codina_2001,Olshanskii_2002,Gelhard_2005}
arises naturally in \gls{GLS} formulations of the \gls{INS} equations, and a number of authors
have reported that its use leads to improved accuracy and conditioning of the linear systems
associated with both \gls{LBB}-stable and equal-order discretizations.
As mentioned when the viscous stress tensor was first introduced
in~\eqref{eqn:stress_tensor}, a contribution of the form
$-\frac{2\mu}{3} \left(\Div\vec{u}\right) \bm{I}$ was neglected. Had
we retained this term, the weak formulation~\eqref{eqn:a_definition}
would have had an additional term of the form
\begin{align}
  \int_{\Omega} -\frac{2\mu}{3} \left(\Div \vec{u}\right) \left(\Div \vec{v}\right) \text{d}\Omega
\end{align}
Note the sign of this term, which is the same as that of the pressure term. On the other hand, the
$\tau_{\text{LSIC}}$ term implied by~\eqref{eqn:stabilized_formulation} is
\begin{align}
  \int_{\Omega'} \rho\, \tau_{\text{LSIC}} \left(\Div \vec{u}_h\right) \left(\Div \vec{v}_h\right) \text{d}\Omega'
\end{align}
Therefore, we can think of neglecting the $-\frac{2\mu}{3} \left(\Div\vec{u}\right) \bm{I}$ contribution as effectively
introducing a least-squares incompressibility stabilization contribution ``for free'' with $\tau_{\text{LSIC}} = \frac{2\nu}{3}$.

Implementing the general LSIC stabilization term gives one the
flexibility of choosing $\tau_{\text{LSIC}}$ independently of the
physical value of $\nu$, but, as discussed by
Olshanskii~\cite{Olshanskii_2002}, choosing too large a value can
cause the corresponding linear algebraic systems to become
ill-conditioned. Furthermore, it appears that the ``optimal'' value of
$\tau_{\text{LSIC}}$ depends on the \gls{LBB} stability
constant, and is therefore problem-dependent and difficult to compute
in general.  The coefficient $\tau_{\text{LSIC}}$ has physical units of (length)$^2$/time, the same as e.g.\
kinematic viscosity, and therefore the coefficient $\rho$
in~\eqref{eqn:stabilization_operator} ensures dimensional consistency
when this term is added to the momentum conservation equation.  Unlike
the \gls{PSPG} and \gls{SUPG} stabilization contributions, the
\gls{LSIC} term has not yet been implemented in the \ns module,
however we plan to do so in the near future.


\subsection{Time discretization\label{sec:time_discretization}}
To complete the description of the discrete problem, we introduce the usual polynomial basis functions
$\{ \varphi_i \}_{i=1}^N$ which span the velocity component spaces of $\mathcal{V}^h$, and $\{ \psi_i \}_{i=1}^M$
which span the pressure part of $\mathcal{V}^h$. These are the so-called ``global'' basis functions that are induced by the discretization of $\Omega$ into elements.
Taking, alternately, $\vec{V}_h = (\varphi_i \unit{e}_k, 0)$, $k=1,\ldots,n_s$ and $\vec{V}_h=(\vec{0}, \psi_i)$, the stabilized weak
formulation~\eqref{eqn:stabilized_formulation} and \eqref{eqn:stabilization_operator} leads to the component equations
\begin{align}
  \nonumber
  &\int_{\Omega'} \left[ \rho \left(\frac{\partial \hcomp{u}{k}}{\partial t}
  + \vec{u}_h \cdot \grad \hcomp{u}{k}\! \right)\varphi_i
  + \bm{\sigma}_h\colon \grad\! \left(\varphi_i \unit{e}_k\right) -f_k\varphi_i  \right] \text{d}\Omega'
    -\int_{\Gamma_N} \neumannscalar_{k} \varphi_i  \;\text{d}\Gamma
  \\
  \nonumber
  & \qquad + \int_{\Omega'}
    \tau_{\text{SUPG}} \, (\vec{u}_h \cdot \grad \varphi_i )
    \left[ \rho\left(\frac{\partial \hcomp{u}{k}}{\partial t} + \vec{u}_h \cdot \grad \hcomp{u}{k}\!\right) - \unit{e}_k \cdot \Div \bm{\sigma}_h - f_k \right]
    \text{d}\Omega'
  \\
  \label{eqn:semi_discrete_momentum_k}
  & \qquad + \int_{\Omega'} \rho\, \tau_{\text{LSIC}} \left(\Div \varphi_i \unit{e}_k \right) \left( \Div \vec{u}_h \right) \text{d}\Omega'  = 0,
    \hspace{8em}
    \left\{
    \begin{array}{c}
      i=1,\ldots,N \\
      k=1,\ldots,n_s
    \end{array}
    \right.
  \\
  \nonumber
  &\int_{\Omega'} \rho^{-1}\, \tau_{\text{PSPG}}\, \grad \psi_i \cdot \left[ \rho\left(\frac{\partial \vec{u}_h}{\partial t} + (\vec{u}_h \cdot \grad\!)\vec{u}_h\right) - \Div \bm{\sigma}_h - \vec{f}\, \right] \text{d}\Omega'
  \\
  \label{eqn:semi_discrete_mass}
  &\qquad -\int_{\Omega'} \psi_i\, \Div \vec{u}_h \;\text{d}\Omega' = 0, \hspace{16em} i=1,\ldots,M
\end{align}
where $\hcomp{u}{k}$ indicates the $k$th component of the vector
$\vec{u}_h$ (similarly for $\neumannscalar_k$ and $f_k$), and $\bm{\sigma}_h$ is the
total stress tensor evaluated at the approximate solution $(\vec{u}_h, p_h)$.
Since we are now in the discrete setting, it is also convenient to express all of the integrals as discrete sums over the elements.
Equations~\eqref{eqn:semi_discrete_momentum_k} and \eqref{eqn:semi_discrete_mass} also remain valid
in any of the standard coordinate systems since we have retained the generalized definitions
of $\Div\!$ and $\grad\!$. The component equations for the Laplace version of~\eqref{eqn:stabilized_formulation} and \eqref{eqn:stabilization_operator}
are derived analogously.

Equations~\eqref{eqn:semi_discrete_momentum_k} and \eqref{eqn:semi_discrete_mass}
are also still in semi-discrete form due to the presence of the continuous time derivative terms, so
the next step in our development is to discretize them in time using a finite difference method.
Writing $\hcomp{u}{k}(\vec{x},t_n)\equiv u_k^n$ and $p_h(\vec{x},t_n) \equiv
p^n$ for the finite element solutions at time $t=t_n$ and applying the $\theta$-method~\cite{Iserles_2009}
to~\eqref{eqn:semi_discrete_momentum_k} and \eqref{eqn:semi_discrete_mass}
results in the system of equations:
%
\begin{align}
  \label{eqn:momentum_k_fully_discrete_theta}
  \int_{\Omega'} \left[ \rho \left( \frac{u_k^{n+1} - u_k^n}{\Delta t} \right) \widetilde{\varphi}_i^{\,n+\theta} \right] \text{d}\Omega'
    + b_{k}(\vec{U}^{n+\theta}, \varphi_i) &= 0,
    \hspace{1em}
    \left\{
    \begin{array}{c}
      i=1,\ldots,N \\
      k=1,\ldots,n_s
    \end{array}
    \right.
  \\[6pt]
  \label{eqn:mass_fully_discrete_theta}
  \int_{\Omega'} \left[ \tau^{n+\theta}_{\text{PSPG}}\, \grad \psi_i \cdot \left( \frac{\vec{u}^{\,n+1} - \vec{u}^{\,n}}{\Delta t} \right)  \right] \text{d}\Omega'
  + c(\vec{U}^{n+\theta}, \psi_i) &= 0, \hspace{2.5em} i=1,\ldots,M
\end{align}
to be solved for $\vec{U}^{n+1}$, where
\begin{align}
  \widetilde{\varphi}_i^{\, n} \equiv \varphi_i + \tau^n_{\text{SUPG}} \, (\vec{u}^{\,n} \cdot \grad \varphi_i )
\end{align}
is the $i$th upwind velocity test function evaluated at time $t_n$, and
\begin{align}
  \nonumber
  b_{k}(\vec{U}^n, \varphi_i) &\equiv \int_{\Omega'}
  \Big[
  \left( \rho \vec{u}^{\,n} \cdot \grad u_k^n - f_k^n \right) \widetilde{\varphi}_i^{\, n}
  + \bm{\sigma}^n \colon \grad\! \left(\varphi_i \unit{e}_k\right)
  \Big]\, \text{d}\Omega'
  \\
  \nonumber
  &- \int_{\Omega'}
     \tau^n_{\text{SUPG}} \, (\vec{u}^{\,n} \cdot \grad \varphi_i )
     \left(  \unit{e}_k \cdot \Div \bm{\sigma}^n \right)
     \text{d}\Omega'
  \\
  \label{eqn:bi}
  &
    + \int_{\Omega'} \rho\, \tau^n_{\text{LSIC}} \left(\Div \varphi_i \unit{e}_k \right) \left( \Div \vec{u}^{\,n} \right) \text{d}\Omega'
    -\int_{\Gamma_N} \neumannscalar^n_k \varphi_i  \;\text{d}\Gamma
    \\
  \label{eqn:ci}
  c (\vec{U}^n, \psi_i) &\equiv \int_{\Omega'} \left[ \rho^{-1}\, \tau^n_{\text{PSPG}}\, \grad \psi_i \cdot \left( \rho (\vec{u}^{\,n} \cdot \grad\!)\vec{u}^{\,n} - \Div \bm{\sigma}^n - \vec{f}^{\,n} \right)
                    -\psi_i\, \Div \vec{u}^{\,n} \right] \text{d}\Omega'
\end{align}
The intermediate time $t_{n+\theta}$ and solution $\vec{U}^{n+\theta}$ are defined as
\begin{align}
  t_{n+\theta} &\equiv t_n + \theta \Delta t
  \\
  \vec{U}^{n+\theta} &\equiv \theta \vec{U}^{n+1} + (1-\theta) \vec{U}^{n}
\end{align}
Setting the parameter $\theta=1$ in~\eqref{eqn:momentum_k_fully_discrete_theta} and \eqref{eqn:mass_fully_discrete_theta} results in the
first-order accurate implicit Euler method, while setting $\theta=1/2$ corresponds
to the second-order accurate implicit midpoint method. Setting $\theta=0$ results in the explicit Euler method, but because this scheme
introduces severe stability-related timestep restrictions, we do not regularly use it in the \ns module. Other time discretizations, such
as diagonally-implicit Runge-Kutta schemes, are also available within the \gls{MOOSE} framework and can be used in the \ns module,
but they have not yet been tested or rigorously verified to work with the \gls{INS} equations.

If we express the finite element solutions at time $t_{n+1}$ as linear combinations of their respective basis functions according to:
\begin{align}
  u_k^{n+1} &= \sum_{j=1}^N u_{k_j}(t_{n+1}) \varphi_j(\vec{x}), \qquad k=1,\ldots,n_s
  \\
  p^{n+1} &= \sum_{j=1}^M p_j(t_{n+1}) \psi_j(\vec{x})
\end{align}
where $u_{k_j}(t_{n+1}),\, p_j(t_{n+1}) \in \mathbb{R}$ are unknown coefficients,
then~\eqref{eqn:momentum_k_fully_discrete_theta} and \eqref{eqn:mass_fully_discrete_theta}
represent a fully-discrete (when the integrals are approximated via numerical
quadrature) nonlinear system of equations in $O(n_s N + M)$
unknowns (depending on the boundary conditions). This nonlinear system
of equations must be solved at each timestep, typically using
some variant of the inexact Newton method~\cite{Dembo_1982}, in order to advance the approximate solution
forward in time.

\subsection{Preconditioning and solvers\label{sec:preconditioning}}
Computationally efficient implementation of the inexact Newton method
requires parallel, sparse, preconditioned iterative
solvers~\cite{Saad_2003} such as those implemented in the
high-performance numerical linear algebra library,
PETSc~\cite{PETSc_2016}. The preconditioning operator is related
to the Jacobian of the system~\eqref{eqn:momentum_k_fully_discrete_theta} and \eqref{eqn:mass_fully_discrete_theta},
and the \ns module gives the user the flexibility to choose between
computing and storing the full Jacobian matrix, or approximating its
action using the Jacobian-free Newton--Krylov~\cite{Brown_1990} method.
The full details of selecting the preconditioning method,
Krylov solver, relative and absolute convergence tolerances, etc.\ are
beyond the scope of this discussion, but the \ns module contains a
number of working examples (see~\S\ref{sec:verification_tests} and~\S\ref{sec:results}) that provide a
useful starting point for practitioners. In addition, we briefly describe
a field-split preconditioning method which is available in the \ns module below.

For convenience, let us write the discrete system of nonlinear
equations~\eqref{eqn:momentum_k_fully_discrete_theta} and \eqref{eqn:mass_fully_discrete_theta} as:
\begin{align}
  \label{eqn:nonlinear_algebra_function}
  \vec{F}(\vec{y}) = \vec{0}
\end{align}
where $\vec{F}$ is a vector of nonlinear equations in the generic unknown $\vec{y}$.
When~\eqref{eqn:nonlinear_algebra_function} is solved via Newton's
method, a linear system of equations involving the Jacobian of $\vec{F}$
must be solved at each iteration. These linear systems are typically
solved using Krylov subspace methods such as GMRES~\cite{Saad_1986,Saad_2003},
since these methods have been shown to scale well in parallel
computing environments, but the effectiveness of Krylov methods depends strongly
on the availability of a good preconditioner.

In Newton's method, the new iterate $\vec{y}^{\,(\ell+1)}$ is obtained by an update of the form
\begin{align}
  \label{eqn:nonlinear_algebra_newton_update}
  \vec{y}^{\,(\ell+1)} = \vec{y}^{\,(\ell)} + \alpha^{(\ell)} \delta \vec{y}^{\,(\ell)}
\end{align}
where $\alpha^{(\ell)} \in \mathbb{R}$ is a scale factor computed by a line
search method~\cite{Jorge_2006}. In \gls{MOOSE},
different line search schemes can be selected by setting the \texttt{line{\us}search{\us}type}
parameter to one of \{\texttt{basic}, \texttt{bt}, \texttt{cp},
\texttt{l2}\}. These options correspond to the line search options available in PETSc;
the reader should consult Section 5.2.1 of the PETSc Users Manual~\cite{PETSc_2016}
for the details of each.
The update vector $\delta \vec{y}^{\,(\ell)}$ is computed by solving the linear system
\begin{align}
  \label{eqn:nonlinear_algebra_jacobian_system}
  \mat{J}(\vec{y}^{\,(\ell)}) \delta \vec{y}^{\,(\ell)} = -\vec{F}(\vec{y}^{\,(\ell)})
\end{align}
where $\mat{J}(\vec{y}^{\,(\ell)})$ is the Jacobian of $\vec{F}$ evaluated at iterate $\vec{y}^{\,(\ell)}$.
Since the accuracy of $\mat{J}$ (we shall henceforth omit the $\vec{y}^{\,(\ell)}$ argument for simplicity) has a
significant impact on the convergence behavior of Newton's method, multiple
approaches for computing $\mat{J}$ are supported in \gls{MOOSE}.
One approach is to approximate the action of the Jacobian without explicitly
storing it. The action of $\mat{J}$ on an arbitrary vector $\vec{v}$ can
be approximated via the finite difference formula
\begin{align}
  \label{eqn:nonlinear_algebra_jacobian_free}
  \mat{J}( \vec{y} ) \vec{v} \approx \frac{\vec{F}(\vec{y} + \varepsilon \vec{v}) - \vec{F}(\vec{y})}{\varepsilon}
\end{align}
where $\varepsilon \in \mathbb{R}$ is a finite difference parameter whose value is adjusted
dynamically based on several factors~\cite{Pernice_1998}.
As mentioned previously, a good preconditioner is required for efficiently
solving~\eqref{eqn:nonlinear_algebra_jacobian_system}.  The (right-) preconditioned
version of~\eqref{eqn:nonlinear_algebra_jacobian_system}, again omitting the dependence on $\vec{y}^{\,(\ell)}$, is
\begin{align}
  \label{eqn:eqn:nonlinear_algebra_jacobian_preconditioner}
  \mat{J} \mat{B}^{-1}  \mat{B} \delta \vec{y} = -\vec{F}
\end{align}
where $\mat{B}$ is a preconditioning matrix.  When the Jacobian-free
approach~\eqref{eqn:nonlinear_algebra_jacobian_free} is employed
(\texttt{solve{\us}type = PJFNK} in \gls{MOOSE}), then $\mat{B}$ need
only be an approximation to $\mat{J}$ which is ``good enough'' to
ensure convergence of the Krylov subspace method.  When the Jacobian is
explicitly formed ({\texttt{solve{\us}type = NEWTON}} in \gls{MOOSE}), it
makes sense to simply use $\mat{B}=\mat{J}$.
We note that it is generally not necessary to explicitly form $\mat{B}^{-1}$
in order to compute its action on a vector, although this approach
is supported in \gls{MOOSE}, and may be useful in some scenarios.

While there is no single preconditioner which is ``best'' for
all \gls{INS} applications, the preconditioning approaches which are
known to work well can be broadly categorized into the ``fully-coupled''
and ``field-split'' types. These two types are distinguished by how the velocity
and pressure ``blocks'' of the Jacobian matrix are treated.
In the fully-coupled preconditioner category, there are several popular
preconditioners, such as the Additive Schwarz ({\texttt{-pc{\us}type
asm}})~\cite{Smith_2004,Kong_2016,Kong_2017}, Block Jacobi
({\texttt{-pc{\us}type bjacobi}}), and Incomplete LU ({\texttt{-pc{\us}type
ilu}})~\cite{Saad_2003} methods, which are directly accessible from within the \ns module.

Field-split preconditioners, on the other hand, take the specific
structure of $\mat{B}$ into account. The linear system~\eqref{eqn:nonlinear_algebra_jacobian_system}
has the block structure:
\begin{align}
  \label{eqn:eqn:nonlinear_algebra_field_split}
  \begin{bmatrix}
    \mat{B}_{uu} & \mat{B}_{up}
    \\[6pt]
    \mat{B}_{pu} & \mat{B}_{pp}
  \end{bmatrix}
  \begin{bmatrix}
    \delta \vec{y}_u
    \\[6pt]
    \delta \vec{y}_p
  \end{bmatrix}
  = -
  \begin{bmatrix}
    \vec{F}_u
    \\[6pt]
    \vec{F}_p
  \end{bmatrix}
\end{align}
where the subscripts $u$ and $p$ refer to (all components of) the
velocity and the pressure, respectively, and $\mat{B}_{pp}$ is zero
unless the \gls{PSPG} formulation is employed.
$\mat{B}_{uu}$ corresponds to the diffusive, advective, and
time-dependent terms of the momentum conservation equation,
$\mat{B}_{up}$ corresponds to the momentum equation pressure term, and
$\mat{B}_{pu}$ corresponds to the continuity equation divergence term.
Several variants of the field-split preconditioner for the \gls{INS}
equations are derived from the ``LDU'' block factorization of $\mat{B}$,
which is given by
\begin{align}
  \label{eqn:eqn:nonlinear_algebra_LDU}
  \begin{bmatrix}
    \mat{B}_{uu}   & \mat{B}_{up}
    \\[6pt]
    \mat{B}_{pu}   & \mat{B}_{pp}
  \end{bmatrix}
  =
  \begin{bmatrix}
    \mat{I}                     & \mat{0}
    \\[6pt]
    \mat{B}_{pu} \mat{B}_{uu}^{-1} & \mat{I}
  \end{bmatrix}
  \begin{bmatrix}
    \mat{B}_{uu}  & \mat{0}
    \\[6pt]
    \mat{0} & \mat{S}
  \end{bmatrix}
  \begin{bmatrix}
    \mat{I}  & \mat{B}_{uu}^{-1} \mat{B}_{up}
    \\[6pt]
    \mat{0} & \mat{I}
  \end{bmatrix}
\end{align}
where $\mat{S} \equiv \mat{B}_{pp} - \mat{B}_{pu} \mat{B}_{uu}^{-1}
\mat{B}_{up}$ is the Schur complement. The first matrix on the
right-hand side of~\eqref{eqn:eqn:nonlinear_algebra_LDU} is block-lower-triangular,
the second matrix is block-diagonal, and the third is block-upper-triangular.

The primary challenge in applying the field-split preconditioner
involves approximating $\mat{S}$, which is generally dense due to its
dependence on $\mat{B}_{uu}^{-1}$, a matrix which is usually not
explicitly formed.  A comprehensive description of advanced techniques
for approximating $\mat{S}$ is beyond the scope of this work;
interested readers should refer to~\cite{Elman_2008,Quarteroni_2000,Wesseling_2009}
for more information. The \ns module currently supports several choices for approximating
$\mat{S}$ via the PETSc command line option
\begin{quote}
  \centering
  \texttt{-pc{\us}fieldsplit{\us}schur{\us}precondition} \{\texttt{a11, selfp, user, full}\}
\end{quote}
The different options are:
\begin{itemize}
\item {\texttt{a11}}: $\mat{S} = \mat{B}_{pp}$. This approach requires \gls{PSPG} stabilization, otherwise $\mat{B}_{pp}$ is zero.
\item {\texttt{selfp}}: $\mat{S} = \mat{B}_{pp} - \mat{B}_{pu} \text{diag}(\mat{B}_{uu})^{-1} \mat{B}_{up}$. Only considers the diagonal when approximating $\mat{B}_{uu}^{-1}$.
\item {\texttt{user}}: Use the user-provided matrix as $\mat{S}$. Supports the construction of application-specific preconditioners.
\item {\texttt{full}}: $\mat{S}$ is formed exactly. This approach is typically only used for debugging due to the expense of computing $\mat{B}_{uu}^{-1}$.
\end{itemize}
By default, all factors on the right-hand side of~\eqref{eqn:eqn:nonlinear_algebra_LDU}
are used in constructing the preconditioner, but we can instead compute a simplified representation,
which ignores some of the factors, by using the PETSc command line option
\begin{quote}
  \centering
  \texttt{-pc{\us}fieldsplit{\us}schur{\us}fact{\us}type} \{\texttt{diag, lower, upper, full}\}
\end{quote}
An advantage of the field-split
preconditioner approach is that one of the fully-coupled
preconditioners mentioned previously can be applied to compute the
action of $\mat{B}^{-1}_{uu}$ or $\mat{S}$ directly, and is likely to
perform better than if it were applied directly to $\mat{B}$.
An \gls{INS} example which demonstrates the use of field-split
preconditioners for a simple problem can be found in the
\texttt{pressure{\us}channel} test directory of the \ns module in the
\texttt{open{\us}bc{\us}pressure{\us}BC{\us}fieldSplit.i} input file.

\section{Software implementation\label{sec:software}}
The \gls{MOOSE} framework is a general purpose toolkit in the vein of other
popular, customizable \cpp \gls{FEM} libraries such as
deal.II~\cite{Bangerth_2007}, DUNE~\cite{Bastian_2008},
GRINS~\cite{Bauman_2016}, \freefem\cite{Hecht_2012},
OpenFOAM~\cite{Jasak_2007}, and the computational back end of FEniCS,
DOLFIN~\cite{Logg_2012}. Generally speaking, these libraries are
designed to enable practitioners who are familiar with the engineering
and applied mathematics aspects of their application areas (that is,
who already have a mathematical model and associated variational
statement) to translate those methods from ``pen and paper''
descriptions into portable, extensible, well supported, high
performance software.  The libraries differ in their approaches to
this translation; various techniques include user-developed/library-assisted hand
written \cpp, automatic low-level code generation based on domain
specific languages, and reconfigurable, extensible pre-developed
software modules for specific applications.

The \gls{MOOSE} framework employs a combination of hand written code
and pre-developed modules in its ``method translation'' approach. Rather than developing
custom software configuration management code and low-level meshing,
finite element, and numerical linear algebra library interfaces,
\gls{MOOSE} relies on the libMesh~\cite{Kirk_2006} finite element
library's implementation of these features. The \gls{MOOSE} application
programming interface focuses on several high level
``Systems'' which map to standard concepts in finite element
programming, including: weak form representation (\texttt{Kernels}),
auxiliary variable calculation (\texttt{AuxKernels}), boundary/initial
conditions, program flow control (\texttt{Executioners}), material
properties, data transfers, and postprocessors. In the rest of this
section, we describe how these systems are leveraged within the \ns
module using examples based on the numerical method for the \gls{INS}
equations described in \S\ref{sec:numerical_method}.

\subsection{Sample Kernel: \texttt{INSMomentumTimeDerivative}\label{sec:ins_time_derivative}}
To connect the \gls{INS} finite element formulation of \S\ref{sec:numerical_method} to the software in
the \ns module, we consider, in detail, the time-dependent term
from the semi-discrete momentum component equation~\eqref{eqn:semi_discrete_momentum_k}:
\begin{align}
  \label{eqn:time_derivative_alone}
  F_i \equiv \int_{\Omega'} \rho \frac{\partial \hcomp{u}{k}}{\partial t} \varphi_i \; \text{d}\Omega'
\end{align}
In order to compute $F_i$, several key ingredients are required:
\begin{enumerate}
  \setlength\itemsep{.1em}
  \item A loop over the finite elements.
  \item A loop over the quadrature points on each element.
  \item Access to the constant value $\rho$, which can be set by the user.
  \item A characterization of the $\frac{\partial \hcomp{u}{k}}{\partial t}$ term consistent with the chosen time integration scheme.
\end{enumerate}
Without much loss of generality, in this discussion we will assume
the approximate time derivative can be written as a linear function of the finite
element solution, i.e.
\begin{align}
  \label{eqn:time_derivative_linear_function}
  \frac{\partial \hcomp{u}{k}}{\partial t} = \sigma_1 u_k^h + \sigma_2
\end{align}
where $\sigma_1$ and $\sigma_2$ are coefficients that may depend on
the current timestep, the solution from the previous timestep, and other factors.
Then, the term~\eqref{eqn:time_derivative_alone} implies the nested summation
\begin{align}
  \label{eqn:time_derivative_alone_sums}
  F_i = \sum_{e} \sum_q |J_e(\vec{x}_q)| w_q \Big( \rho \left(\sigma_1 u_k^h(\vec{x}_q) + \sigma_2 \right) \varphi^e_i(\vec{x}_q)   \Big)
\end{align}
over the elements $e$ and quadrature points $q$ on each element, where
$w_q$ is a quadrature weight,
$\varphi^e_i$ represents the restriction of global basis function
$\varphi_i$ to element $\Omega_e$, $|J_e|$ is the determinant of the
Jacobian of the mapping between the physical element $\Omega_e$ and
the reference element $\hat{\Omega}_e$, and each of the terms is
evaluated at the quadrature point $\vec{x}_q$. The corresponding
Jacobian contribution for this term is given by:
\begin{align}
  \label{eqn:time_derivative_alone_jac}
  J_{ij} = \sum_{e} \sum_q |J_e(\vec{x}_q)| w_q \, \rho\, \sigma_1 \varphi^e_j(\vec{x}_q) \, \varphi^e_i(\vec{x}_q)
\end{align}

\begin{figure}[htpb]
  \centering
  \includegraphics[width=.9\linewidth]{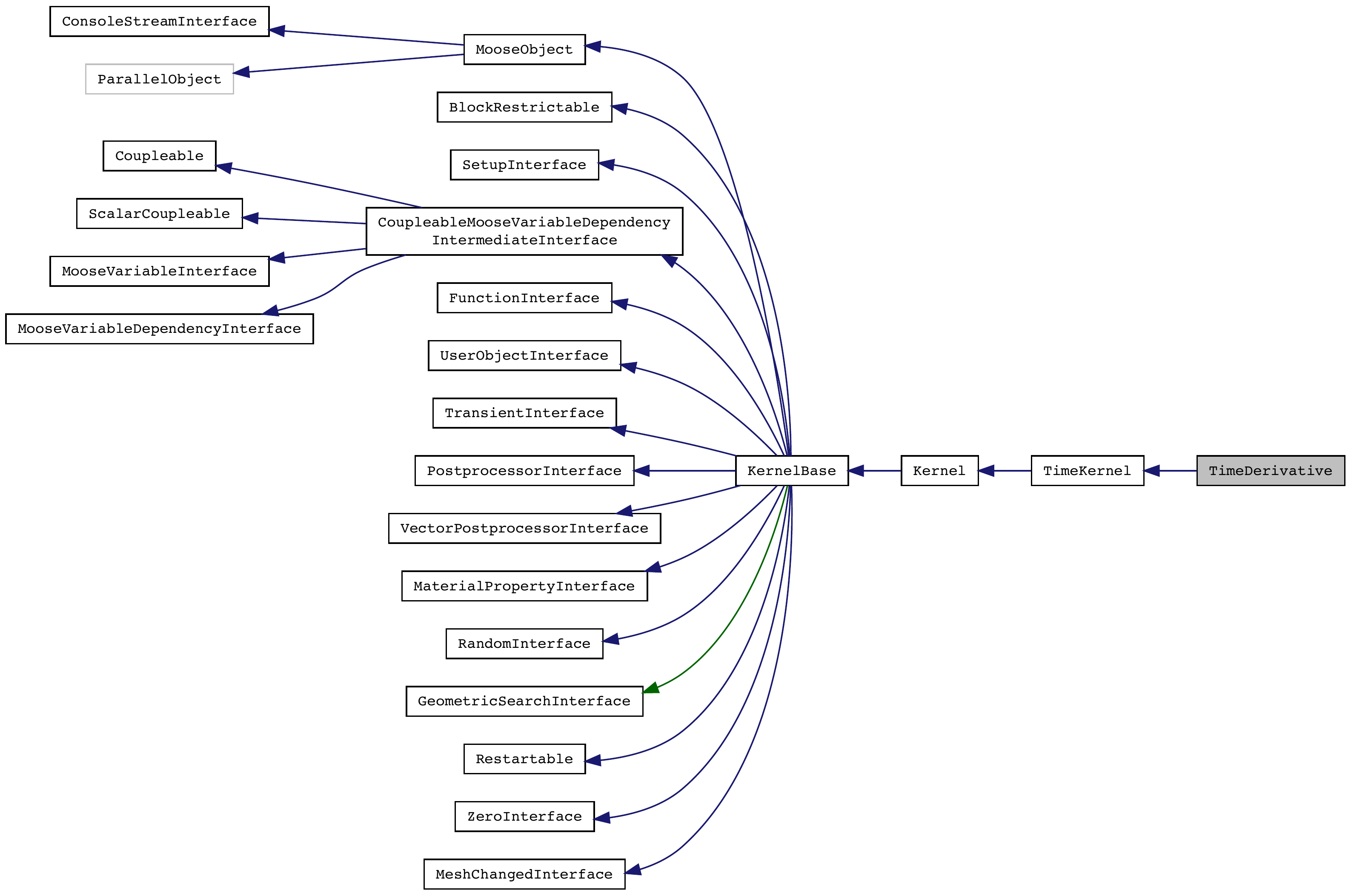}
  \caption{Inheritance diagram for the \texttt{TimeDerivative}
    \texttt{Kernel}.  Arrows point from child (derived) classes to parent
    (base) classes. Multiply-inherited classes with the \texttt{Interface}
    suffix are used to \emph{append} functionality to derived classes,
    while singly-inherited classes are used to \emph{override} base
    class functionality, two standard object-oriented programming
    techniques.\label{fig:time_derivative_hierarchy}}
\end{figure}

In \gls{MOOSE}, the element loop, quadrature loop, and
multiplication by the element Jacobian and quadrature weight
in~\eqref{eqn:time_derivative_alone_sums} and~\eqref{eqn:time_derivative_alone_jac} are handled by the
framework, and the user is responsible for writing \cpp code which
multiplies together the remaining terms. That is,
the hand written code is effectively the body of a loop, and the
framework both controls when the loop is called and prepares the values
which are to be used in the calculation. In \gls{MOOSE}
terminology, this loop body is referred to as a \texttt{Kernel}, and
the various types of \texttt{Kernels}, such as the time-dependent
\texttt{Kernel} shown here, are all derived
(in the sense of \cpp inheritance) from a common base class which
resides in the framework. This relationship is depicted graphically
in Fig.~\ref{fig:time_derivative_hierarchy} for \gls{MOOSE}'s \texttt{TimeDerivative}
class, which is the base class for the \ns module's \texttt{INSMomentum\-Time\-Derivative} class.

The diagram in Fig.~\ref{fig:time_derivative_hierarchy} also shows the
extensive use of inheritance, especially multiple inheritance for
interfaces, employed by the \gls{MOOSE} framework.  In this context,
the ``interface'' classes have almost no data, and are used to
provide access to the various \gls{MOOSE} systems (\texttt{Functions},
\texttt{UserObjects}, \texttt{Postprocessors}, pseudorandom numbers, etc.)
within \texttt{Kernels}. Each \texttt{Kernel} also ``is a'' (in the
object-oriented sense) \texttt{MooseObject}, which allows it to be
stored polymorphically in collections of other \texttt{MooseObjects},
and used in a generic manner. Finally, we note that inheriting from
the \texttt{Coupleable} interface allows \texttt{Kernels}
for a given equation to depend on the other variables in the system
of equations, which is essential for multiphysics applications.

\noindent\begin{minipage}{\linewidth}
\begin{lstlisting}[language=C++,caption={\texttt{TimeDerivative} base class public and protected interfaces.},label={lst:time_derivative_h}]
// framework/include/kernels/TimeDerivative.h
class TimeDerivative : public TimeKernel
{
public:
  TimeDerivative(const InputParameters & parameters);
  virtual void computeJacobian() override;
protected:
  virtual Real computeQpResidual() override;
  virtual Real computeQpJacobian() override;
};
\end{lstlisting}
\end{minipage}

The public and protected interfaces for the \texttt{TimeDerivative} base class are shown in Listing~\ref{lst:time_derivative_h}.
There are three virtual interfaces which override the behavior of the parent \texttt{TimeKernel} class:
\texttt{computeJacobian()}, \texttt{computeQpResidual()}, and \texttt{computeQpJacobian()}. The latter two
functions are designed to be called at each quadrature point (as the names suggest) and have default
implementations which are suitable for standard finite element problems. As will be discussed momentarily,
these functions can also be customized and/or reused in specialized applications.

The \texttt{computeJacobian()} function, which is shown in abbreviated form
in Listing~\ref{lst:time_derivative_C}, consists of the loops over shape functions and quadrature
points necessary for assembling a single element's Jacobian contribution into the \texttt{ke} variable.
This listing shows where the \texttt{computeQpJacobian()} function is called in the body of the loop,
as well as where the element Jacobian/quadrature weight (\texttt{{\us}JxW}), and coordinate transformation
value (\texttt{{\us}coord}, used in e.g.\ axisymmetric simulations) are all multiplied together.
Although this assembly loop is not typically overridden in \gls{MOOSE} codes, the possibility
nevertheless exists to do so.

\noindent\begin{minipage}{\linewidth}
\begin{lstlisting}[language=C++,caption={\texttt{Time\-Derivative::compute\-Jacobian()} definition.},label={lst:time_derivative_C}]
// framework/src/kernels/TimeDerivative.C
void
TimeDerivative::computeJacobian()
{
  DenseMatrix<Number> & ke = _assembly.jacobianBlock(_var.number(), _var.number());

  for (_i = 0; _i < _test.size(); _i++)
    for (_j = 0; _j < _phi.size(); _j++)
      for (_qp = 0; _qp < _qrule->n_points(); _qp++)
        ke(_i, _i) += _JxW[_qp] * _coord[_qp] * computeQpJacobian();
}
\end{lstlisting}
\end{minipage}

Moving on to the \ns module, we next consider the implementation of
the \texttt{INSMomentumTimeDerivative} class itself, which is shown in
Listing~\ref{lst:ins_mom_time_deriv_h}. From the listing, we observe
that this class inherits the \texttt{TimeDerivative} class from the
framework, and overrides three interfaces from that class:
\texttt{computeQpResidual()}, \texttt{computeQpJacobian()}, and
\texttt{computeQpOffDiagJacobian()}.  These overrides are the
mechanism by which the framework implementations are customized to
simulate the \gls{INS} equations.
Since~\eqref{eqn:time_derivative_alone_sums}
and~\eqref{eqn:time_derivative_alone_jac} also involve the fluid
density $\rho$, this class manages a constant reference to a
\texttt{MaterialProperty}, \texttt{{\us}rho}, that can be used at each
quadrature point.

\noindent\begin{minipage}{\linewidth}
\begin{lstlisting}[language=C++,caption={\texttt{INSMomentumTimeDerivative} class methods.},label={lst:ins_mom_time_deriv_h}]
// modules/navier_stokes/include/kernels/INSMomentumTimeDerivative.h
class INSMomentumTimeDerivative : public TimeDerivative
{
public:
  INSMomentumTimeDerivative(const InputParameters & parameters);

protected:
  virtual Real computeQpResidual() override;
  virtual Real computeQpJacobian() override;
  virtual Real computeQpOffDiagJacobian(unsigned jvar) override;

  // Parameters
  const MaterialProperty<Real> & _rho;
};
\end{lstlisting}
\end{minipage}

The implementation of the specialized \texttt{computeQpResidual()} and \texttt{computeQpJacobian()}
functions of the \texttt{INSMomentumTimeDerivative} class are shown in Listing~\ref{lst:ins_mom_time_deriv_C}.
The only difference between the base class and specialized implementations is multiplication by
the material property $\rho$, so the derived class implementations are able to reuse code by calling the
base class methods directly. Although this example is particularly simple, the basic ideas extend to the
other \texttt{Kernels}, boundary conditions, \texttt{Postprocessors}, etc.\ used in the \ns module.
In \S\ref{sec:ns_kernel_design}, we go into a bit more detail about the overall design of the
\gls{INS} \texttt{Kernels}, while continuing to highlight the manner in which \cpp language features are used
to encourage code reuse and minimize code duplication.

\noindent\begin{minipage}{\linewidth}
\begin{lstlisting}[language=C++,caption={\texttt{INSMomentumTimeDerivative} implementations showing relation to the base \texttt{TimeDerivative} class methods.},label={lst:ins_mom_time_deriv_C}]
// modules/navier_stokes/src/kernels/INSMomentumTimeDerivative.C
Real
INSMomentumTimeDerivative::computeQpResidual()
{
  return _rho[_qp] * TimeDerivative::computeQpResidual();
}

Real
INSMomentumTimeDerivative::computeQpJacobian()
{
  return _rho[_qp] * TimeDerivative::computeQpJacobian();
}
\end{lstlisting}
\end{minipage}

\subsection{Navier--Stokes module \texttt{Kernel} design\label{sec:ns_kernel_design}}
All \gls{INS} kernels, regardless of whether they contribute to the
mass or momentum equation residuals, inherit from
the \texttt{INSBase} class whose interface is partially shown in
Listing~\ref{lst:ins_base_h}. This approach reduces code duplication
since both the momentum and mass equation stabilization terms
(\gls{SUPG} and \gls{PSPG} terms, respectively) require access to the
strong form of the momentum equation residual.

Residual contributions
which have been integrated by parts are labeled with the
``\texttt{weak}'' descriptor in the \texttt{INSBase} class, while the
non-integrated-by-parts terms have the ``\texttt{strong}''
descriptor. Methods with a ``\texttt{d}'' prefix are used for
computing Jacobian contributions. Several of the \texttt{INSBase}
functions return \texttt{RealVectorValue} objects, which are
vectors of length $n_s$ that have mathematical operations
(inner products, norms, etc.) defined on them.

\begin{lstlisting}[language=C++,caption={Partial listing of \texttt{INSBase} class methods.},label={lst:ins_base_h}]
// modules/navier_stokes/include/kernels/INSBase.h
class INSBase : public Kernel
{
public:
  INSBase(const InputParameters & parameters);

protected:
  virtual RealVectorValue convectiveTerm();
  virtual RealVectorValue dConvecDUComp(unsigned comp);

  virtual RealVectorValue strongViscousTermLaplace();
  virtual RealVectorValue strongViscousTermTraction();
  virtual RealVectorValue dStrongViscDUCompLaplace(unsigned comp);
  virtual RealVectorValue dStrongViscDUCompTraction(unsigned comp);

  virtual RealVectorValue weakViscousTermLaplace(unsigned comp);
  virtual RealVectorValue weakViscousTermTraction(unsigned comp);
  virtual RealVectorValue dWeakViscDUCompLaplace();
  virtual RealVectorValue dWeakViscDUCompTraction();

  virtual RealVectorValue strongPressureTerm();
  virtual Real            weakPressureTerm();
  virtual RealVectorValue dStrongPressureDPressure();
  virtual Real            dWeakPressureDPressure();
};
\end{lstlisting}

\begin{figure}[htpb]
  \centering
  \includegraphics[width=.9\linewidth]{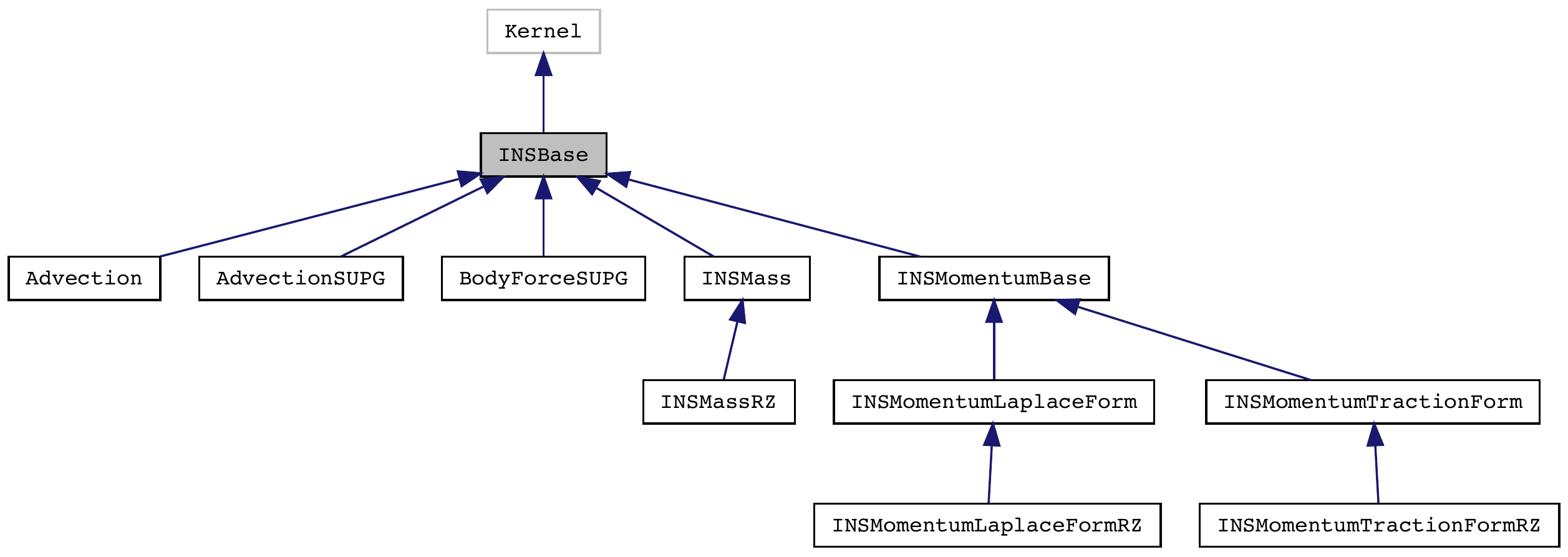}
  \caption{Inheritance diagram for \gls{INS} \texttt{Kernels} in the \ns module.\label{fig:insbase_hierarchy}}
\end{figure}

The inheritance diagram for the \texttt{INSBase} class is shown in
Fig.~\ref{fig:insbase_hierarchy}.  The \texttt{Advection},
\texttt{AdvectionSUPG}, and \texttt{BodyForceSUPG} \texttt{Kernels}
are specific to the scalar advection equation, while the
\texttt{INSMass} and \texttt{INSMomentumBase} classes are used in
simulations of the \gls{INS} equations. \texttt{INSMomentumBase} is
not invoked directly in applications, instead one of the four
subclasses, \texttt{INSMomentum\-Laplace\-Form},
\texttt{INSMomentum\-Traction\-Form},
\texttt{INSMomentum\-Laplace\-FormRZ}, and
\texttt{INSMomentum\-Traction\-FormRZ} must be used, depending on what
coordinate system and viscous term form is employed.
\texttt{INSMomentumBase} is a so-called abstract base class, which is indicated
by the \texttt{=0} syntax on various functions in Listing~\ref{lst:ins_mom_base_h}.

\noindent\begin{minipage}{\linewidth}
\begin{lstlisting}[language=C++,caption={Partial listing of \texttt{INSMomentumBase} class members.},label={lst:ins_mom_base_h}]
// modules/navier_stokes/include/kernels/INSMomentumBase.h
class INSMomentumBase : public INSBase
{
public:
  INSMomentumBase(const InputParameters & parameters);

protected:
  virtual Real computeQpResidual();
  virtual Real computeQpJacobian();
  virtual Real computeQpOffDiagJacobian(unsigned jvar);
  virtual Real computeQpResidualViscousPart() = 0;
  virtual Real computeQpJacobianViscousPart() = 0;
  virtual Real computeQpOffDiagJacobianViscousPart(unsigned jvar) = 0;
};
\end{lstlisting}
\end{minipage}

Recalling~\eqref{eqn:a_definition} and~\eqref{eqn:aL_definition}, the main difference between
the Laplace and traction forms of the momentum equation are in the viscous term.
Therefore, \texttt{compute\-Qp\-Residual\-Viscous\-Part()}, \texttt{compute\-Qp\-Jacobian\-Viscous\-Part()},
and \texttt{compute\-QpOff\-Diag\-Jacobian\-Viscous\-Part()} are left unspecified in the
\texttt{INS\-Momentum\-Base} class. The basic idea behind this design is to let
the derived classes, \texttt{INSMomentum\-Laplace\-Form} and \texttt{INSMomentum\-Traction\-Form},
specialize these functions.  The specialized functions for the Laplace version are shown in Listing~\ref{lst:ins_mom_lap_form_C}.
We note in particular that the off-diagonal contribution (Jacobian contributions
due to variables other than the one the current \texttt{Kernel} is acting on) is zero for this form
of the viscous term.

\noindent\begin{minipage}{\linewidth}
\begin{lstlisting}[language=C++,caption={Pure virtual \texttt{INSMomentumBase} members
    overridden by the \texttt{INSMomentum\-LaplaceForm} subclass.},label={lst:ins_mom_lap_form_C}]
// modules/navier_stokes/src/kernels/INSMomentumLaplaceForm.C
Real
INSMomentumLaplaceForm::computeQpResidualViscousPart()
{
  return _mu[_qp] * (_grad_u[_qp] * _grad_test[_i][_qp]);
}

Real
INSMomentumLaplaceForm::computeQpJacobianViscousPart()
{
  return _mu[_qp] * (_grad_phi[_j][_qp] * _grad_test[_i][_qp]);
}

Real
INSMomentumLaplaceForm::computeQpOffDiagJacobianViscousPart(unsigned /*jvar*/)
{
  return 0.;
}
\end{lstlisting}
\end{minipage}

Finally, as will be discussed in \S\ref{sec:axi}, the
\texttt{INSMassRZ}, \texttt{INS\-Momentum\-Laplace\-FormRZ}, and
\texttt{INS\-Momentum\-Traction\-FormRZ} \texttt{Kernels} implement
additional terms which appear in the axisymmetric form of the
governing equations. In this case, the object of the inheritance
structure shown in Fig.~\ref{fig:insbase_hierarchy} is not for derived
classes to override their base class functionality, but instead to append
to it. The implementation is therefore analogous to the way in which
the \texttt{INS\-Momentum\-Time\-Derivative} class explicitly calls
methods from the base \texttt{TimeDerivative} class.

\subsection{Sample input file\label{sec:sample_input}}
The primary way in which users construct and customize MOOSE-based
simulations is by creating and modifying text-based input files. These
files can be developed using either a text editor (preferably one which
can be extended to include application-specific syntax highlighting
and keyword suggestion/auto completion such as
Atom\footnote{\url{https://atom.io}}) or via the MOOSE GUI, which is
known as Peacock. In this section, we describe in some detail the
structure of an input file that will be subsequently used in
\S\ref{sec:sphere} to compute the three-dimensional flow field
over a sphere.

The first part of the input file in question is shown in
Listing~\ref{lst:global_params}.  The file begins by defining the
\texttt{mu} and \texttt{rho} variables which correspond to the
viscosity and density of the fluid, respectively, and can be used to
fix the Reynolds number for the problem. These variables, which do not
appear in a named section (regions demarcated by square
brackets) are employed in a preprocessing step by the input
file parser based on a syntax known as ``dollar bracket expressions''
or DBEs. As we will see subsequently in
Listing~\ref{lst:functions_materials}, the expressions
\texttt{\$\{mu\}} and \texttt{\$\{rho\}} are replaced by the specified
numbers everywhere they appear in the input file, before any other
parsing takes place.

\noindent\begin{minipage}{\linewidth}
\begin{lstlisting}[language=bash,caption={\texttt{[GlobalParams]} section for the ``flow over a sphere'' problem
solved in \S\ref{sec:sphere}. This initial part of the input file is used to
customize the weak form that will be solved and set constants common to all other
input file sections.},label={lst:global_params}]
# Automatic substitution variables
mu=4e-3
rho=1

[GlobalParams]
  # Variable coupling and naming
  u = vel_x
  v = vel_y
  w = vel_z
  p = p

  # Stabilization parameters
  supg = true
  pspg = true
  alpha = 1e0

  # Problem coefficients
  gravity = '0 0 0'

  # Weak form customization
  convective_term = true
  integrate_p_by_parts = true
  transient_term = true
  laplace = true
[]
\end{lstlisting}
\end{minipage}

The substitution variables are followed by the \texttt{[GlobalParams]}
section, which contains key/value pairs that will be used (if
applicable) in all of the other input file sections. In this example,
we use the \texttt{[GlobalParams]} section to turn on the SUPG and
PSPG stabilization contributions, declare variable coupling which is
common to all \texttt{Kernels}, and customize the weak formulation by
e.g.\ turning on/off the convective and transient terms, and toggling
whether the pressure gradient term in the momentum equations is
integrated by parts. Specifying the variable name mapping,
i.e.\ \texttt{u = vel{\us}x}, gives the user control over the way
variables are named in output files, and the flexibility to couple
together different sets of equations while avoiding name collisions and
the need to recompile the code.

%

\noindent\begin{minipage}{\linewidth}
\begin{lstlisting}[language=bash,caption={Kernels for the ``flow over a sphere'' problem described in \S\ref{sec:sphere}.},label={lst:kernels}]
[Kernels]
  # continuity equation
  [./mass]
    type = INSMass
    variable = p
  [../]

  # x1-momentum equation
  [./x_time]
    type = INSMomentumTimeDerivative
    variable = vel_x
  [../]
  [./x_momentum_space]
    type = INSMomentumLaplaceForm
    variable = vel_x
    component = 0
  [../]

  # additional kernels...
[]
\end{lstlisting}
\end{minipage}

The next input file section, \texttt{[Kernels]}, is shown in
Listing~\ref{lst:kernels}. This section controls both which
\texttt{Kernels} are built for the simulation, and which
equation each \texttt{Kernel} corresponds to. In the \ns module, the
temporal and spatial parts of the momentum equations are split into
separate \texttt{Kernels} to allow the user to flexibly define either
transient or steady state simulations. In each case, the \texttt{type}
key corresponds to a \cpp class name that has been registered with MOOSE
at compile time, and the \texttt{variable} key
designates the equation a particular \texttt{Kernel} applies to.

The \texttt{INSMomentumLaplaceForm} \texttt{Kernel} also has a
\texttt{component} key which is used to specify whether the
\texttt{Kernel} applies to the $x_1$, $x_2$, or $x_3$-momentum
equation. In this input file snippet, only the \texttt{Kernel} for
$x_1$-momentum equation is shown; the actual input file has similar
sections corresponding to the other spatial directions. This design
allows multiple copies of the same underlying \texttt{Kernel} code to
be used in a dimension-agnostic way within the simulation.  Finally,
we note that the size of the input file can be reduced
through the use of custom MOOSE \texttt{Actions} which are capable
of adding \texttt{Variables}, \texttt{Kernels}, \texttt{BCs}, etc.\
to the simulation programmatically, although this is an advanced
approach which is not discussed in detail here.

Listing~\ref{lst:functions_materials} shows the \texttt{[Functions]},
\texttt{[BCs]}, and \texttt{[Materials]} blocks from the same input
file. In the \texttt{[Functions]} block, a \texttt{ParsedFunction}
named ``\texttt{inlet{\us}func}'' is declared. In MOOSE, \texttt{ParsedFunctions}
are constructed from strings of standard mathematical expressions
such as \texttt{sqrt}, \texttt{*}, \texttt{/}, \texttt{\^}, etc.
These strings are parsed, optimized, and compiled at runtime,
and can be called at required spatial locations and times during
the simulation (the characters \texttt{x}, \texttt{y}, \texttt{z},
and \texttt{t} in such strings are automatically treated as
spatiotemporal coordinates in MOOSE). In the spherical flow application,
\texttt{inlet{\us}func} is used in the \texttt{vel{\us}z{\us}inlet}
block of the \texttt{[BCs]} section to specify a parabolic inflow profile
in the $x_3$-direction. The boundary condition which imposes this
inlet profile is of type \texttt{FunctionDirichletBC}, and it acts
on the \texttt{vel{\us}z} variable, which corresponds to the $x_3$
component of the velocity.

\noindent\begin{minipage}{\linewidth}
  \begin{lstlisting}[language=bash,caption={\texttt{[Functions]}, \texttt{[BCs]}, and \texttt{[Materials]} blocks
      for the ``flow over a sphere'' problem described in \S\ref{sec:sphere}.},label={lst:functions_materials}]
[Functions]
  [./inlet_func]
    type = ParsedFunction
    value = 'sqrt((x-2)^2 * (x+2)^2 * (y-2)^2 * (y+2)^2) / 16'
  [../]
[]

[BCs]
  [./vel_z_inlet]
    type = FunctionDirichletBC
    function = inlet_func
    variable = vel_z
    boundary = inlet
  [../]

  # additional bcs...
[]

[Materials]
  [./const]
    type = GenericConstantMaterial
    prop_names = 'rho mu'
    prop_values = '${rho}  ${mu}'
  [../]
[]
\end{lstlisting}
\end{minipage}

The last section of Listing~\ref{lst:functions_materials} shows the
\texttt{[Materials]} block for the spherical flow problem. In MOOSE,
material properties can be understood generically as ``quadrature
point quantities,'' that is, values that are computed independently at
each quadrature point and used in the finite element assembly
routines. In general they can depend on the current solution and other
auxiliary variables computed during the simulation, but in this
application the material properties are simply constant. Therefore,
they are implemented using the correspondingly simple
\texttt{Generic\-Constant\-Material} class which is built in to
MOOSE. \texttt{Generic\-Constant\-Material} requires the user to
provide two input parameter lists named \texttt{prop{\us}names} and
\texttt{prop{\us}values}, which contain, respectively, the names of
the material properties being defined and their numerical values.  In
this case, the numerical values are provided by the dollar bracket
expressions \texttt{\$\{rho\}} and \texttt{\$\{mu\}}, which the input
file parser replaces with the numerical values specified at the top of
the file.

Listing~\ref{lst:executioner} details the \texttt{[Executioner]}
section of the sphere flow input file. This is the input file section
where the time integration scheme is declared and customized, and the
solver parameters, tolerances, and iteration limits are defined. For
the sphere flow problem, we employ a \texttt{Transient}
\texttt{Executioner} which will perform \texttt{num{\us}steps = 100}
timesteps using an initial timestep of $\Delta t = 0.5$. Since no
specific time integration scheme is specified, the MOOSE default,
first-order implicit Euler time integration, will be used. If the
nonlinear solver fails to converge for a particular timestep, MOOSE
will automatically decrease $\Delta t$ by a factor of $1/2$ and
retry the most recent solve until a timestep smaller than \texttt{dtmin}
is reached, at which point the simulation will exit with an appropriate
error message.

Adaptive timestep selection is controlled by the
\texttt{IterationAdaptiveDT} \texttt{TimeStepper}. For this
\texttt{TimeStepper}, the user specifies an
\texttt{optimal{\us}iterations} number corresponding to their desired
number of Newton steps per timestep. If the current timestep requires
fewer than \texttt{optimal{\us}iterations} to converge, the timestep
is increased by an amount proportional to the
\texttt{growth{\us}factor} parameter, otherwise the timestep is shrunk
by an amount proportional to the \texttt{cutback{\us}factor}. Use of
this \texttt{TimeStepper} helps to efficiently and robustly drive the
simulation to steady state, but does not provide any rigorous
guarantees of local or global temporal error control.

\noindent\begin{minipage}{\linewidth}
  \begin{lstlisting}[language=bash,caption={The \texttt{[Executioner]} blocks
      for the ``flow over a sphere'' problem described in \S\ref{sec:sphere}.},label={lst:executioner}]
[Executioner]
  # TimeIntegrator and TimeStepper customization
  type = Transient
  num_steps = 100
  trans_ss_check = true
  ss_check_tol = 1e-10
  dtmin = 5e-4
  dt = .5
  [./TimeStepper]
    dt = .5
    type = IterationAdaptiveDT
    cutback_factor = 0.4
    growth_factor = 1.2
    optimal_iterations = 5
  [../]

  # Solver tolerances and iteration limits
  nl_rel_tol = 1e-8
  nl_abs_tol = 1e-12
  nl_max_its = 10
  l_tol = 1e-6
  l_max_its = 10
  line_search = 'none'

  # Options passed directly to PETSc
  petsc_options = '-snes_converged_reason -ksp_converged_reason'
  petsc_options_iname = '-pc_type -pc_factor_shift_type -pc_factor_mat_solver_package'
  petsc_options_value = 'lu NONZERO superlu_dist'
[]
\end{lstlisting}
\end{minipage}

Finally, we discuss a few of the linear and nonlinear solver
tolerances that can be controlled by parameters specified in the
\texttt{[Executioner]} block.  The Newton solver's behavior is
primarily governed by the \texttt{nl{\us}rel{\us}tol},
\texttt{nl{\us}abs{\us}tol}, and \texttt{nl{\us}max{\us}its}
parameters, which set the required relative and absolute residual norm
reductions (the solve stops if either one of these tolerances is met)
and the maximum allowed number of nonlinear iterations, respectively,
for the nonlinear solve.  The \texttt{l{\us}tol} and
\texttt{l{\us}max{\us}its} parameters set the corresponding values for
the linear solves which occur at each nonlinear iteration, and the
\texttt{line{\us}search} parameter can be used to specify the
line searching algorithm employed by the nonlinear solver; the
possible options are described in~\S\ref{sec:preconditioning}.

The \texttt{petsc{\us}options},
\texttt{petsc{\us}options{\us}iname}, and
\texttt{petsc{\us}options{\us}value} parameters are used to specify command line
options that are passed directly to PETSc. The first is used to
specify PETSc command line options that don't have a corresponding
value, while the second two are used to specify lists of key/value
pairs that must be passed to PETSc together. This particular input
file gives the PETSc command line arguments required to set up a
specific type of direct solver, but because almost all of PETSc's
behavior is controllable via the command line, these parameters
provide a very flexible and traceable (since most input files are
maintained with version control software) approach for exerting
fine-grained control over the solver.

\section{Verification tests\label{sec:verification_tests}}
In this section, we discuss several verification tests that are
available in the \ns module. These tests were used in the development
of the software, and versions of them are routinely run in support of
continuous integration testing. The tests described here include
verification of the \gls{SUPG} formulation for the scalar advection
equation (\S\ref{sec:scalar_advection}), the \gls{SUPG}/\gls{PSPG}
stabilized formulation of the full \gls{INS} equations
(\S\ref{sec:verif}) via the \gls{MMS}, and finally, in \S\ref{sec:jh},
both the \gls{PSPG}-stabilized and unstabilized formulations of the
\gls{INS} equations based on the classical Jeffery--Hamel exact
solution for two-dimensional flow in a wedge-shaped region.  The
images in this section and \S\ref{sec:results} were created using the
Paraview~\cite{Henderson_2004} visualization tool.

\subsection{SUPG stabilization: Scalar advection equation\label{sec:scalar_advection}}
To verify the \gls{SUPG} implementation in a simplified setting before tackling the
full \gls{INS} equations, a convergence study was first conducted for the scalar advection equation:
\begin{alignat}{2}
  \label{eq:scalar_advection}
  \vec{a} \cdot \nabla u &= f && \; \in \Omega
  \\
  u &= g && \; \in \Gamma
\end{alignat}
where $\vec{a}$ is a constant velocity vector and $f$ is a forcing
function. For the one-dimensional case, $\Omega=[0,1]$, and $\vec{a}$ and $f$ were chosen to be:
\begin{align}
  \label{eq:scalar_one_d_a}
  \vec{a} &= (1,0,0)
  \\
  \label{eq:scalar_one_d_f}
  f &= 1 - x^2
\end{align}
In two-dimensions, $\Omega=[0,1]^2$, and
\begin{align}
  \label{eq:scalar_two_d_a}
  \vec{a} &= (1,1,0)
  \\
  \label{eq:scalar_two_d_f}
  f &= \frac{1}{10} \left(4 \sin \left(\frac{\pi x}{2}\right) + 4\sin(\pi y) + 7\sin\left(\frac{\pi xy}{5}\right) + 5 \right)
\end{align}

\begin{figure}[htpb]
  \centering
  \begin{subfigure}[t]{.49\linewidth}
    \includegraphics[width=\linewidth]{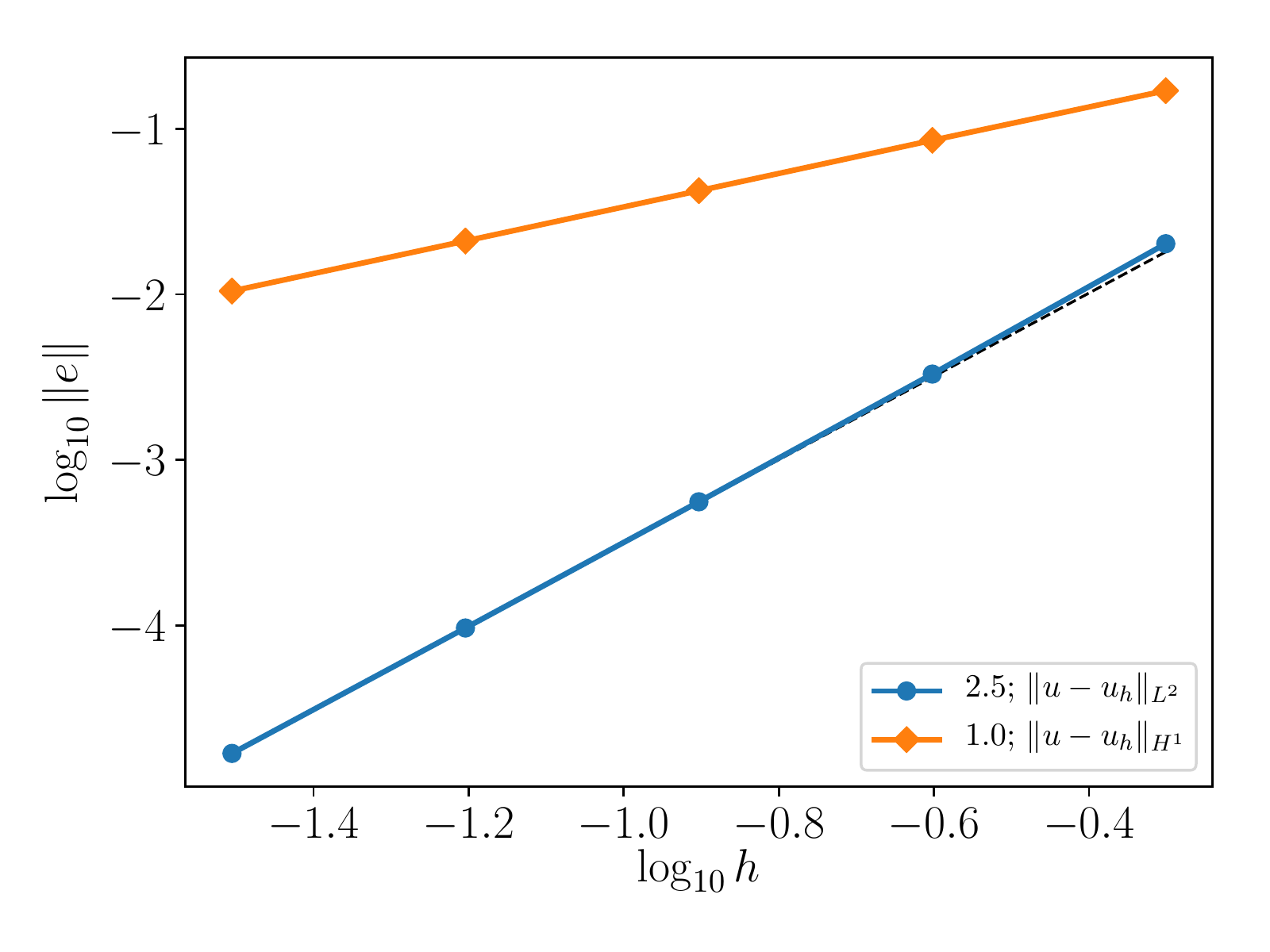}
    \caption{Linear elements, 1D.\label{fig:advection_eqn_P1}}
  \end{subfigure}
  \begin{subfigure}[t]{.49\linewidth}
    \includegraphics[width=\linewidth]{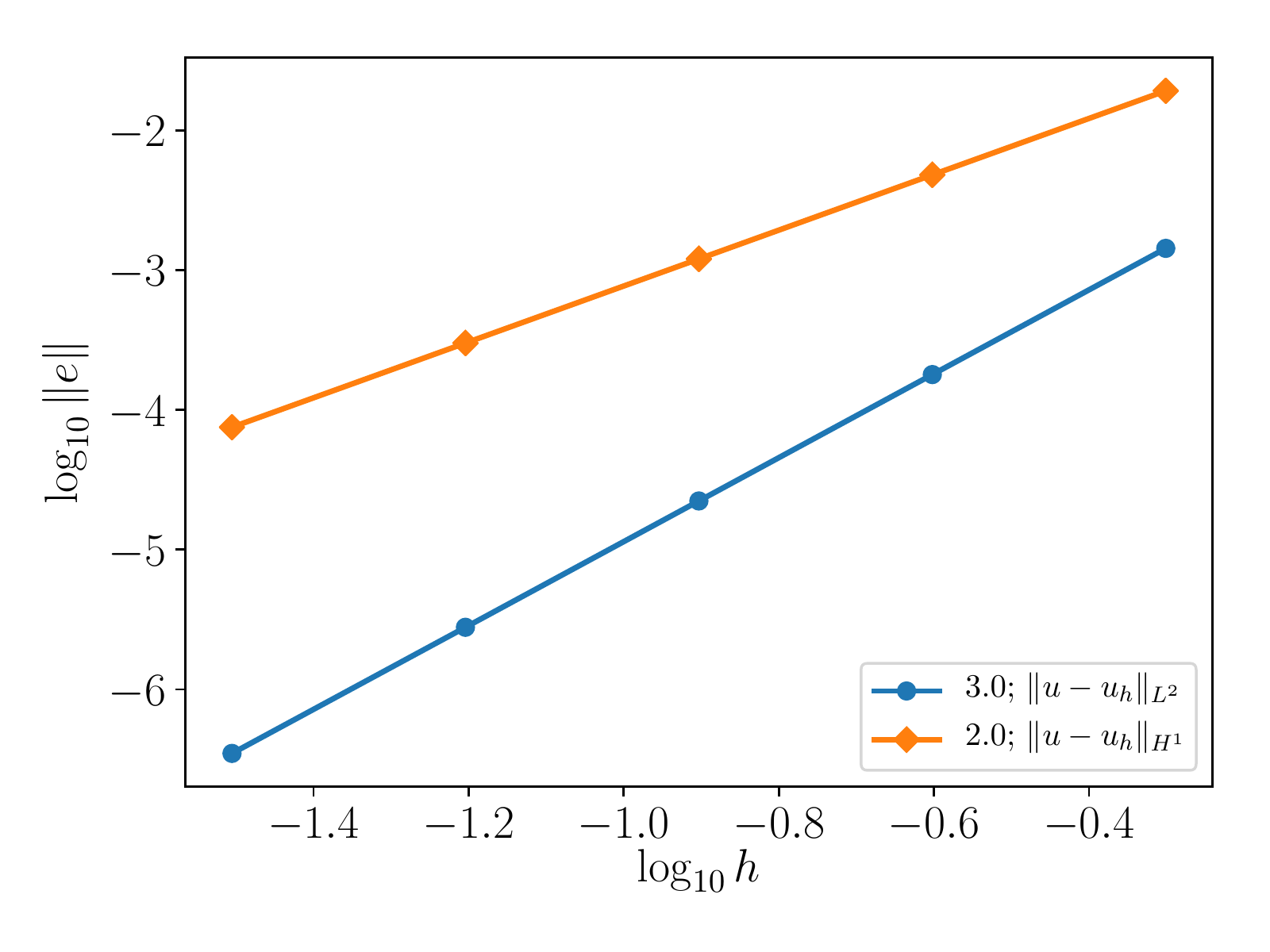}
    \caption{Quadratic elements, 1D.\label{fig:advection_eqn_P2}}
  \end{subfigure}
  \caption{Convergence rates for (\subref{fig:advection_eqn_P1})~$P^1$ and
    (\subref{fig:advection_eqn_P2})~$P^2$ approximate solutions of the pure advection
    equation using the problem parameters given in~\eqref{eq:scalar_one_d_a} and \eqref{eq:scalar_one_d_f}.
    \label{fig:advection_eqn_1D}}
\end{figure}

For linear ($P^1$) elements in 1D, nodally-exact solutions
(superconvergent in the $L^2$-norm) are a well-known characteristic of
the \gls{SUPG} method~\cite{Donea_2003}. We observe a convergence rate
of 2.5 in this norm (see Fig.~\ref{fig:advection_eqn_P1}).  For
quadratic ($P^2$) elements in 1D, convergence in the $L^2$-norm is
generally 3rd-order unless different forms of $\tau$ are used at the
vertex and middle nodes of the elements~\cite{Codina_1992}.
Since this specialized form of $\tau$ is currently not implemented in
the \ns module, we do indeed observe 3rd-order convergence, as shown
in Fig.~\ref{fig:advection_eqn_P2}. The convergence in the $H^1$-norm
is standard for both cases, that is, the superconvergence of the
$L^2$-norm does not carry over to the gradients in the linear element
case.

\begin{table}[htpb]
  \centering
  \caption{Convergence rates, i.e.\ the exponent $p$ in the $O(h^p)$ term, for
    the scalar advection problem with \gls{SUPG} stabilization on
    linear and quadratic elements in 1D ($P^1$, $P^2$) and 2D ($Q^1$, $Q^2$).\label{tab:convergence_scalar}}
  \begin{tabular}[c]{lcc}
    \toprule
    & $\|u - u_h\|_{L^2(\Omega)}$ & $\|u - u_h\|_{H^1(\Omega)}$ \\
    \midrule
    $P^1$ & 2.5 & 1 \\
    $P^2$ & 3 & 2 \\
    \midrule
    $Q^1$ & 2.1 & 1 \\
    $Q^2$ & 3 & 2 \\
    \bottomrule
  \end{tabular}
\end{table}

\begin{figure}[htpb]
  \centering
  \begin{subfigure}[t]{.49\linewidth}
    \includegraphics[width=\linewidth]{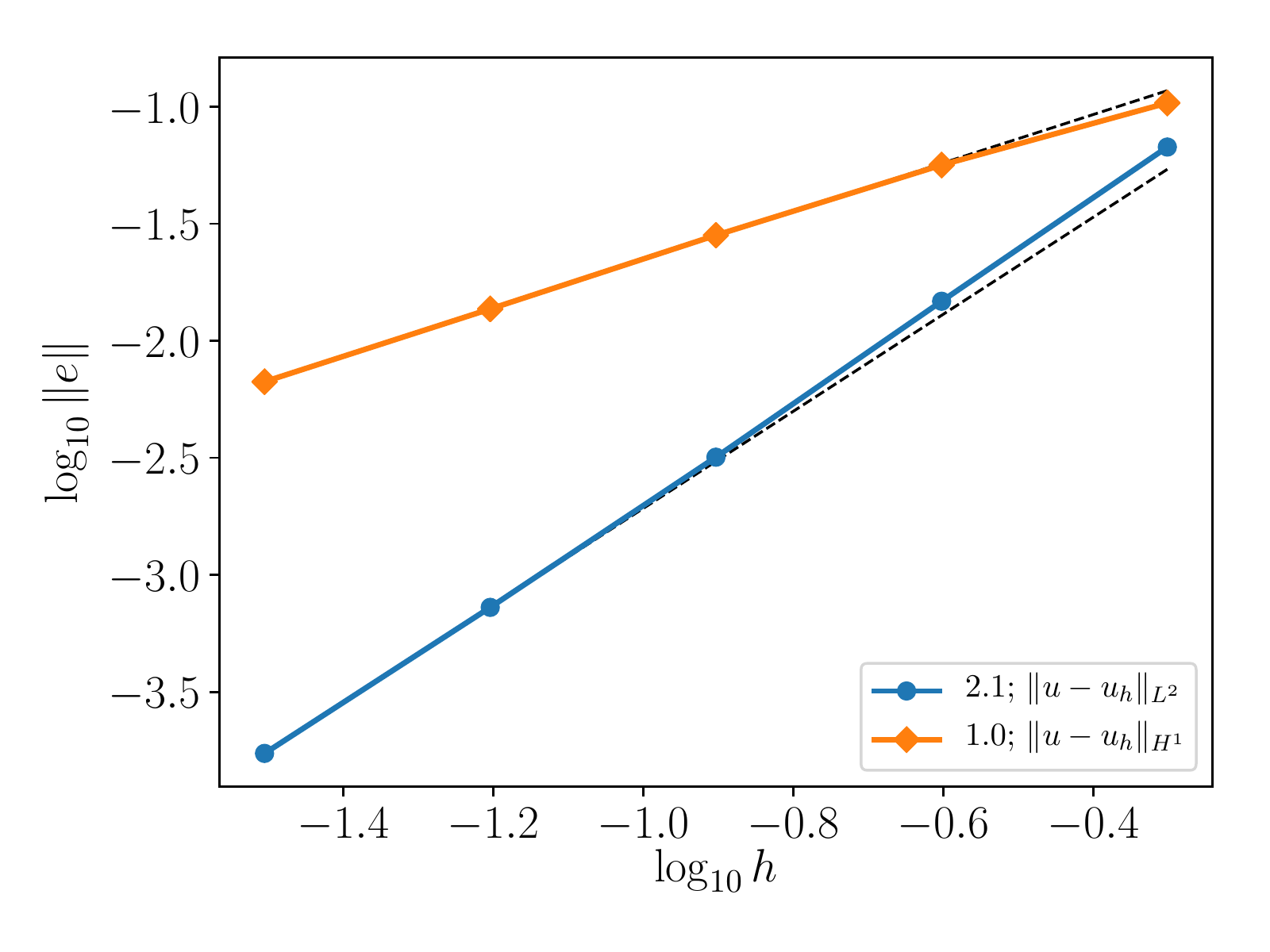}
    \caption{Bilinear elements, 2D.\label{fig:advection_eqn_Q1}}
  \end{subfigure}
  \begin{subfigure}[t]{.49\linewidth}
    \includegraphics[width=\linewidth]{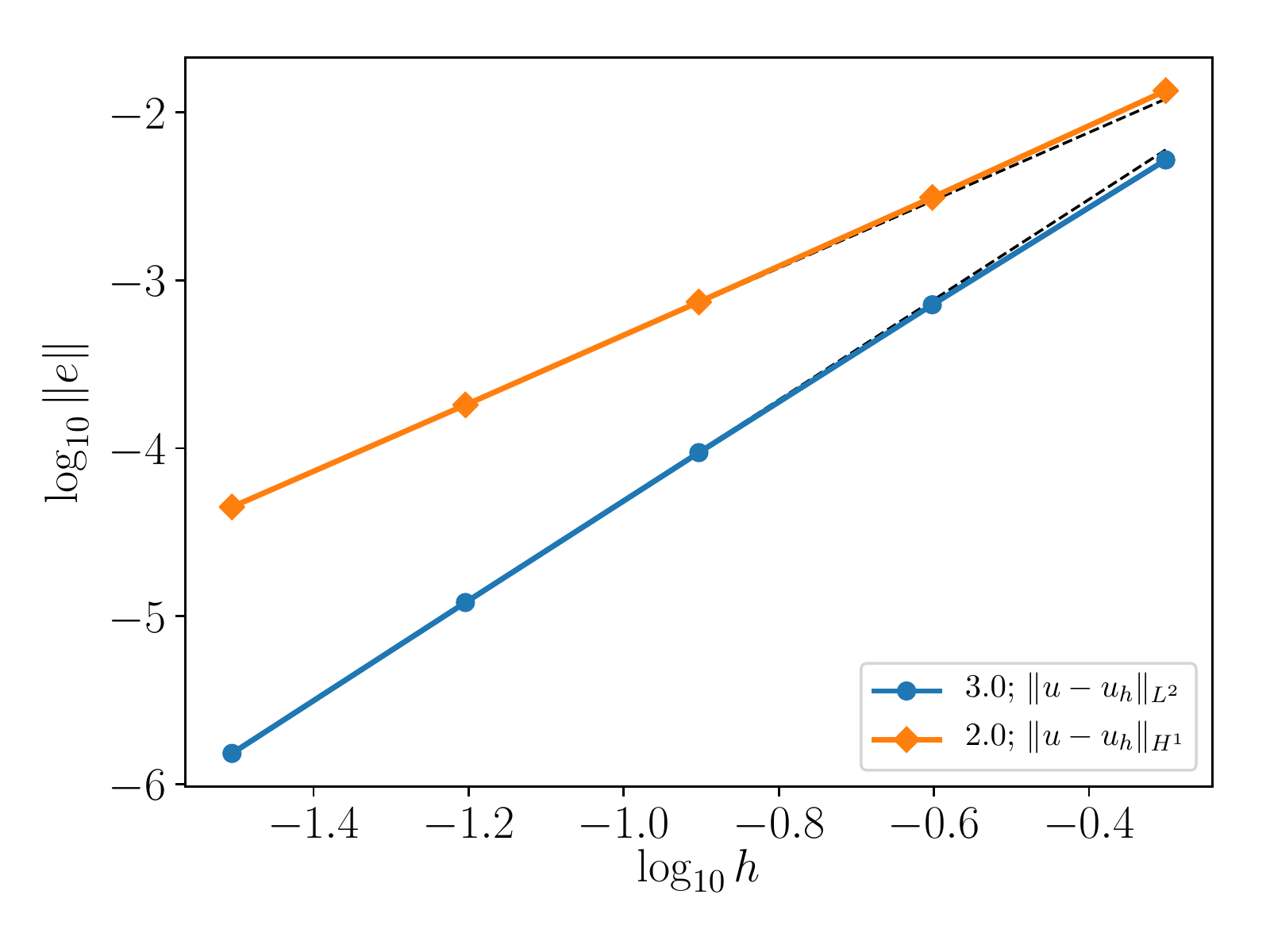}
    \caption{Biquadratic elements, 2D.\label{fig:advection_eqn_Q2}}
  \end{subfigure}
  \caption{Convergence rates for (\subref{fig:advection_eqn_Q1})~$Q^1$ and
    (\subref{fig:advection_eqn_Q2})~$Q^2$ approximate solutions of the pure advection
    equation using the problem parameters given in~\eqref{eq:scalar_two_d_a} and \eqref{eq:scalar_two_d_f}.
    \label{fig:advection_eqn_2D}}
\end{figure}

Nodally-exact/superconvergent solutions are also generally not possible in higher
spatial dimensions, however the optimal rates, $O(h^{k+1})$ in $L^2$ and $O(h^{k})$ in $H^1$, are observed for both
$Q^1$ and $Q^2$ elements, as shown in Fig.~\ref{fig:advection_eqn_2D}.
Table~\ref{tab:convergence_scalar} summarizes the convergence rates
achieved for the scalar advection problem in both 1D and 2D.
Given that each of the tabulated discretizations achieves at least the expected rate of
convergence, we are reasonably confident that the basic approach taken in implementing
this class of residual-based stabilization schemes is correct.

\subsection{SUPG and PSPG stabilization: INS equations\label{sec:verif}}
A 2D convergence verification study of the \gls{SUPG} and \gls{PSPG} methods was conducted via
the \gls{MMS}~\cite{Roache_1998} using the following smooth sinusoidal functions for the
manufactured solutions:
\begin{align}
  u_1 &= \frac{1}{10} \left(4 \sin \left(\frac{\pi x}{2}\right) + 4\sin(\pi y) + 7\sin\left(\frac{\pi xy}{5}\right) + 5 \right)
  \\
  u_2 &= \frac{1}{10} \left( 6\sin \left(\frac{4 \pi x}{5}\right) + 3\sin \left(\frac{3 \pi y}{10}\right) + 2\sin\left(\frac{3\pi x y}{10}\right) + 3 \right)
  \\
  p &= \frac{1}{2} \left( \sin \left(\frac{\pi x}{2}\right) + 2 \sin\left(\frac{3 \pi y}{10}\right) + \sin\left(\frac{\pi x y}{5}\right) + 1\right)
\end{align}
in the unit square domain, $\Omega=[0,1]^2$, using the steady Laplace
form of the governing equations with Dirichlet boundary conditions for
all variables (including $p$) on all sides.

\begin{figure}[htpb]
  \centering
  \begin{subfigure}[t]{.49\linewidth}
    \includegraphics[width=\linewidth]{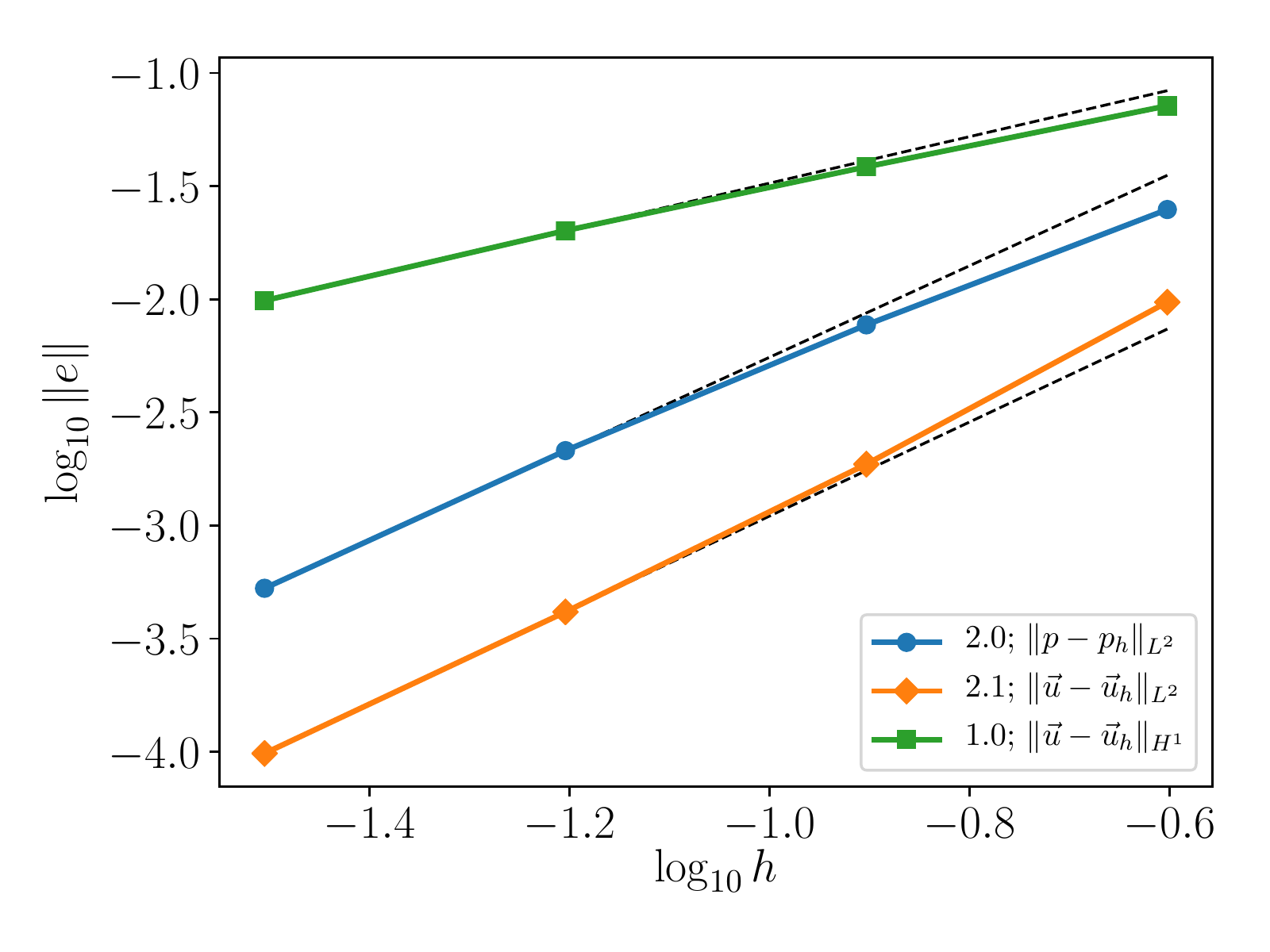}
    \caption{Advection-dominated.\label{fig:supg_pspg_adv}}
  \end{subfigure}
  \begin{subfigure}[t]{.49\linewidth}
    \includegraphics[width=\linewidth]{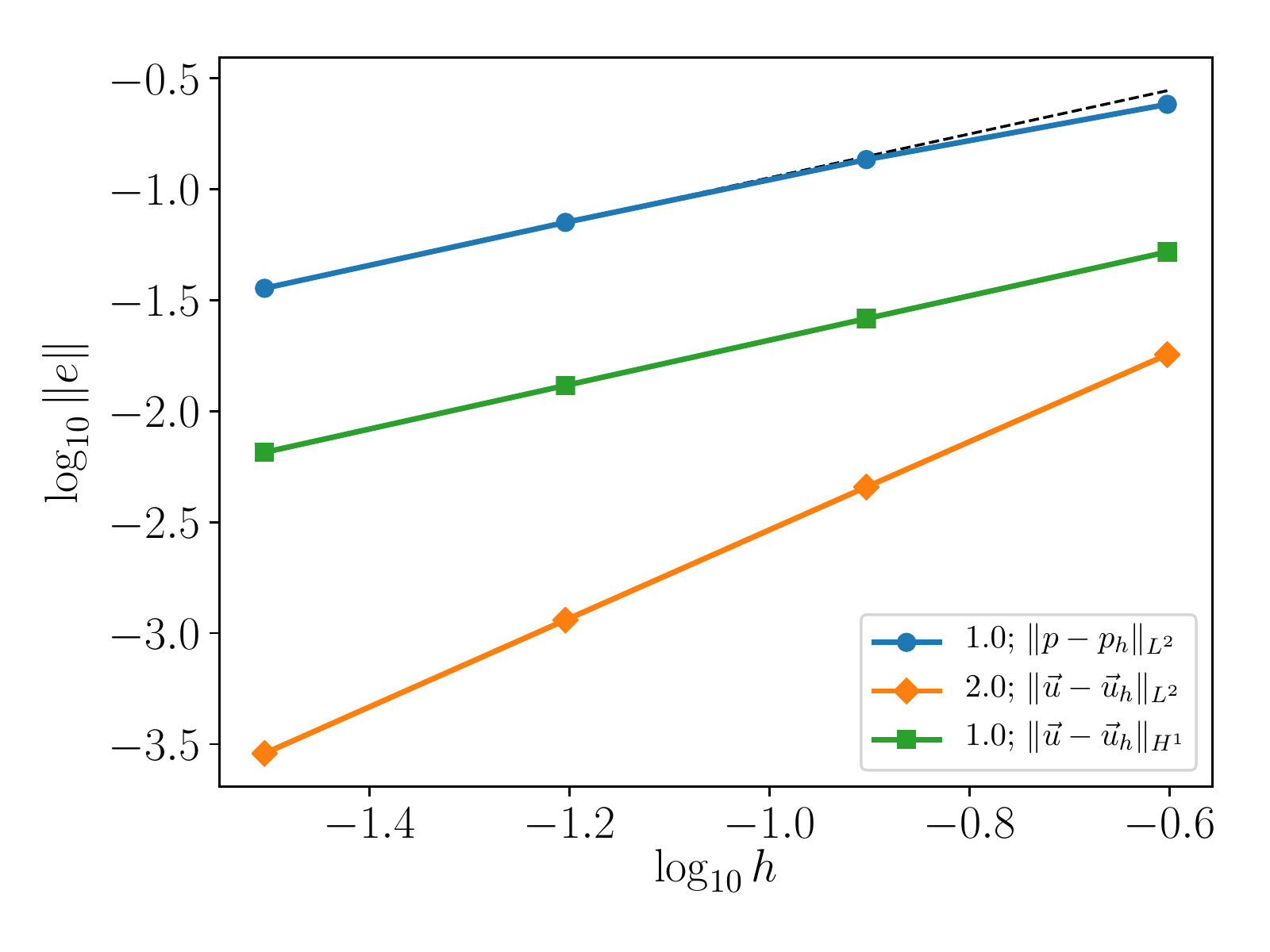}
    \caption{Diffusion-dominated.\label{fig:supg_pspg_diff}}
  \end{subfigure}
  \caption{Convergence rates for the (\subref{fig:supg_pspg_adv})~advection-dominated and
    (\subref{fig:supg_pspg_diff})~diffusion-dominated \gls{INS} equations using \qq{1}{1} elements with both
    \gls{SUPG} and \gls{PSPG} stabilization.\label{fig:supg_pspg_q1q1}}
\end{figure}

Using \qq{1}{1} elements in advection-dominated
flows with combined \gls{SUPG} and \gls{PSPG} stabilization yields optimal
convergence rates, as shown in the first row of Table~\ref{tab:convergence_mms}.
When the flow is diffusion-dominated, the pressure converges at first-order (Table~\ref{tab:convergence_mms}, row 2), which
is consistent with the error analysis conducted
in~\cite{Hughes_1986_PSPG} for \gls{PSPG}-stabilized Stokes
flow, where $\|p-p_h\|_{H^1(\Omega)} = O(h^{k-1})$ was predicted.
Plots of the convergence rates for the advection- and diffusion-dominated
cases can be found in Figs.~\ref{fig:supg_pspg_adv} and~\ref{fig:supg_pspg_diff},
respectively. The advection- and diffusion-dominated cases used
viscosities $\mu=1.5\times 10^{-4}$ and $\mu=15$, respectively. The element Reynolds numbers for the
advection-dominated cases ranged from $7.5 \times 10^3$ down to 469.75 as the mesh was refined,
and from $7.5\times 10^{-2}$ down to $4.6875\times 10^{-3}$ for the diffusion-dominated cases.

\begin{table}[hbt]
  \centering
  \caption{Convergence rates, i.e.\ the exponent $p$ in the $O(h^p)$ term, for
    the \gls{MMS} problem with \gls{SUPG} and \gls{PSPG} stabilization, on
    \qq{1}{1} and \qq{2}{1} elements, for the advection-dominated (A) and
    diffusion-dominated (D) problems. The results with asterisks are
    suboptimal based on the \emph{a priori} error estimates for the linear
    advection-diffusion problem.\label{tab:convergence_mms}}
  \begin{tabular}[c]{lccc}
    \toprule
    & $\|\vec{u} - \vec{u}_h\|_{L^2(\Omega)}$ & $\|\vec{u} - \vec{u}_h\|_{H^1(\Omega)}$ & $\|p - p_h\|_{L^2(\Omega)}$ \\
    \midrule
    \qq{1}{1} (A) & 2 & 1 & 2 \\
    \qq{1}{1} (D) & 2 & 1 & 1 \\
    \midrule
    \qq{2}{1} (A) & 2$^{\ast}$ & 1$^{\ast}$ & 2 \\
    \qq{2}{1} (D) & 3 & 2 & 2 \\
    \bottomrule
  \end{tabular}
\end{table}

\begin{figure}[htpb]
  \centering
  \begin{subfigure}[t]{.49\linewidth}
    \includegraphics[width=\linewidth]{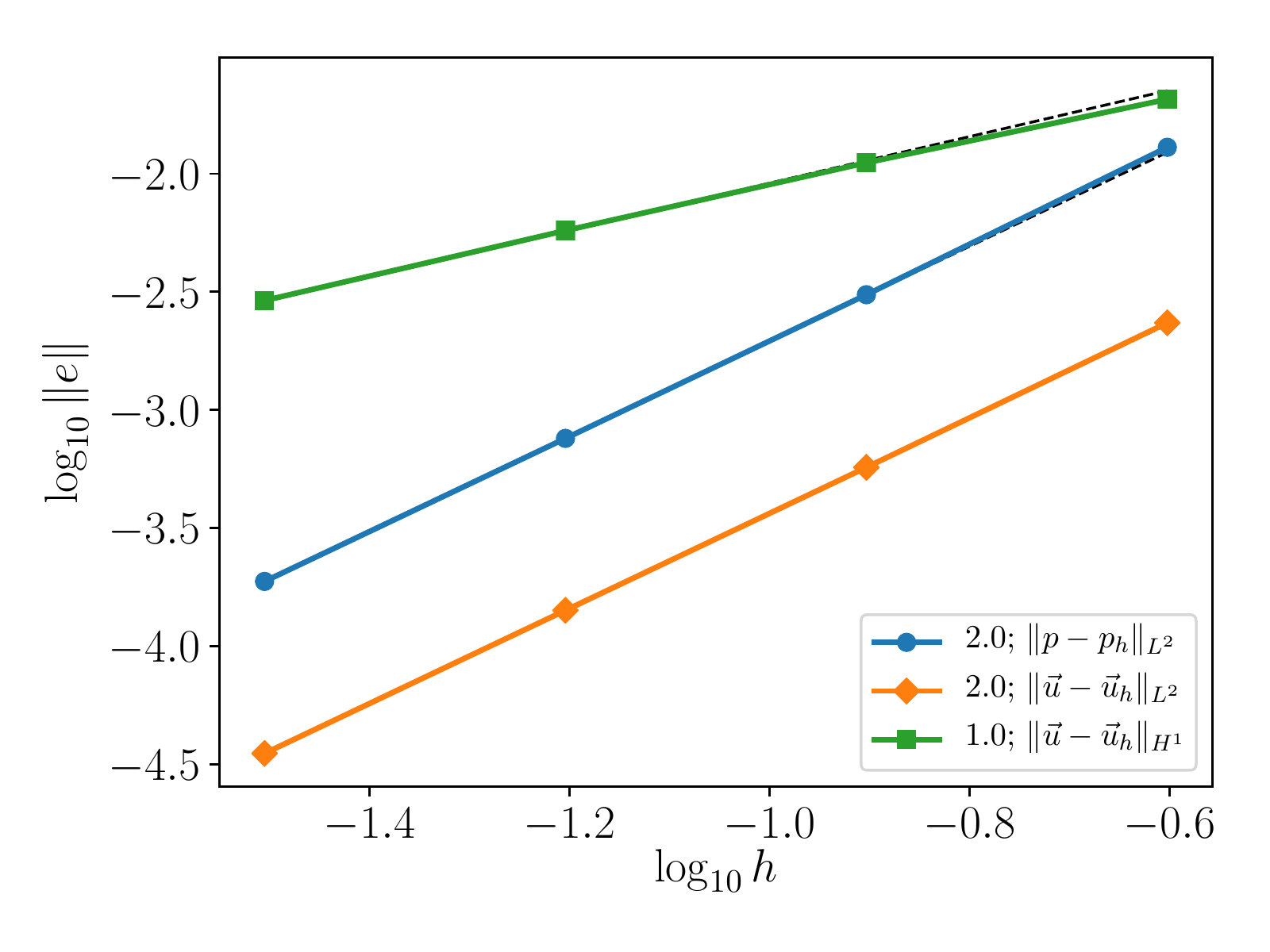}
    \caption{Advection-dominated.\label{fig:supg_adv}}
  \end{subfigure}
  \begin{subfigure}[t]{.49\linewidth}
    \includegraphics[width=\linewidth]{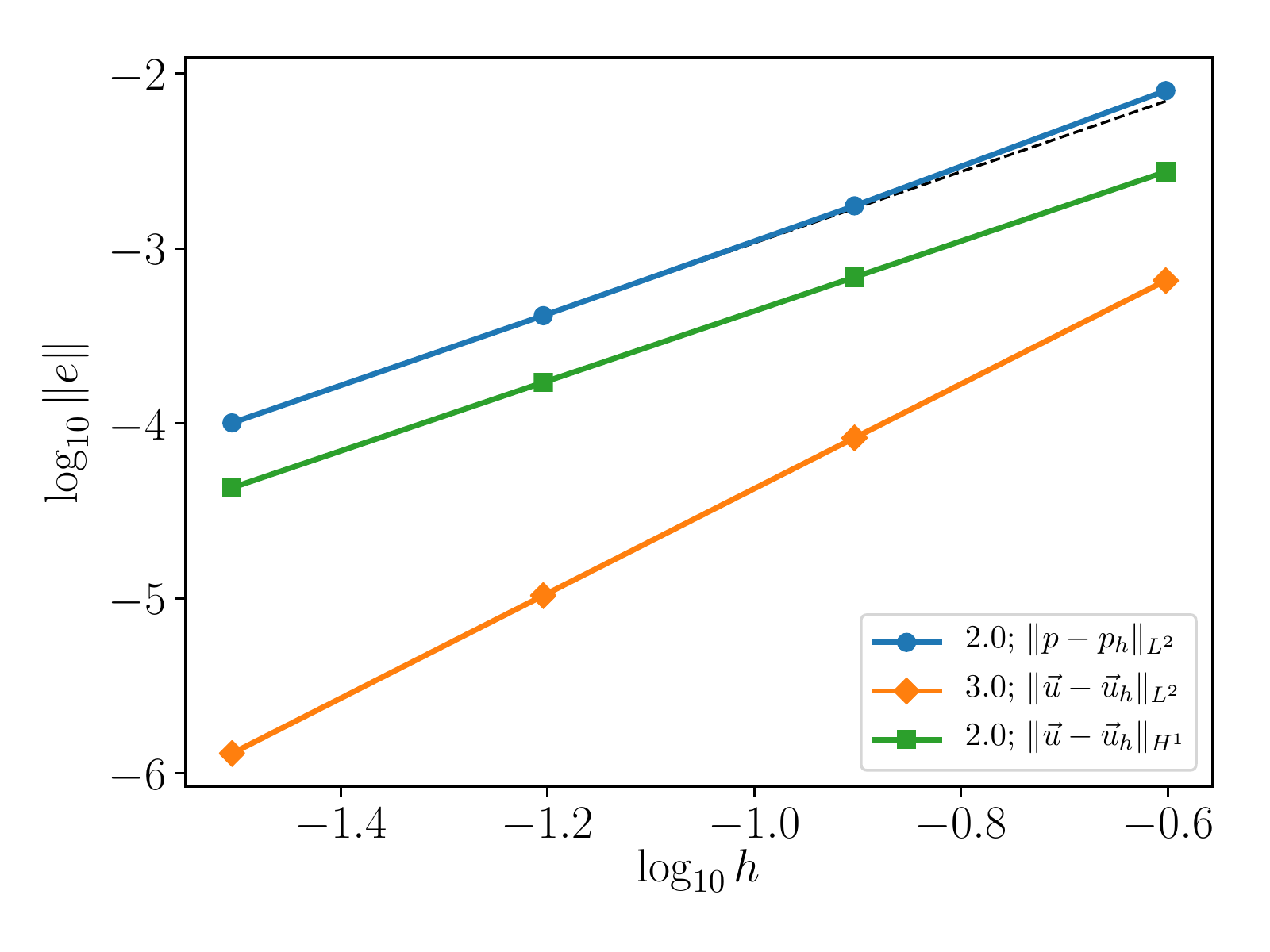}
    \caption{Diffusion-dominated.\label{fig:supg_diff}}
  \end{subfigure}
  \caption{Convergence rates for the (\subref{fig:supg_adv})~advection-dominated and
    (\subref{fig:supg_diff})~diffusion-dominated \gls{INS} equations using \qq{2}{1}
    elements with \gls{SUPG} stabilization. This element pair is \gls{LBB}-stable, and
    therefore \gls{PSPG} stabilization is not required.\label{fig:supg_q2q1}}
\end{figure}

Convergence studies were also conducted using \gls{LBB}-stable
\qq{2}{1} elements with \gls{SUPG} stabilization. For
diffusion-dominated flows, the convergence rates for this element are
optimal (Table~\ref{tab:convergence_mms}, row 4).  For the
advection-dominated flow regime, on the other hand, we observe a drop
in the convergence rates of all norms of $\|\vec{u} - \vec{u}_h\|$ of
a full power of $h$ (Table~\ref{tab:convergence_mms}, row 3).
Johnson~\cite{Johnson_1984} and Hughes~\cite{Hughes_1986} proved
\emph{a priori} error estimates for \gls{SUPG} stabilization of the
scalar advection-diffusion equation showing a half-power reduction
from the optimal rate in the advection-dominated regime, but in cases
such as this where the true solution is sufficiently smooth, optimal
convergence is often obtained for \gls{SUPG} formulations.  We
therefore do not have a satisfactory explanation of these results at
the present time.  Log-log plots for the \qq{2}{1} elements are given
in Fig.~\ref{fig:supg_adv} for the advection-dominated case, and
Fig.~\ref{fig:supg_diff} for the diffusion-dominated case.

\subsection{Jeffery--Hamel flow\label{sec:jh}}
Viscous, incompressible flow in a two-dimensional wedge, frequently
referred to as Jeffery--Hamel flow, is described in many references
including: the original works by Jeffery~\cite{Jeffery_1915} and
Hamel~\cite{Hamel_1916}, detailed analyses of the
problem~\cite{Rosenhead_1940,Haines_2011}, and treatment in fluid
mechanics textbooks~\cite{White_1991,Panton_2012}.
In addition to fundamental fluid mechanics research, the Jeffery--Hamel
flow is also of great utility as a validation tool for
finite difference, finite element, and related numerical codes
designed to solve the \gls{INS} equations under more general conditions.
The \ns module contains a regression test corresponding to a particular
Jeffery--Hamel flow configuration, primarily because it provides a verification method
which is independent from the \gls{MMS} described in the preceding sections.

The Jeffery--Hamel solution corresponds to flow constrained to a
wedge-shaped region defined by: $r_1 \leq r \leq r_2$, $-\alpha \leq
\theta \leq \alpha$.  The origin ($r=0$) is a singular point of the
flow, and is therefore excluded from numerical computations. The
governing equations are the \gls{INS} mass and momentum conservation
equations in cylindrical polar coordinates, under the assumption that
the flow is purely radial ($u_{\theta} = 0$) in nature. The boundary
conditions are no slip on the solid walls ($u_r(r, \pm\alpha)=0$) of
the wedge, and symmetry about the centerline ($\theta=0$) of the
channel. Under these assumptions, it is possible to show~\cite{Peterson_2016}
that the solution of the governing equations is:
\begin{align}
  u_r &= \frac{\lambda}{r}f(\eta)
  \\
  p &= p^{\ast} + \frac{2\mu \lambda}{r^2} (f(\eta) + K)
\end{align}
where $\lambda$ is a constant with units of (length)$^2$/time which is
proportional to the centerline (maximum) velocity, $\eta \equiv \frac{\theta}{\alpha}$
is the non-dimensionalized angular coordinate, $p^{\ast}$ is an arbitrary constant,
$f(\eta)$ is a dimensionless function that satisfies the third-order nonlinear \gls{ODE}
\begin{align}
  \label{eqn:ode3}
  f''' + 2\text{Re}\,\alpha f\!f' + 4\alpha^2f' &= 0 \quad 0 < \eta < 1
  \\
  f(0)  &= 1 \quad \text{(centerline velocity)} \\
  f'(0) &= 0 \quad \text{(centerline symmetry)} \\
  \label{eqn:f1}
  f(1)  &= 0 \quad \text{(no slip)}
\end{align}
where $\text{Re} \equiv \frac{\lambda \alpha}{\nu}$ is the Reynolds number
and $K$ is a dimensionless constant which is given in terms of $f$ by
\begin{align}
  \label{eqn:K_integrated_3}
  K  = \frac{1}{4\alpha^2} \left(\frac{1}{2} f'(1)^2 - \frac{\alpha \text{Re}}{3} - 2\alpha^2\right)
\end{align}

Since~\eqref{eqn:ode3} has no simple closed-form solution in general,
the standard practice is to solve this \gls{ODE} numerically to high accuracy and
tabulate the values for use as a reference solution
in estimating the accuracy of a corresponding 2D finite element
solution. The pressure is only determined up to an arbitrary constant
in both the analytical and numerical formulations of this problem (we impose Dirichlet
velocity boundary conditions on all boundaries in the finite element formulation). Therefore, in practice we
pin a single value of the pressure on the inlet centerline to zero to
constrain the non-trivial nullspace. Finally, we note that although
the Jeffery--Hamel problem is naturally posed in cylindrical coordinates,
we solve the problem using a Cartesian space finite element formulation
in the \ns module.

\begin{figure}[htpb]
  \centering
  \begin{subfigure}[t]{.49\linewidth}
    \includegraphics[width=\linewidth,viewport=400 0 2600 2000,clip=true]{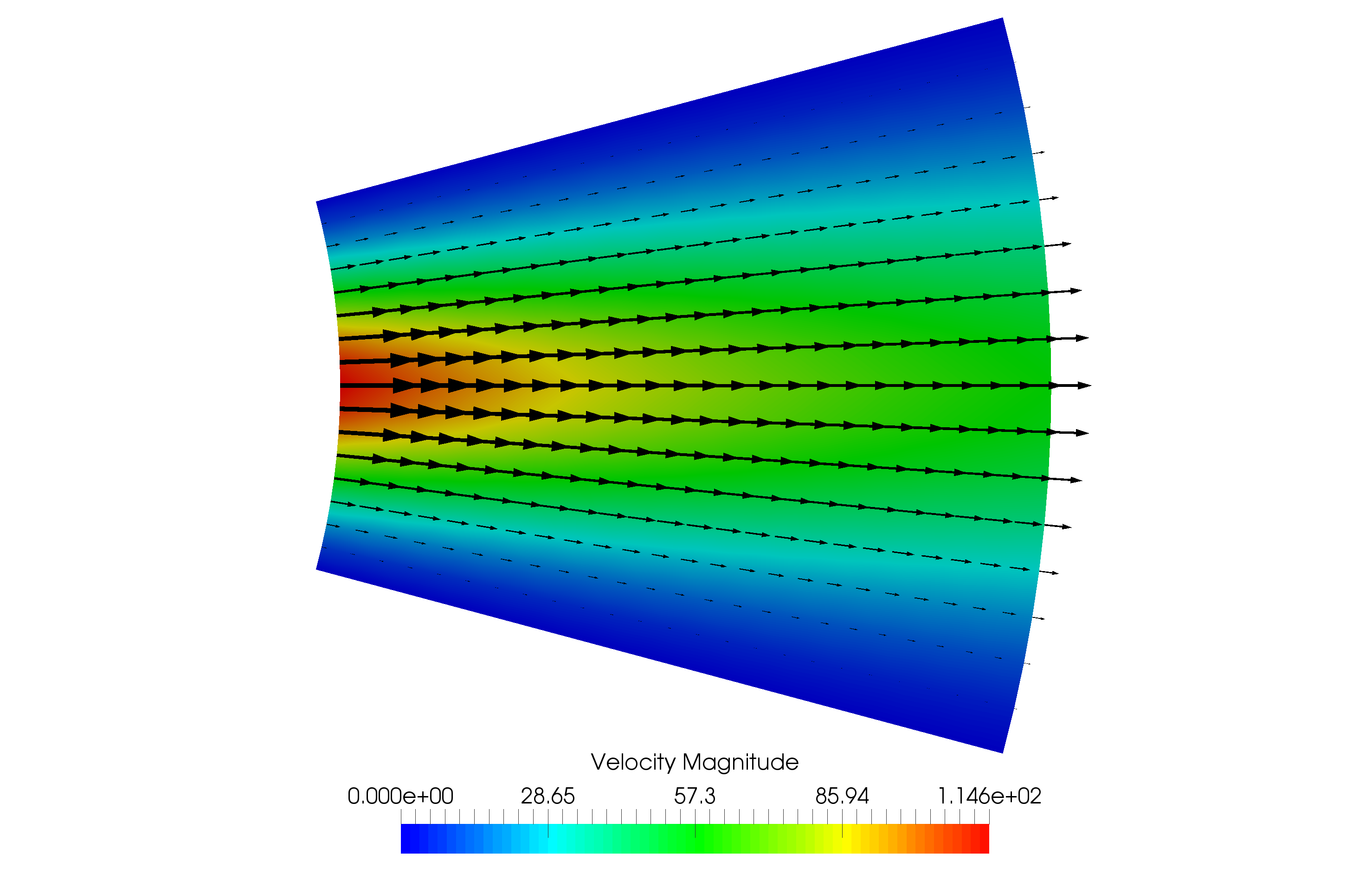}
    \caption{Velocity magnitude and vectors.\label{fig:jh_velocity}}
  \end{subfigure}
  \begin{subfigure}[t]{.49\linewidth}
    \includegraphics[width=\linewidth,viewport=400 0 2600 2000,clip=true]{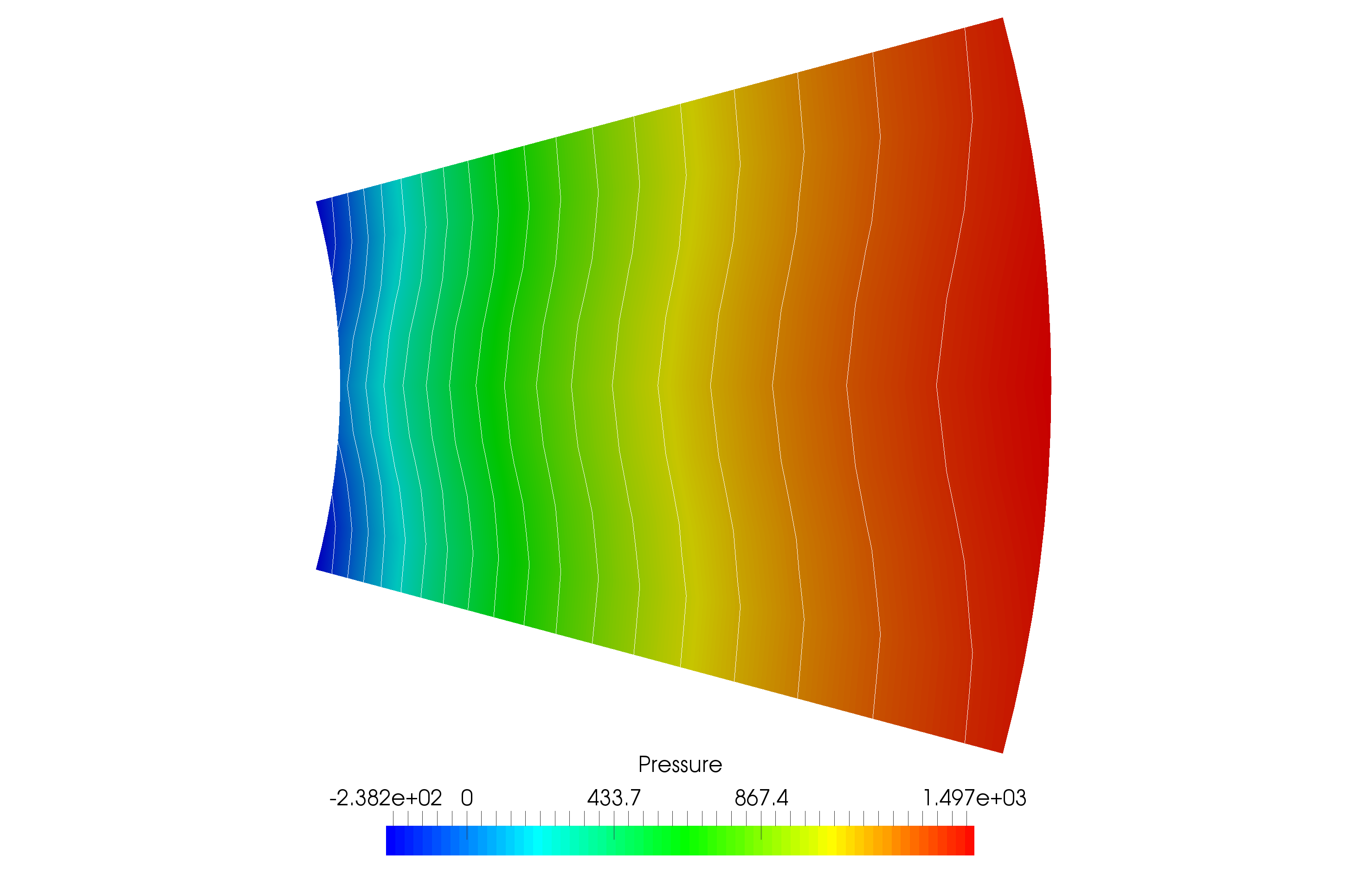}
    \caption{Pressure contours.\label{fig:jh_pressure}}
  \end{subfigure}
  \caption{Plots of Jeffery--Hamel flow showing (\subref{fig:jh_velocity})~velocity magnitude and vectors, and
    (\subref{fig:jh_pressure})~pressure contours for the case $\alpha=15^{\circ}$ and $\text{Re}=30$.\label{fig:jh_results}}
\end{figure}

The Jeffery--Hamel velocity and pressure solutions for the particular
case $\alpha=15^{\circ}$, $\text{Re}=30$, $r_1=1$, and $r_2=2$ are shown in
Fig.~\ref{fig:jh_results}. As expected from the exact solution,
the flow is radially self-similar, and the centerline maximum value is proportional to $1/r$. The pressure field has a corresponding
$-1/r^2$ dependence on the radial distance, since $K\approx -9.7822146449$ is
negative in this configuration. We note that the radial flow solution is linearly
unstable for configurations in which $\text{Re}\,\alpha > 10.31$, since the
boundary layers eventually separate and recirculation regions form due to the
increasingly adverse pressure gradient. Our test case is in the linearly stable regime.

\begin{figure}[htpb]
  \centering
  \begin{subfigure}[t]{.49\linewidth}
    \includegraphics[width=\linewidth]{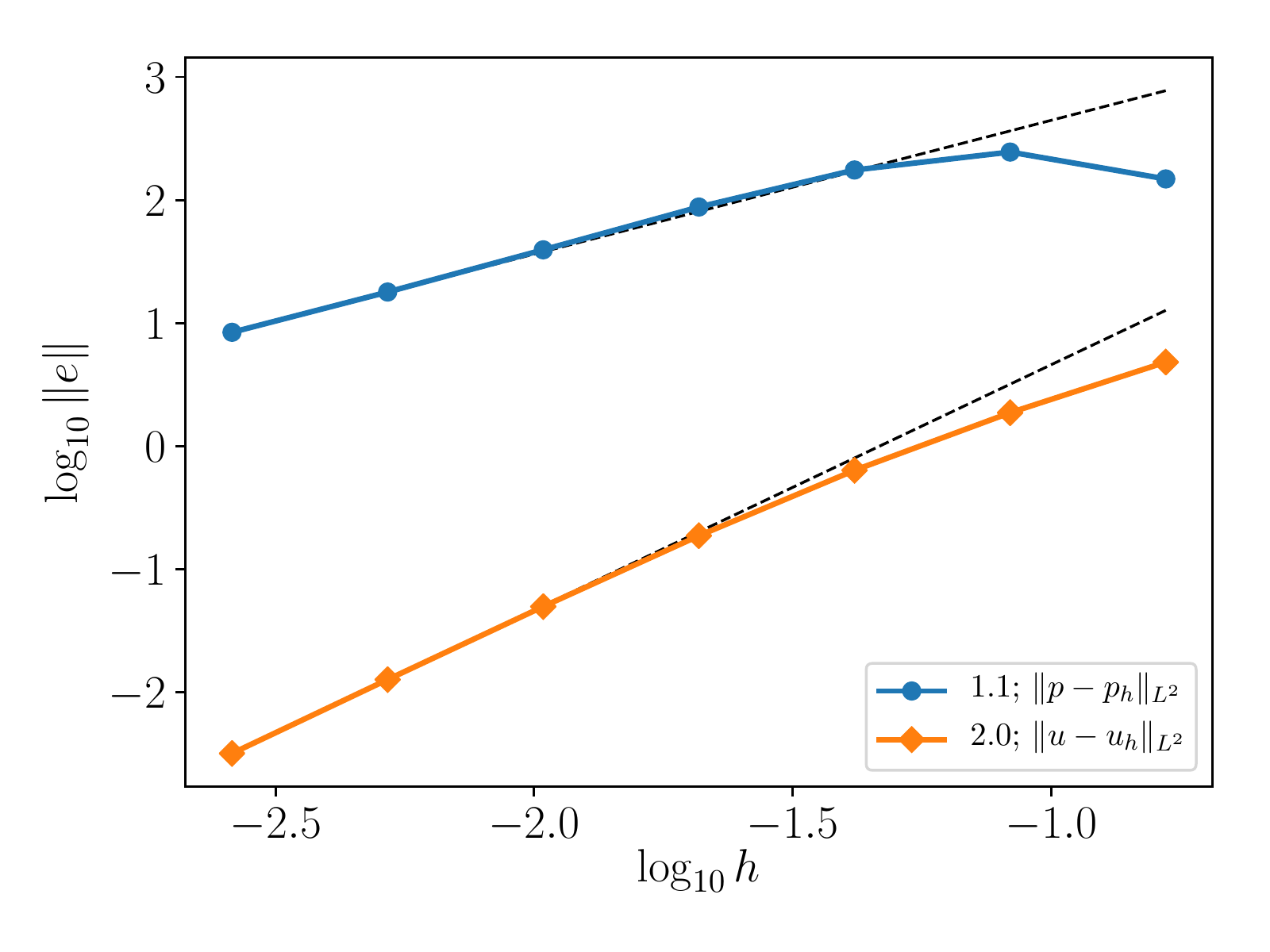}
    \caption{\gls{PSPG}-stabilized \qq{1}{1} elements.\label{fig:jh_q1q1_convergence}}
  \end{subfigure}
  \begin{subfigure}[t]{.49\linewidth}
    \includegraphics[width=\linewidth]{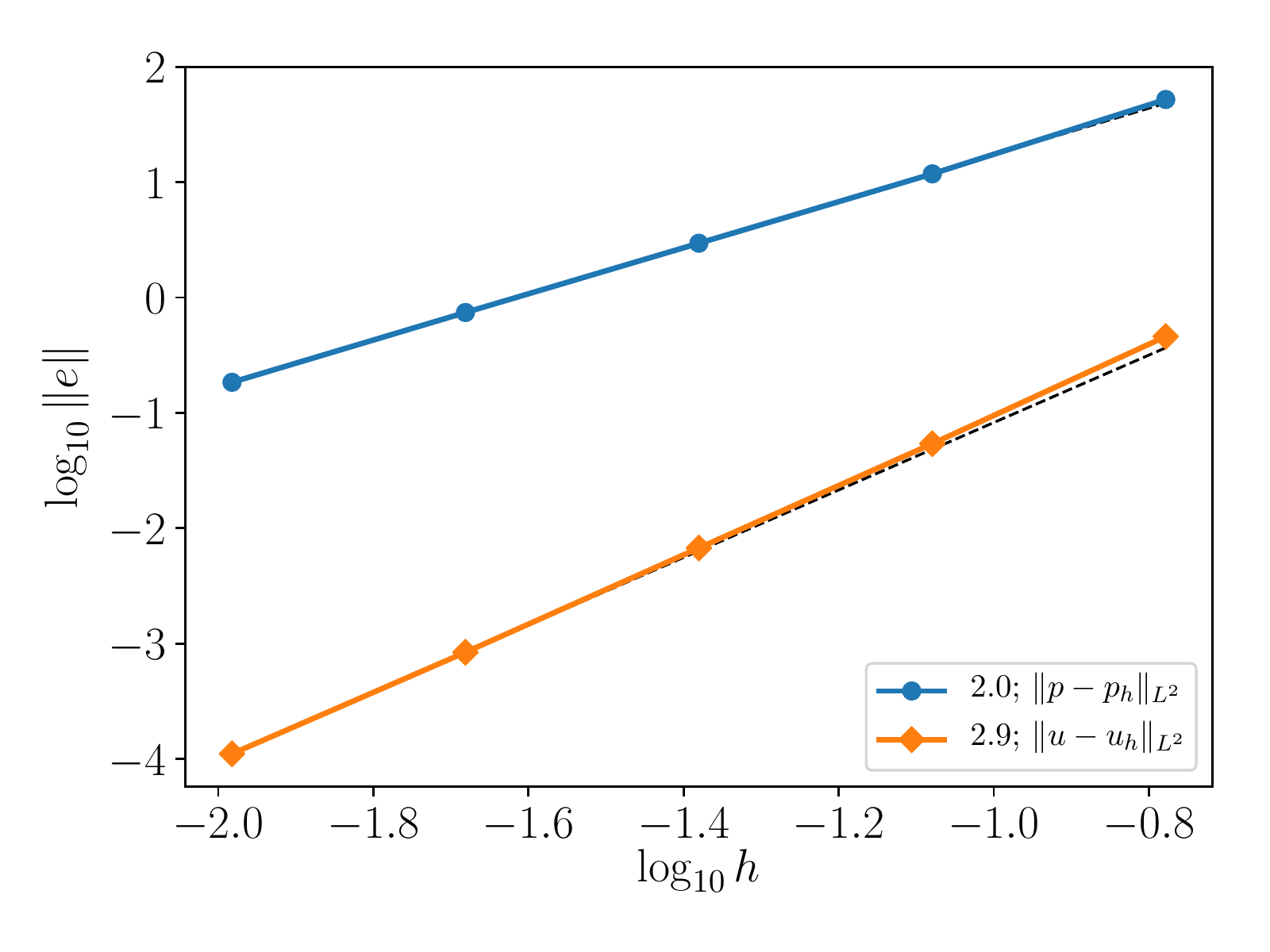}
    \caption{\qq{2}{1} elements.\label{fig:jh_q2q1_convergence}}
  \end{subfigure}
  \caption{Convergence rates of (\subref{fig:jh_q1q1_convergence})~\gls{PSPG}-stabilized \qq{1}{1} and
    (\subref{fig:jh_q2q1_convergence})~\qq{2}{1} discretizations of the Jeffery--Hamel problem for
    $\alpha=15^{\circ}$ and $\text{Re}=30$.\label{fig:jh_convergence_all}}
\end{figure}

The convergence results for \gls{PSPG}-stabilized \qq{1}{1} and
unstabilized \qq{2}{1} finite element formulations on a sequence of
uniformly-refined grids for this problem are shown in
Fig.~\ref{fig:jh_convergence_all}. Optimal rates in both the velocity
and pressure are obtained for the unstabilized case, which is consistent
with the diffusion-dominated results of \S\ref{sec:verif}. In the \gls{PSPG}-stabilized
case, the velocity converges at the optimal order, while the $\|p-p_h\|_{L^2}$
norm once again exhibits first-order convergence as in \S\ref{sec:verif} for
the diffusion-dominated problem. An interesting extension to
these results would be to consider \emph{converging} flow
(i.e.\ choosing $\lambda$ and therefore $\text{Re}< 0$). The converging
flow configuration is linearly stable for all (negative) Reynolds numbers, and
thus, for large enough negative Reynolds numbers, the numerical solution may require
\gls{SUPG} stabilization.

\section{Representative results\label{sec:results}}
In this section, we move beyond the verification tests of
\S\ref{sec:verification_tests} and employ the various finite element
formulations of the \ns module on several representative applications.
In \S\ref{sec:lid}, we report results for a ``non-leaky'' version of the
lid-driven cavity problem at Reynolds numbers 10$^3$ and $5\times 10^3$, which are
in good qualitative agreement with benchmark results. Then, in
\S\ref{sec:axi}, we introduce a non-Cartesian coordinate system and
investigate diffusion and advection-dominated flows in a 2D expanding
nozzle. Finally, in \S\ref{sec:sphere}, we solve the classic problem
of 3D advection-dominated flow over a spherical obstruction in a
narrow channel using a hybrid mesh of prismatic and tetrahedral
elements.

\subsection{Lid-driven cavity\label{sec:lid}}
In the classic lid-driven cavity problem, the flow is induced by a moving ``lid'' on top of
a ``box'' containing an incompressible, viscous fluid. In the continuous (non-leaky) version of the problem, the lid velocity is set to zero
at the top left and right corners of the domain for consistency with the no-slip boundary conditions on the
sides of the domain. In the present study, the computational domain is taken to be $\Omega=[0,1]^2$,
and the lid velocity is set to:
\begin{align}
  \label{eq:lid_function}
  u_1 = 4x(1-x)
\end{align}
while the viscosity $\mu$ is varied to simulate flows at different
Reynolds numbers. Since Dirichlet boundary conditions are imposed on
the entire domain, a single value of the pressure is pinned to zero in
the bottom left-hand corner.  We consider two representative cases of
$\text{Re} = 10^3$ and $\text{Re} = 5\times 10^3$ in the present work, and
both simulations employ an \gls{SUPG}/\gls{PSPG}-stabilized finite
element discretization on a structured grid of \qq{1}{1} elements.

\begin{figure}[htpb]
  \centering
  \begin{subfigure}[t]{.49\linewidth}
    \includegraphics[width=\linewidth,viewport=600 100 2250 1900,clip=true]{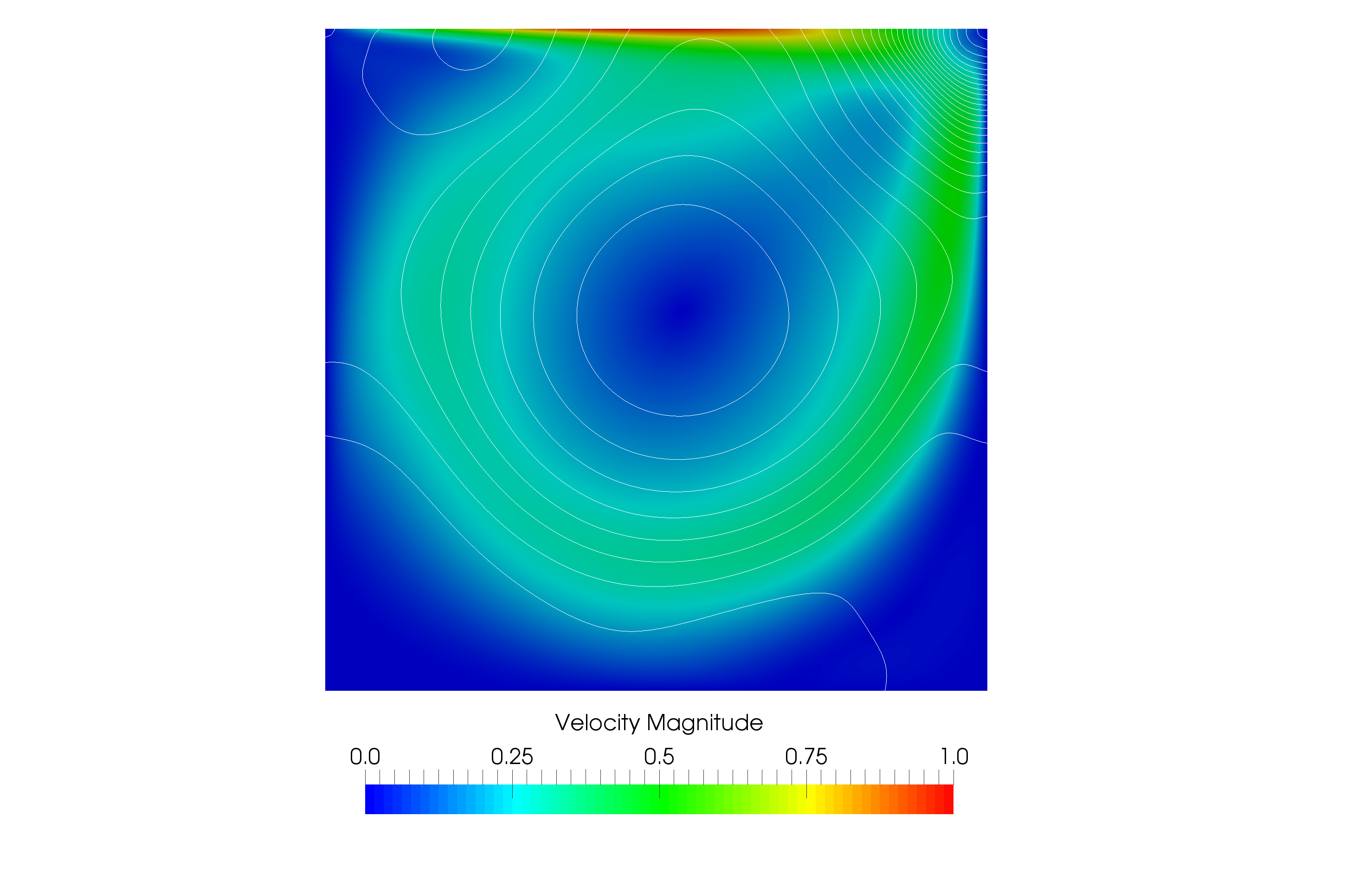}
    \caption{$\text{Re}=10^3$, $h=1/128$.\label{fig:lid_re1k_pressure_contours}}
  \end{subfigure}
  \begin{subfigure}[t]{.49\linewidth}
    \includegraphics[width=\linewidth,viewport=600 100 2250 1900,clip=true]{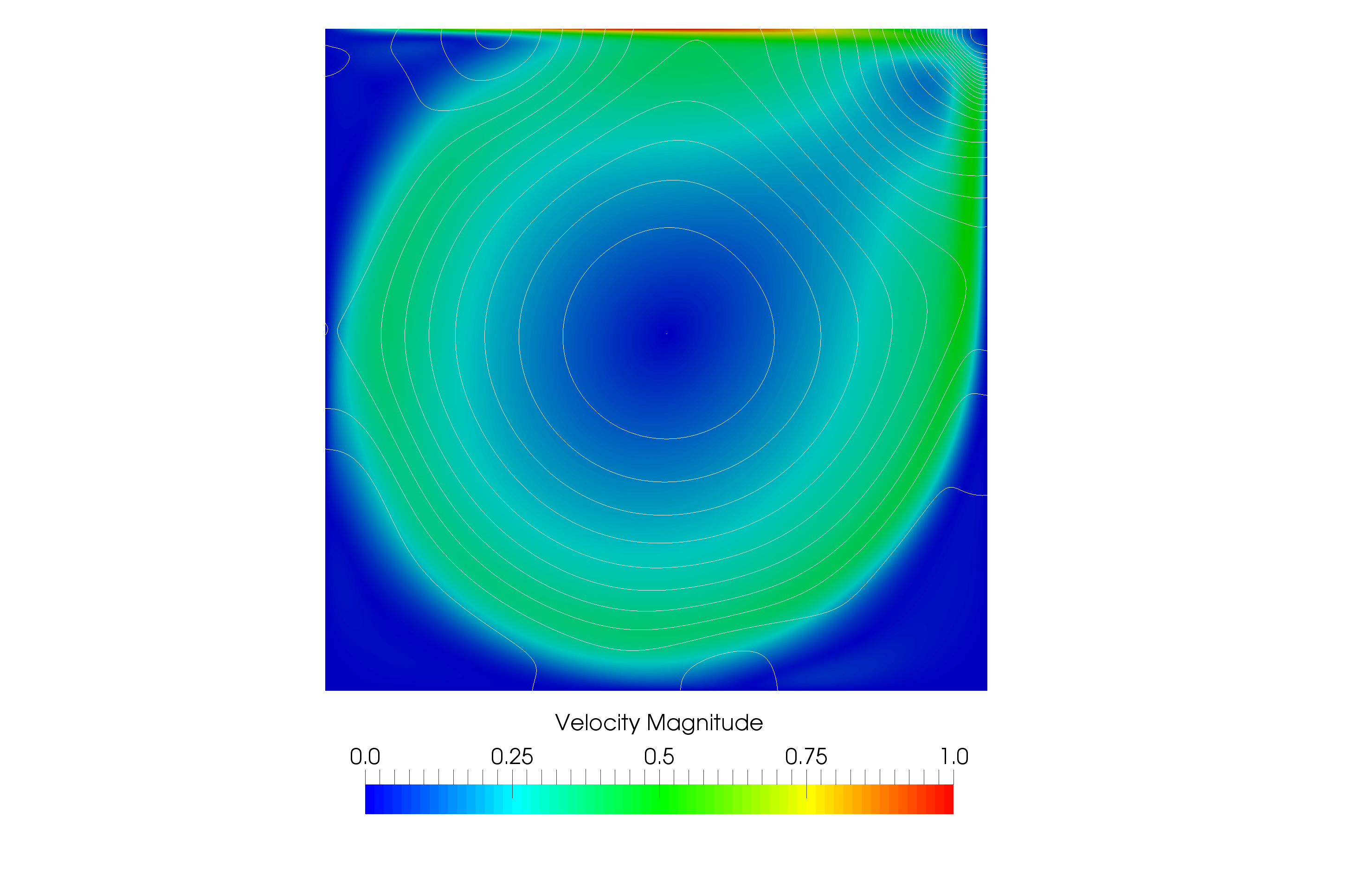}
    \caption{$\text{Re}=5\times 10^3$, $h=1/256$.\label{fig:lid_re5k_pressure_contours}}
  \end{subfigure}
  \caption{Pressure isobars overlaid on velocity magnitude contour plots for the
    lid-driven cavity problem with
    (\subref{fig:lid_re1k_pressure_contours})~$\text{Re}=10^3$, $h=1/128$, and
    (\subref{fig:lid_re5k_pressure_contours})~$\text{Re}=5\times 10^3$, $h=1/256$. These
    results are qualitatively similar to those discussed
    in~\cite{Olshanskii_2002} and other reference
    works.\label{fig:lid_driven_pressure_contours}}
\end{figure}

\begin{figure}[htpb]
  \centering
  \begin{subfigure}[t]{.49\linewidth}
    \includegraphics[width=\linewidth,viewport=350 100 1100 800,clip=true]{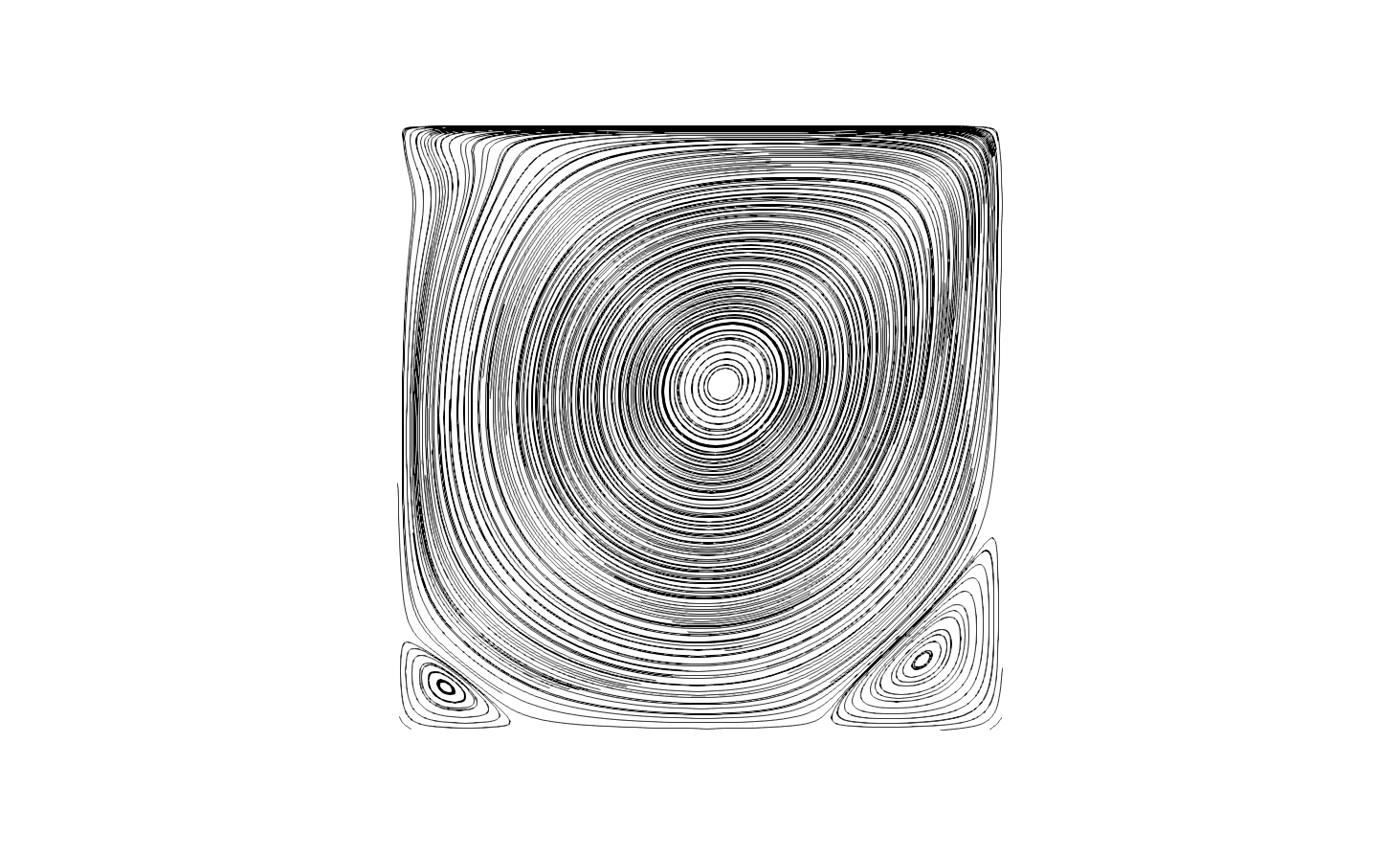}
    \caption{$\text{Re}=10^3$, $h=1/128$.\label{fig:lid_re1k_streamlines}}
  \end{subfigure}
  \begin{subfigure}[t]{.49\linewidth}
    \includegraphics[width=\linewidth,viewport=400 100 1150 800,clip=true]{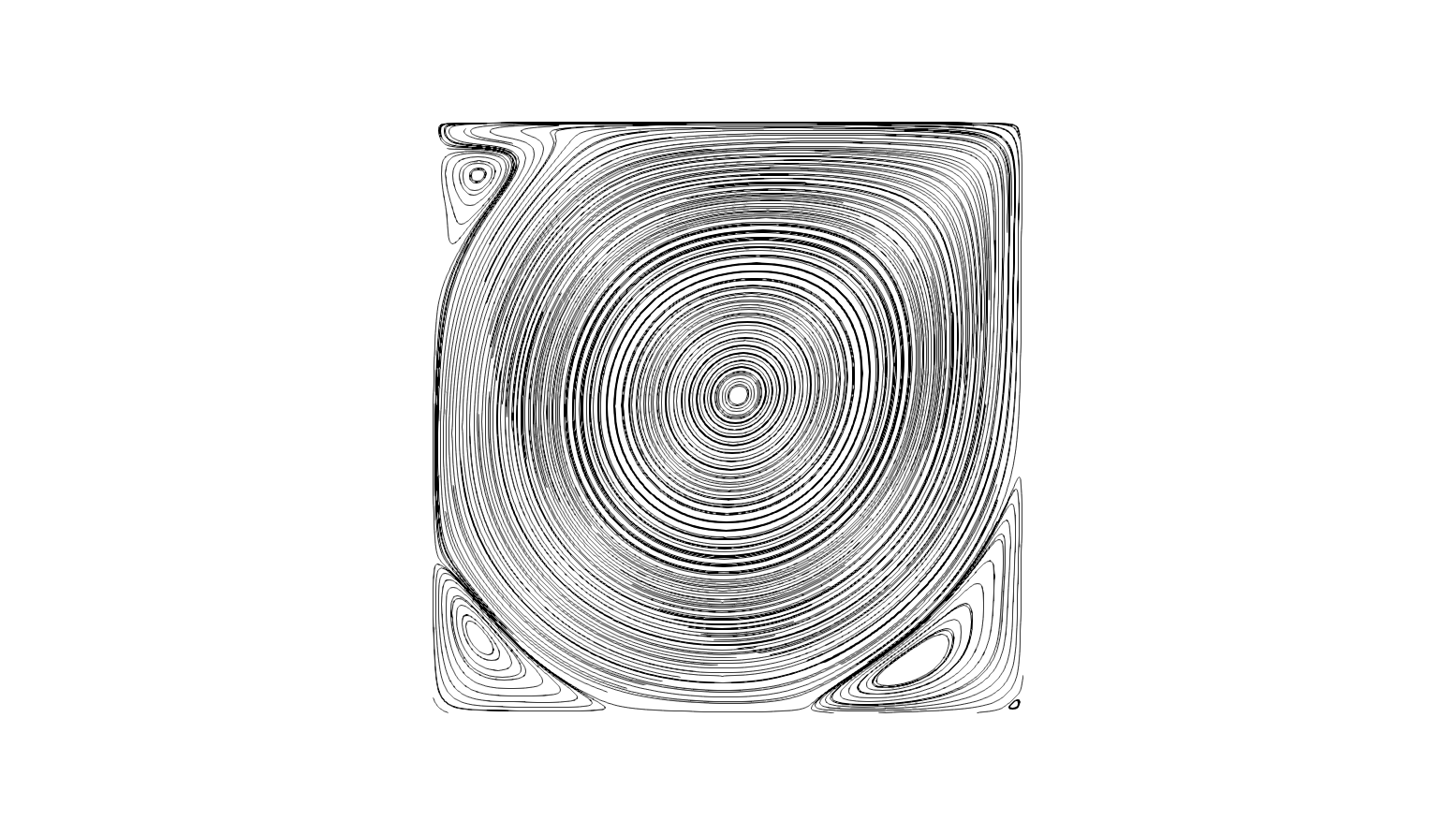}
    \caption{$\text{Re}=5\times 10^3$, $h=1/256$.\label{fig:lid_re5k_streamlines}}
  \end{subfigure}
  \caption{Streamlines for the lid-driven cavity case with (\subref{fig:lid_re1k_streamlines})~$\text{Re}=10^3$, $h=1/128$, and
    (\subref{fig:lid_re5k_streamlines})~$\text{Re}=5\times 10^3$, $h=1/256$. The third region of recirculating flow present in the
    $\text{Re}=5\times 10^3$ result is consistent with other lid-driven cavity reference solutions.\label{fig:lid_driven_streamlines}}
\end{figure}

Figs.~\ref{fig:lid_re1k_pressure_contours}
and~\ref{fig:lid_re5k_pressure_contours} show the velocity magnitude
and pressure contours at steady state for the two cases, using a mesh
size of $h=1/128$ in the $\text{Re} = 10^3$ case, and $h=1/256$ in the
$\text{Re} = 5 \times 10^3$ case. Solutions of the leaky version of the
lid-driven cavity problem exhibit strong pressure singularities (and
velocity discontinuities) in the upper corners of the domain. In
contrast, for the non-leaky version of the problem solved here, the
velocity field is continuous and therefore the solution has higher
global regularity. Nevertheless, the pressure still varies rapidly in
the upper right-hand corner of the domain (as evidenced by the
closely-spaced contour lines), and fine-scale flow structures develop
at high Reynolds numbers.  Flow streamlines for the two problems are
shown in Figs.~\ref{fig:lid_re1k_streamlines}
and~\ref{fig:lid_re5k_streamlines}. The locations and strengths of the
vortices, including the appearance of a third recirculation region at
$\text{Re} = 5\times 10^3$, are consistent with the results reported
in~\cite{Olshanskii_2002} and other reference works.

\subsection{Axisymmetric channel\label{sec:axi}}
To derive the axisymmetric version of the governing equations, we start by defining
the vector $\vec{v}$ in cylindrical coordinates as
\begin{align}
  \label{eqn:cylindrical_vector}
  \vec{v} = v_r \unit{e}_r + v_{\theta} \unit{e}_{\theta} + v_z \unit{e}_z
\end{align}
for unit vectors $\unit{e}_r$, $\unit{e}_{\theta}$, $\unit{e}_z$.
The divergence and gradient of $\vec{v}$ in cylindrical coordinates are then given by:
\begin{align}
  \label{eqn:cylindrical_divergence}
  \Div \vec{v} &\equiv
  \frac{1}{r}\frac{\partial \left(r v_r\right)}{\partial r}
  +
  \frac{1}{r}\frac{\partial v_{\theta}}{\partial \theta}
  +
  \frac{\partial v_z}{\partial z}
  \\[6pt]
  \label{eqn:vec_gradient}
  \grad \vec{v} &\equiv
  \begin{bmatrix}
    \frac{\partial v_r}{\partial r}      & \frac{1}{r}\left(\frac{\partial v_r}{\partial \theta} - v_{\theta} \right) & \frac{\partial v_r}{\partial z} \\[6pt]
    \frac{\partial v_{\theta}}{\partial r} & \frac{1}{r}\left(\frac{\partial v_{\theta}}{\partial \theta} + v_{r} \right) & \frac{\partial v_{\theta}}{\partial z} \\[6pt]
    \frac{\partial v_{z}}{\partial r}     & \frac{1}{r}\frac{\partial v_{z}}{\partial \theta} & \frac{\partial v_{z}}{\partial z}
  \end{bmatrix}
\end{align}
while the gradient and Laplacian of a scalar $g$ are given by
\begin{align}
  \label{eqn:cylindrical_gradient}
  \grad g &\equiv \frac{\partial g}{\partial r} \unit{e}_r
  +
  \frac{1}{r}\frac{\partial g}{\partial \theta} \unit{e}_{\theta}
  +
  \frac{\partial g}{\partial z} \unit{e}_{z}
  \\
  \label{eqn:cylindrical_lap_scalar}
  \Lap g &\equiv \frac{1}{r} \frac{\partial}{\partial r} \left(r \frac{\partial g}{\partial r} \right)
  +
  \frac{1}{r^2}\frac{\partial^2 g}{\partial \theta^2}
  +
  \frac{\partial^2 g}{\partial z^2}
\end{align}
respectively. Combining these definitions allows us to define the vector Laplacian operator
in cylindrical coordinates as:
\begin{align}
  \label{eqn:cylindrical_vector_laplacian}
  \text{Lap}\, \vec {u} \equiv
  \Div (\grad \vec{u})^T =
  \begin{bmatrix}
    \text{Lap}\, u_r \\[6pt]
    \text{Lap}\, u_{\theta} \\[6pt]
    \text{Lap}\, u_z
  \end{bmatrix}
  +
  \begin{bmatrix}
    - \frac{2}{r^2}\frac{\partial u_{\theta}}{\partial \theta} - \frac{u_r}{r^2}
    \\[6pt]
    \phantom{-}\frac{2}{r^2}\frac{\partial u_r}{\partial \theta} - \frac{u_{\theta}}{r^2}
    \\[6pt]
    0
  \end{bmatrix}
\end{align}

The axisymmetric formulation then proceeds by
substituting~\eqref{eqn:cylindrical_vector}--\eqref{eqn:cylindrical_vector_laplacian}
into e.g.\ the semi-discrete component
equations~\eqref{eqn:semi_discrete_momentum_k} and \eqref{eqn:semi_discrete_mass}
and then assuming the velocity, pressure and body force terms satisfy
\begin{align}
  u_{\theta} = \frac{\partial u_r}{\partial \theta} = \frac{\partial u_z}{\partial \theta} = \frac{\partial p}{\partial \theta} = f_{\theta} = 0
\end{align}
Under these assumptions, we can ignore the $\theta$-component of~\eqref{eqn:semi_discrete_momentum_k} since it is trivially satisfied.
The remaining component equations in Laplace form are (neglecting the stabilization terms and
assuming natural boundary conditions for brevity):
\begin{alignat}{2}
  \nonumber
  &\int_{\Omega'} \left[
  \rho\left(\frac{\partial u_r}{\partial t} +
    u_r \frac{\partial u_r}{\partial r} + u_z \frac{\partial u_r}{\partial z} \right) \varphi_i
    - p \frac{\partial \varphi_i}{\partial r} \highlight{\mbox{}- p\frac{\varphi_i}{r}}
    - f_r \varphi_i
    \right.
  \\
  \label{eqn:axi_mom_r_simplified}
  &\left.
    \qquad \qquad \mbox{}+
    \mu\left( \frac{\partial u_r}{\partial r}\frac{\partial \varphi_i}{\partial r} +
      \frac{\partial u_r}{\partial z} \frac{\partial \varphi_i}{\partial z}
      \highlight{\mbox{}+ \frac{u_r \varphi_i}{r^2}}\right)
  \right] \text{d} \Omega' = 0, && \quad i=1,\ldots,N
  \\
  \nonumber
  &\int_{\Omega'} \left[
  \rho\left(\frac{\partial u_z}{\partial t} +
    u_r \frac{\partial u_z}{\partial r} + u_z \frac{\partial u_z}{\partial z} \right) \varphi_i
    -p\frac{\partial \varphi_i}{\partial z}
    - f_z \varphi_i
    \right.
  \\
  \label{eqn:axi_mom_z_simplified}
  &\left.
    \qquad \qquad \mbox{}+
    \mu\left(
      \frac{\partial u_z}{\partial r} \frac{\partial \varphi_i}{\partial r}
      + \frac{\partial u_z}{\partial z} \frac{\partial \varphi_i}{\partial z}\right)
    \right] \text{d} \Omega' = 0, && \quad i=1,\ldots,N
  \\
  \label{eqn:axi_mass_simplified}
  &\int_{\Omega'} \left(
    \frac{\partial u_r}{\partial r}
    \highlight{\mbox{} + \frac{u_r}{r}} +
  \frac{\partial u_z}{\partial z}
    \right) \psi_i \; \text{d} \Omega' = 0, && \quad i=1,\ldots,M
\end{alignat}
Upon inspection, we note
that~\eqref{eqn:axi_mom_r_simplified}--\eqref{eqn:axi_mass_simplified}
are equivalent to the two-dimensional Cartesian \gls{INS} equations,
except for the three additional terms which have been highlighted in
red, and the factor of $r$ which arises from the $\text{d} \Omega' = r \, \text{d}r\, \text{d}z$
terms.  Therefore, to simulate axisymmetric flow in the \ns module, we
simply use the two-dimensional Cartesian formulation, treat the $u_1$ and $u_2$
components of the velocity as $u_r$ and $u_z$, respectively, append the
highlighted terms, and multiply by $r$ when computing the element
integrals.

To test the \ns module's axisymmetric simulation capability, we
consider flow through a conical diffuser with inlet radius $R=1/2$ and
specified velocity profile:
\begin{align}
  u_r &= 0
  \\
  \label{eq:inlet_profile_axi}
  u_z &= 1 - 4r^2
\end{align}
The conical section expands to a radius of $R=1$ at $z=1$, and then the
channel cross-section remains constant until the outflow boundary at
$z=4$. No-slip Dirichlet boundary conditions are applied along the outer ($r=R$)
wall of the channel, and a no-normal-flow ($u_r=0$) Dirichlet boundary condition
is imposed along the channel centerline to enforce symmetry about
$r=0$. As is standard for these types of problems, natural boundary
conditions ($\neumannvec=\vec{0}$) are imposed at the outlet.

We note that natural boundary conditions are used here for convenience, and not
because they are particularly well-suited to modeling the ``real'' outflow boundary.
As discussed in \S\ref{sec:bcs}, a more realistic approach would be to employ a longer channel region to ensure that
the chosen outlet boundary conditions do not adversely affect the upstream
characteristics of the flow. Since our intent here is primarily to demonstrate
and exercise the capabilities of the \ns module, no effort is made to minimize
the effects of the outlet conditions. As a final verification step, we observe that
integrating~\eqref{eq:inlet_profile_axi} over the inlet yields a volumetric flow rate of
\begin{align}
  Q = \int_{0}^{2\pi}\!\!\! \int_0^R  \left(1 - 4r^2\right) r\, \text{d}r \, \text{d}\theta = \frac{\pi}{8}
\end{align}
for this case. In this example, we employ a \texttt{VolumetricFlowRate} Postprocessor to ensure that the mass flowing
in the inlet exactly matches (to within floating point tolerances) the mass flowing out the outlet,
thereby confirming the global conservation property of the \gls{FEM} for the axisymmetric case.

The characteristics of the computed solutions depend strongly on the Reynolds
number for this problem, which is similar in nature to the classic backward-facing
step problem~\cite{Gresho_1993}. The flow must decelerate as it flows into
the larger section of the channel, and this deceleration may be accompanied by a
corresponding rise in pressure, also known as an adverse pressure gradient.
If the adverse pressure gradient is
strong enough (i.e.\ at high enough Reynolds numbers) the flow will separate
from the channel wall at the sharp inlet, and regions of recirculating flow will form.
At low Reynolds numbers, on the other hand, the pressure gradient will remain
favorable, and flow separation will not occur.

We consider two distinct cases in the present work: a ``creeping flow'' case
with $\text{Re}=0.5$ ($\mu=1$) and an advection-dominated case with
$\text{Re}=10^3$ ($\mu=0.5 \times 10^{-3}$) based on the average inlet
velocity and inlet diameter. In the former case, we employ an LBB-stable \pp{2}{1}
finite element discretization on TRI6 elements with \gls{SUPG} and \gls{PSPG}
stabilization disabled. In the latter case, an equal-order \pp{1}{1} discretization
on TRI3 elements is employed, and both \gls{SUPG} and \gls{PSPG} stabilization are enabled.
In both cases, unstructured meshes are employed, and no special effort is expended
to grade the mesh into the near-wall/boundary layer regions of the channel.

\begin{figure}[htpb]
  \centering
  \begin{subfigure}[t]{.9\linewidth}
    \includegraphics[width=\linewidth,viewport=200 800 3600 2100,clip=true]{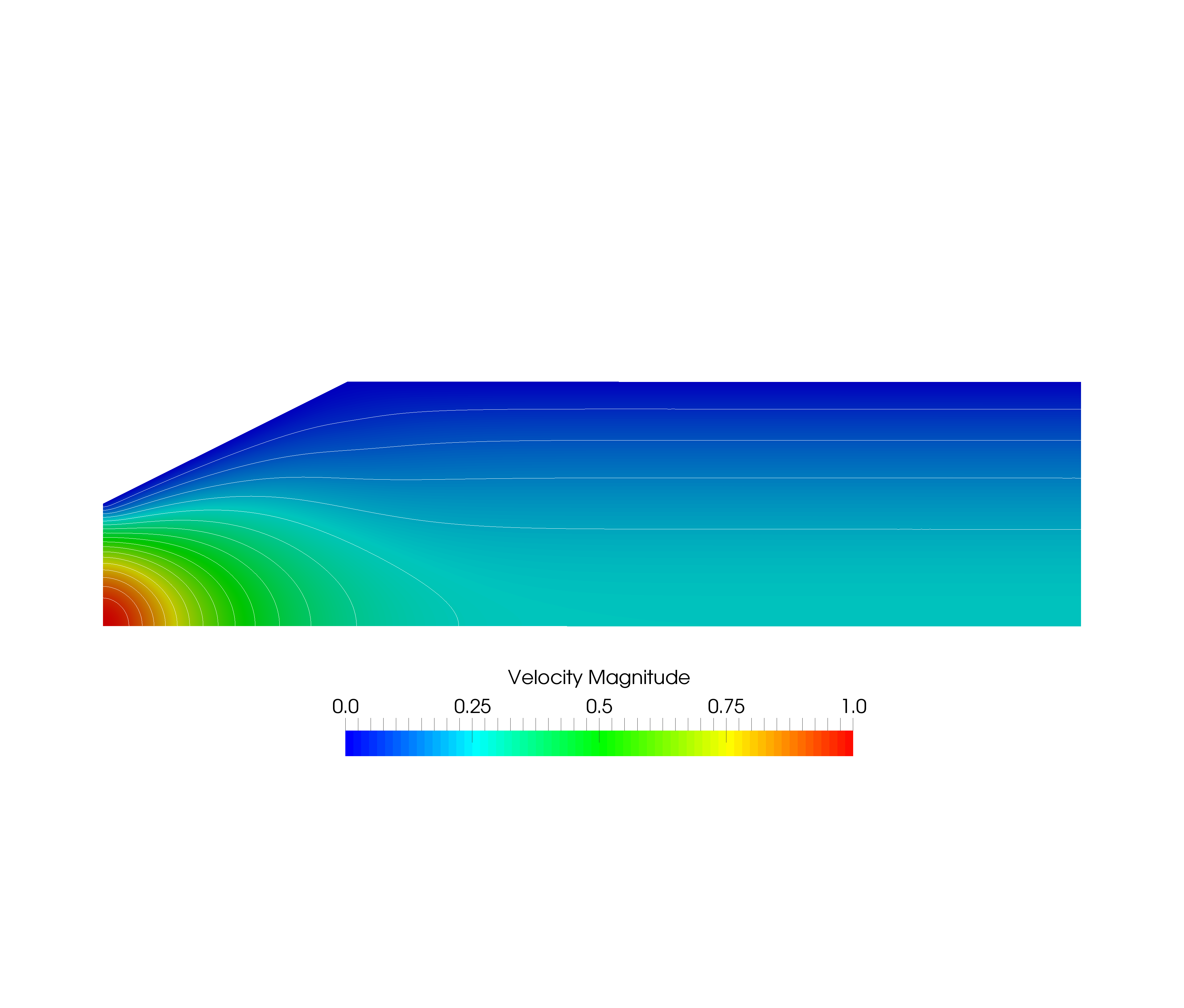}
    \caption{Velocity magnitude contours.\label{fig:cone_vel}}
  \end{subfigure}
  \\[6pt]
  \begin{subfigure}[t]{.9\linewidth}
    \includegraphics[width=\linewidth,viewport=200 800 3600 2100,clip=true]{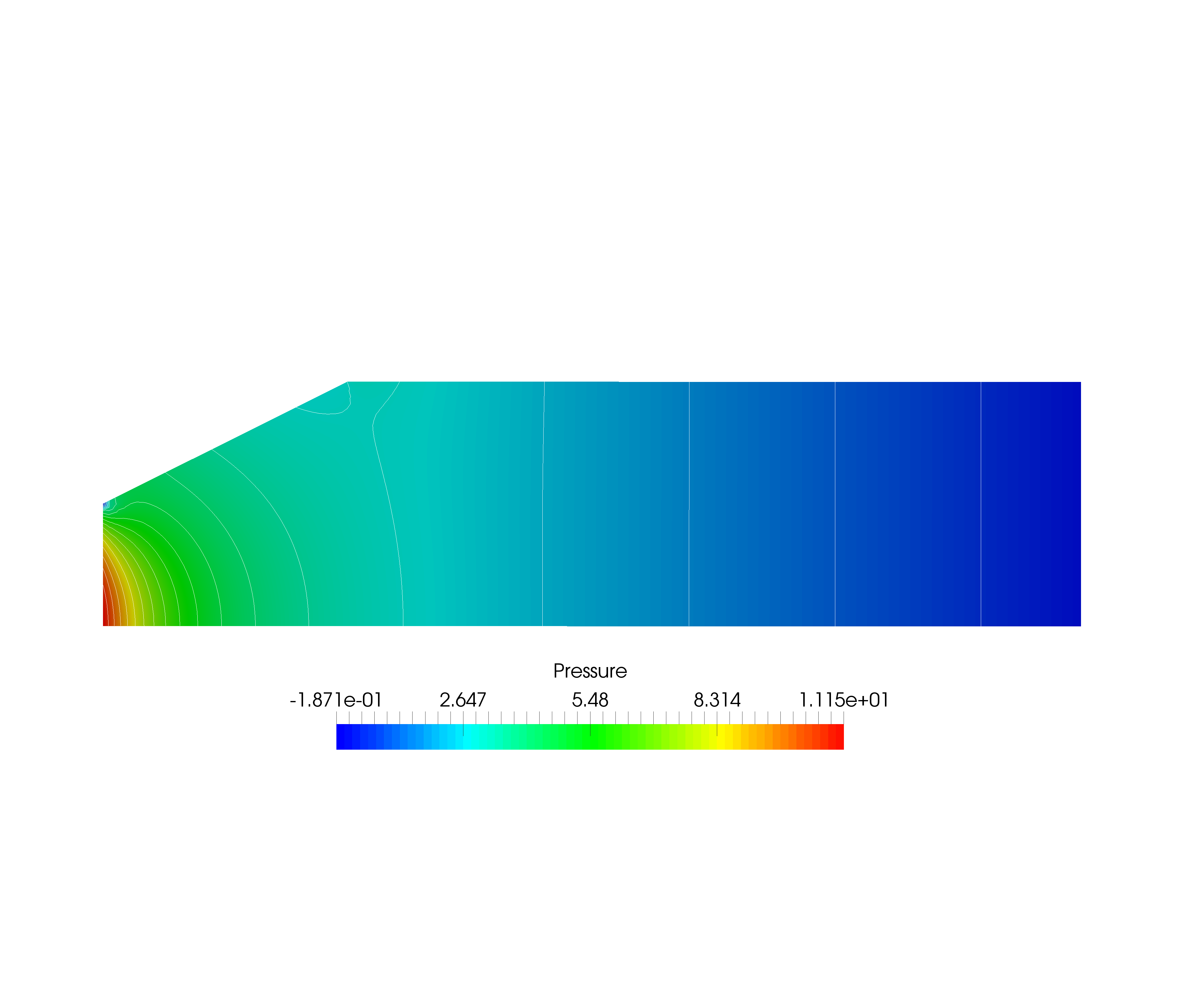}
    \caption{Pressure contours.\label{fig:cone_p}}
  \end{subfigure}
  \caption{Visualization of the axisymmetric channel flow problem
    showing the (\subref{fig:cone_vel})~velocity magnitude and
    (\subref{fig:cone_p})~pressure contours for the creeping flow case ($\text{Re}=0.5$) at steady
    state using an unstabilized \pp{2}{1} finite element discretization on a mesh with 10102 TRI6
    elements and 20513 nodes.\label{fig:axisymmetric_flow}}
\end{figure}

The steady state velocity magnitude and pressure contours for the
diffusion-dominated case are shown in Fig.~\ref{fig:axisymmetric_flow}.
In this case, we observe that the inlet velocity profile diffuses
relatively quickly to a quadratic profile at the outlet, and that
the flow remains attached throughout the channel. No boundary
layers form, and there is a favorable, i.e.\ $\frac{\partial p}{\partial z} < 0$,
pressure gradient throughout the channel, which is disturbed only
by the pressure singularity which forms at the sharp inlet. This pressure singularity
is also a characteristic of backward-facing step flows, and although it reduces the
overall regularity of the solution, it does not seem to have any serious
negative effects on the convergence of the numerical scheme. We also remark
that the downstream natural boundary conditions seem to have little effect, if
any, on this creeping flow case, since the flow becomes fully-developed well
before reaching the outlet.

\begin{figure}[htpb]
  \centering
  \begin{subfigure}[t]{.9\linewidth}
    \includegraphics[width=\linewidth,viewport=200 200 3600 1700,clip=true]{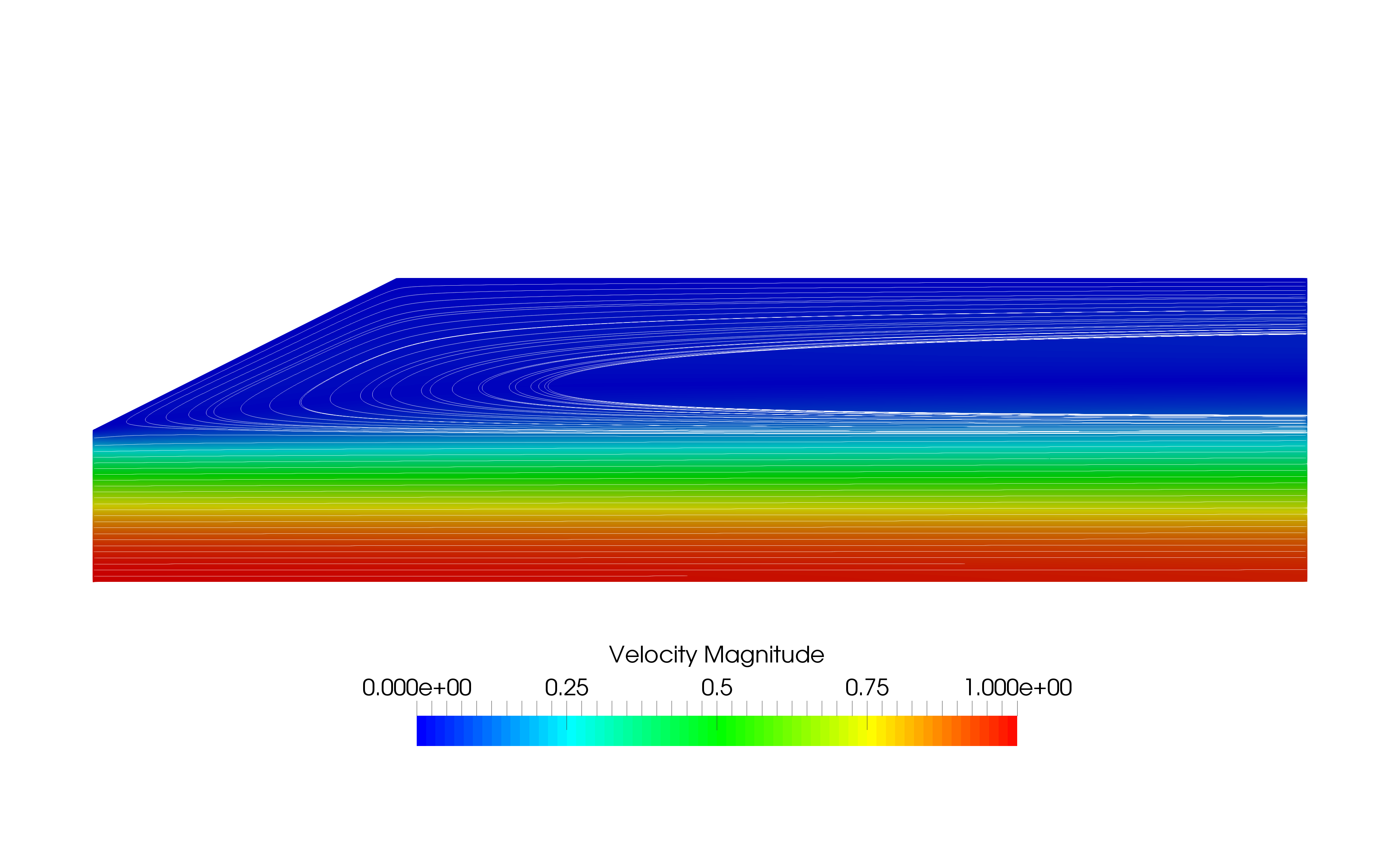}
    \caption{Velocity magnitude and streamtraces.\label{fig:RZ_cone_1k_streamtraces_86212_elem}}
  \end{subfigure}
  \\[6pt]
  \begin{subfigure}[t]{.9\linewidth}
    \includegraphics[width=\linewidth,viewport=70 250 3750 1750,clip=true]{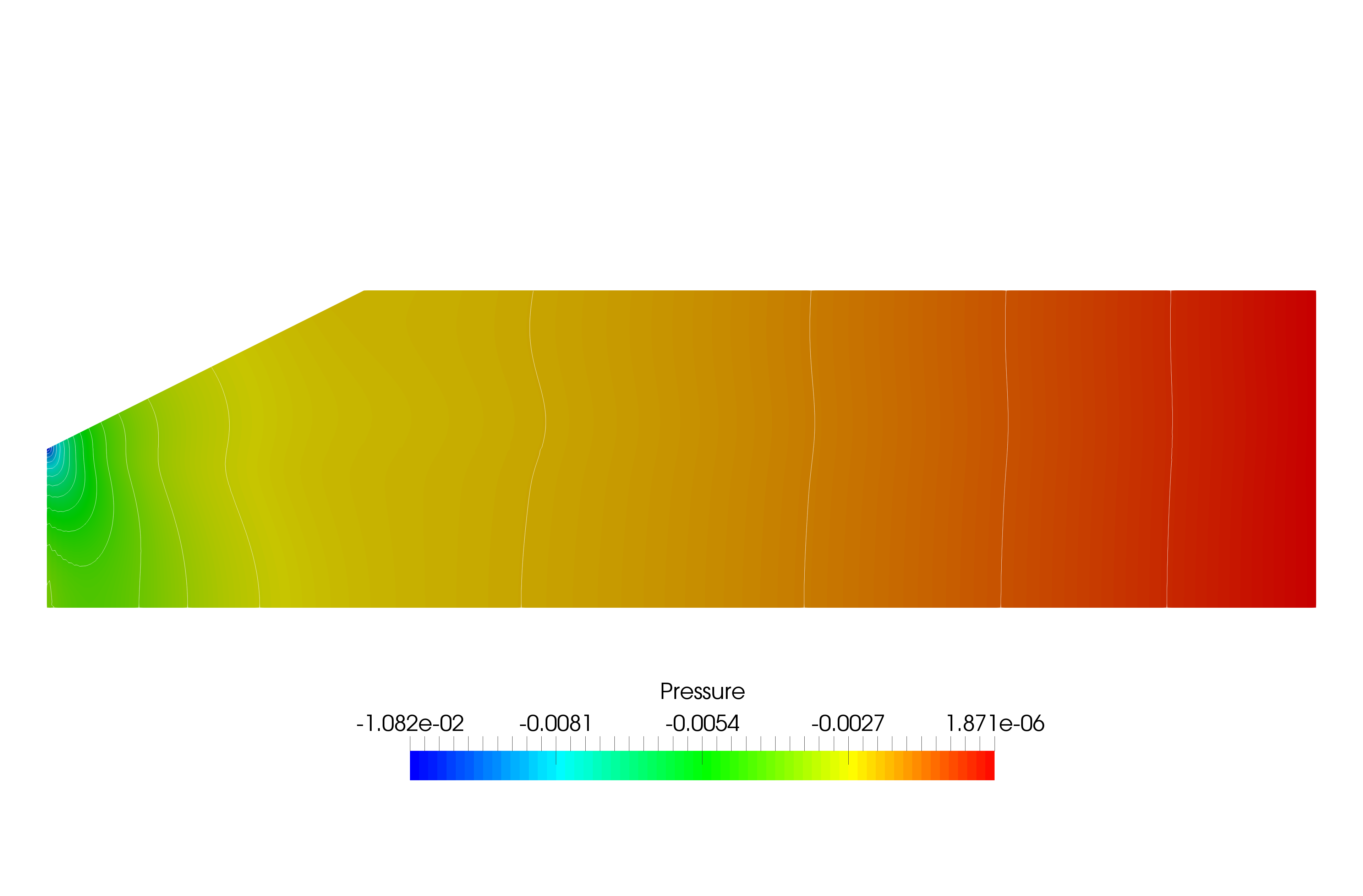}
    \caption{Pressure contours.\label{fig:RZ_cone_1k_pressure_86212_elem}}
  \end{subfigure}
  \caption{Advection-dominated ($\text{Re}=10^3$) axisymmetric channel flow problem
    visualization of (\subref{fig:RZ_cone_1k_streamtraces_86212_elem})~velocity magnitude and streamtraces, and
    (\subref{fig:RZ_cone_1k_pressure_86212_elem})~pressure contours at steady
    state on a mesh with 86212 TRI3 (\pp{1}{1}) elements and 43588 nodes. The region of slow-moving reverse
    flow near the no-slip wall is clearly visible.\label{fig:adv_dom_axisymmetric_flow}}
\end{figure}

\begin{figure}[htpb]
  \centering
  \begin{subfigure}[t]{.9\linewidth}
    \includegraphics[width=\linewidth,viewport=200 400 3600 1800,clip=true]{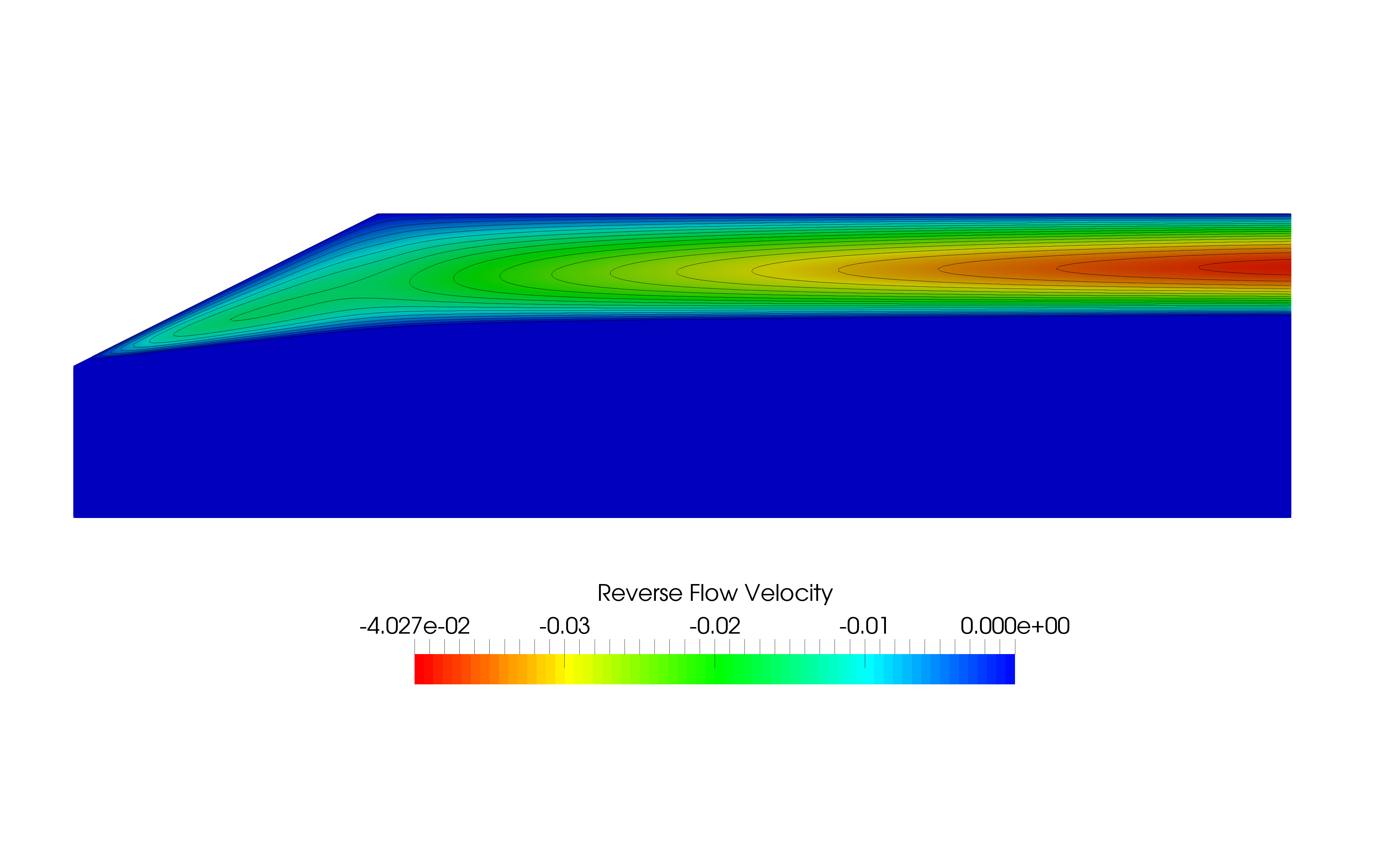}
    \caption{Contours of reverse longitudinal velocity component.\label{fig:RZ_cone_1k_reverse_flow_velocity_86212_elem}}
  \end{subfigure}
  \\[6pt]
  \begin{subfigure}[t]{.9\linewidth}
    \includegraphics[width=\linewidth,viewport=150 0 3600 1500,clip=true]{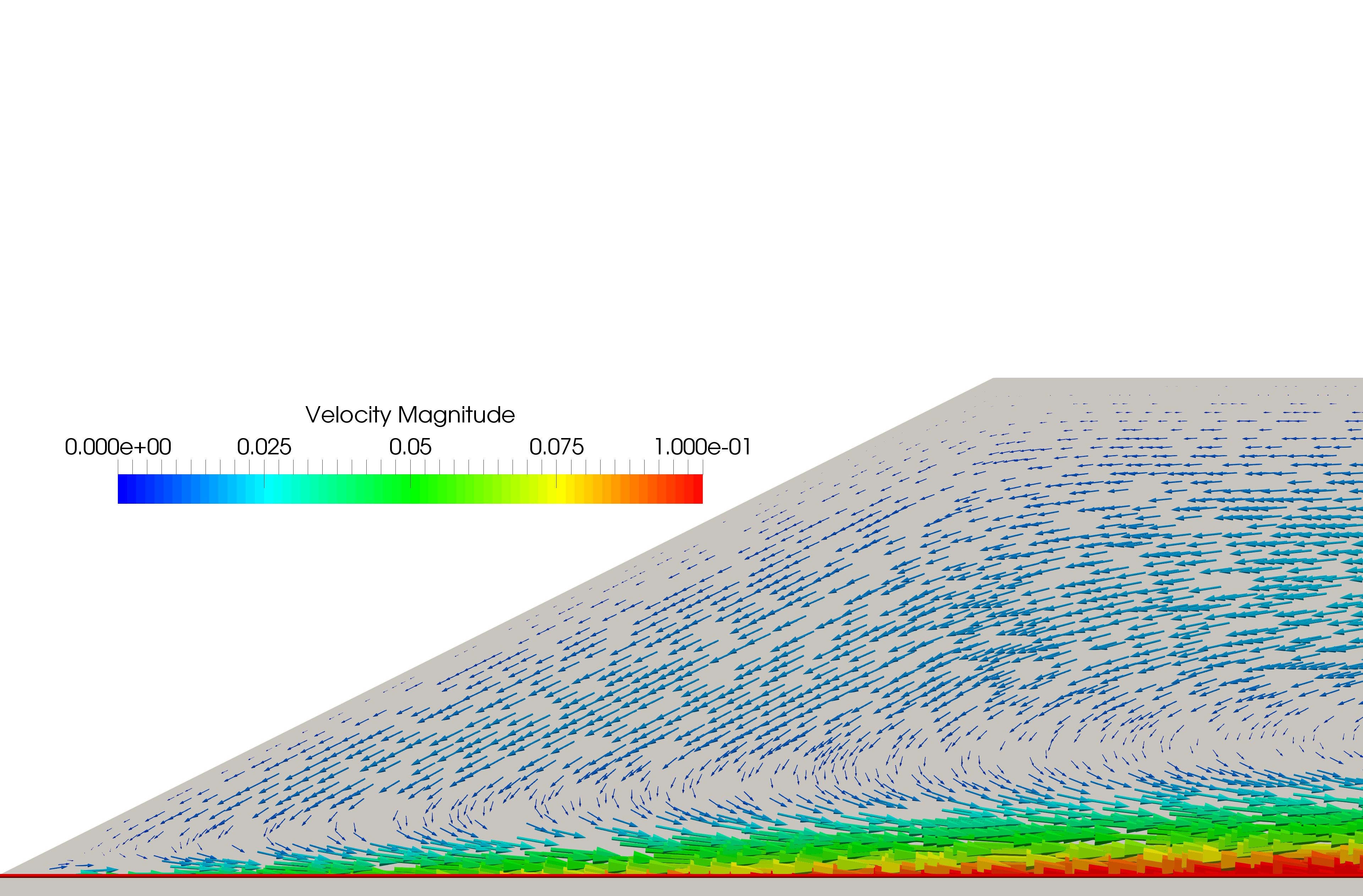}
    \caption{Near-wall velocity vectors.\label{fig:RZ_cone_1k_vectors_86212_elem}}
  \end{subfigure}
  \caption{Advection-dominated ($\text{Re}=10^3$) axisymmetric channel flow problem
    visualization of (\subref{fig:RZ_cone_1k_reverse_flow_velocity_86212_elem})~reverse flow velocity contours, and
    (\subref{fig:RZ_cone_1k_vectors_86212_elem})~velocity vectors near the no-slip wall. Both \gls{SUPG} and \gls{PSPG}
    stabilization are used in this example.\label{fig:adv_dom_axisymmetric_flow_2}}
\end{figure}

The steady state velocity magnitude, streamtraces, and pressure contours for the
advection-dominated ($\text{Re}=10^3$) case are shown in Fig.~\ref{fig:adv_dom_axisymmetric_flow}.
In contrast to the previous case, there is now a large region of slow, recirculating
flow near the channel wall, and an adverse pressure gradient.
In Fig.~\ref{fig:RZ_cone_1k_reverse_flow_velocity_86212_elem}, the region of reversed flow is highlighted
with a contour plot showing the negative longitudinal ($u_z$) velocities, and it seems clear that the downstream boundary condition has a
much larger effect on the overall flow field in this case. In Fig.~\ref{fig:RZ_cone_1k_vectors_86212_elem},
velocity vectors in the upper half of the channel are plotted to aid in visualization of
the reverse-flow region.

\begin{figure}[htpb]
  \centering
  \begin{subfigure}[b]{.9\linewidth}
    \includegraphics[width=\linewidth,viewport=0 0 3800 2000,clip=true]{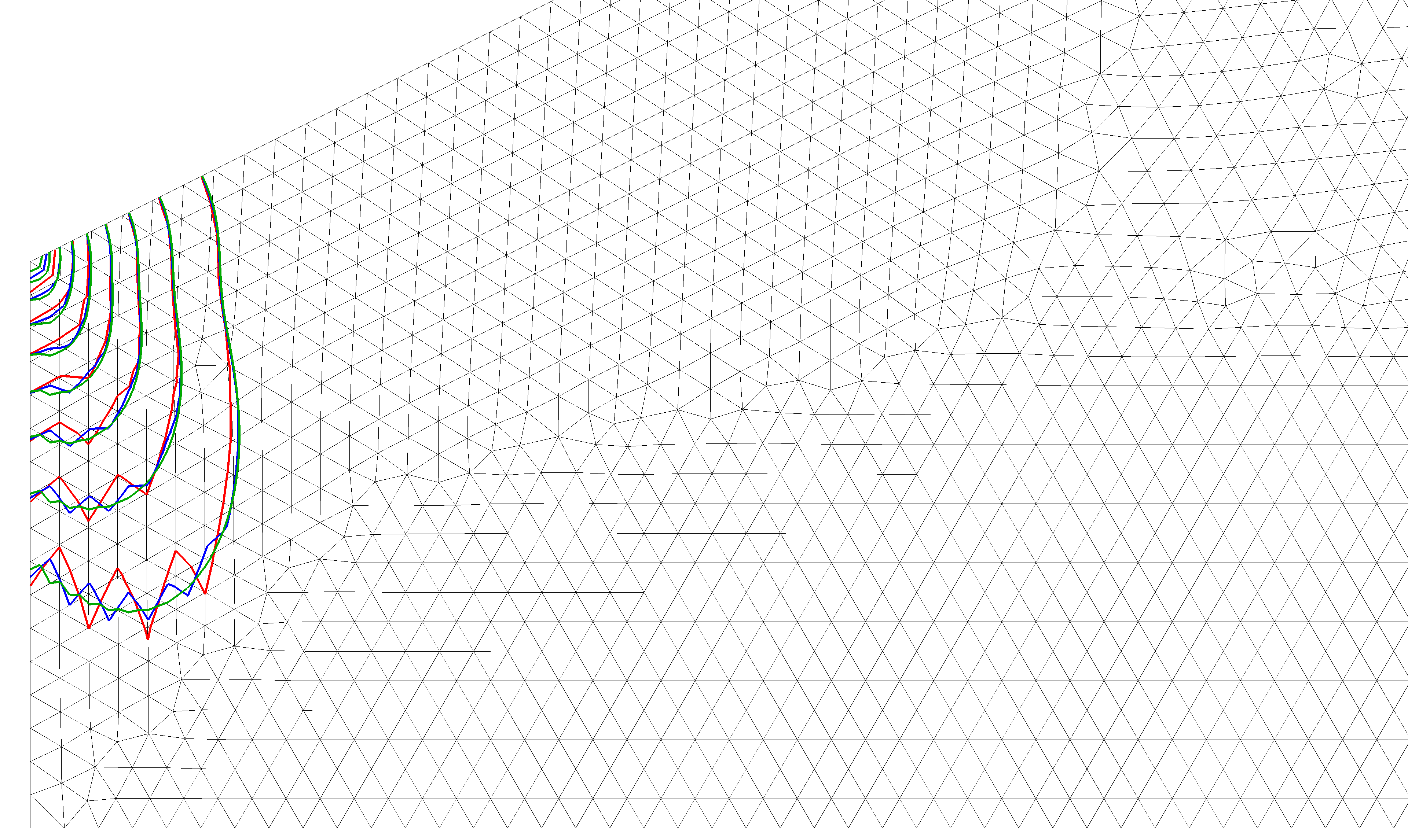}
    \caption{$p$ contours.\label{fig:RZ_cone_1k_zoomed_pressure_comparison_lw6}}
  \end{subfigure}
  \\[6pt]
  \begin{subfigure}[b]{.9\linewidth}
    \includegraphics[width=\linewidth,viewport=0 0 3800 2000,clip=true]{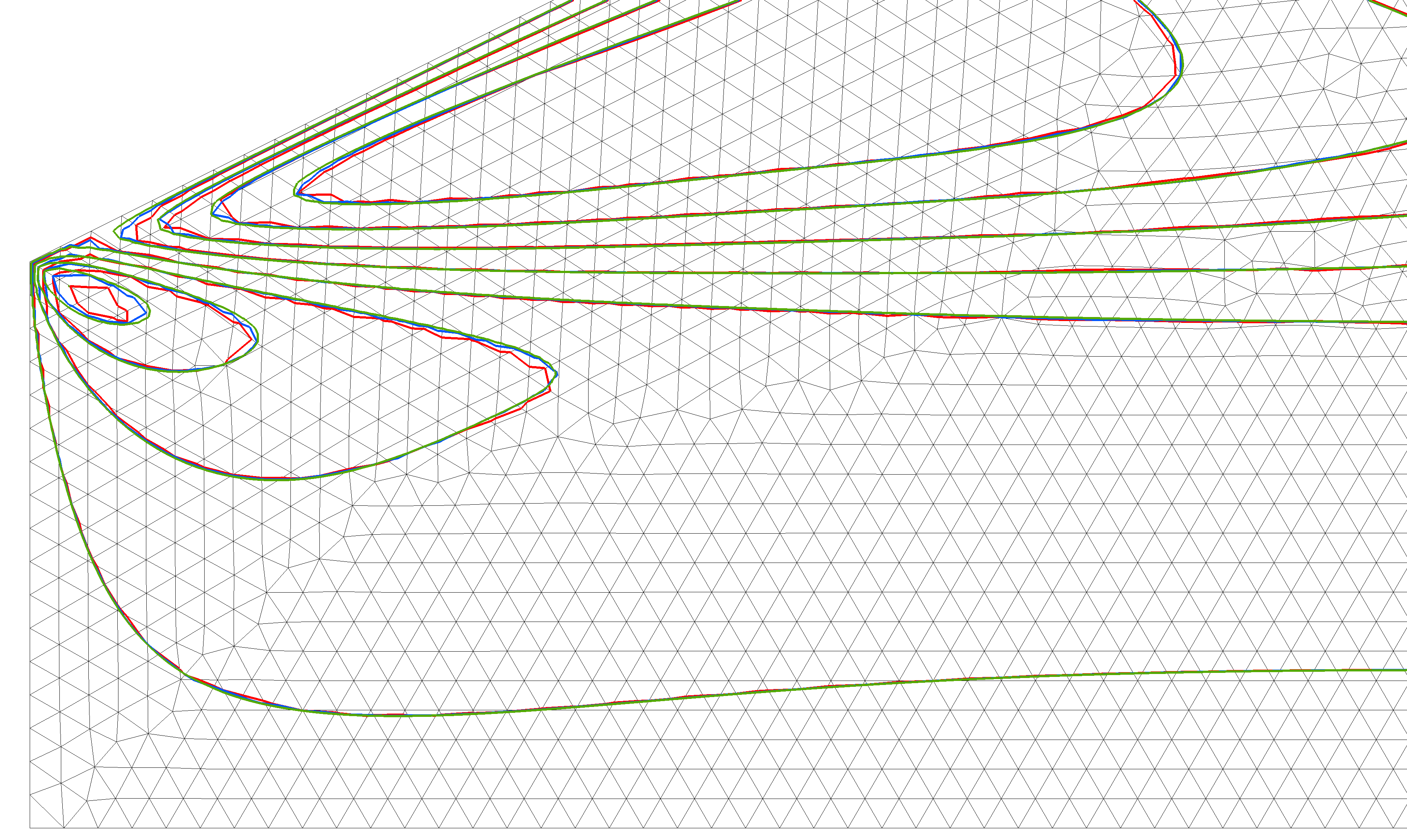}
    \caption{$u_r$ contours.\label{fig:RZ_cone_1k_xvel_comparison}}
  \end{subfigure}
  \caption{Evenly-spaced contours of
    (\subref{fig:RZ_cone_1k_zoomed_pressure_comparison_lw6})~$p \in [-1.08
    \times 10^{-2}, -4.81 \times 10^{-3}]$ and
    (\subref{fig:RZ_cone_1k_xvel_comparison})~$u_r\in [-7.32\times
    10^{-3}, 6.93\times 10^{-3}]$ on meshes with 9456 (red), 21543
    (blue), and 86212 (green) elements. Non-physical oscillations can be
    observed in several of the pressure contours near the sharp inlet
    corner on the coarser grids, but they are mostly absent in the fine
    grid results. The velocity contours are all fairly comparable except
    near the separation point on the upper wall. The background mesh shown
    corresponds to the coarsest grid
    level.\label{fig:RZ_cone_1k_pressure_comparison_lw6}}
\end{figure}

The pressure singularity is still an important characteristic of the
flow in the advection-dominated case, and to investigate the error in
the finite element solution near this point a small mesh convergence
study was conducted on a sequence of grids with 9456, 21543,
and 86212 elements. (The results in Figs.~\ref{fig:adv_dom_axisymmetric_flow}
and~\ref{fig:adv_dom_axisymmetric_flow_2} are from the finest grid in
this sequence.) Contour plots comparing the pressure and radial
velocity near the sharp inlet corner on the different grids are given in
Figs.~\ref{fig:RZ_cone_1k_zoomed_pressure_comparison_lw6}
and~\ref{fig:RZ_cone_1k_xvel_comparison}, respectively. On the coarse
grid (red contours) in particular, we observe pronounced non-physical
oscillations in the pressure field near the singularity. These
oscillations are much less prominent on the finest grid (green
contours). The radial velocity contours are comparable on all three
meshes, with the largest discrepancy occurring on the channel wall
where the boundary layer flow separation initially occurs.

\subsection{Flow over a sphere\label{sec:sphere}}
\begin{figure}[htpb]
  \centering
  \begin{subfigure}[t]{.9\linewidth}
    \includegraphics[width=\linewidth]{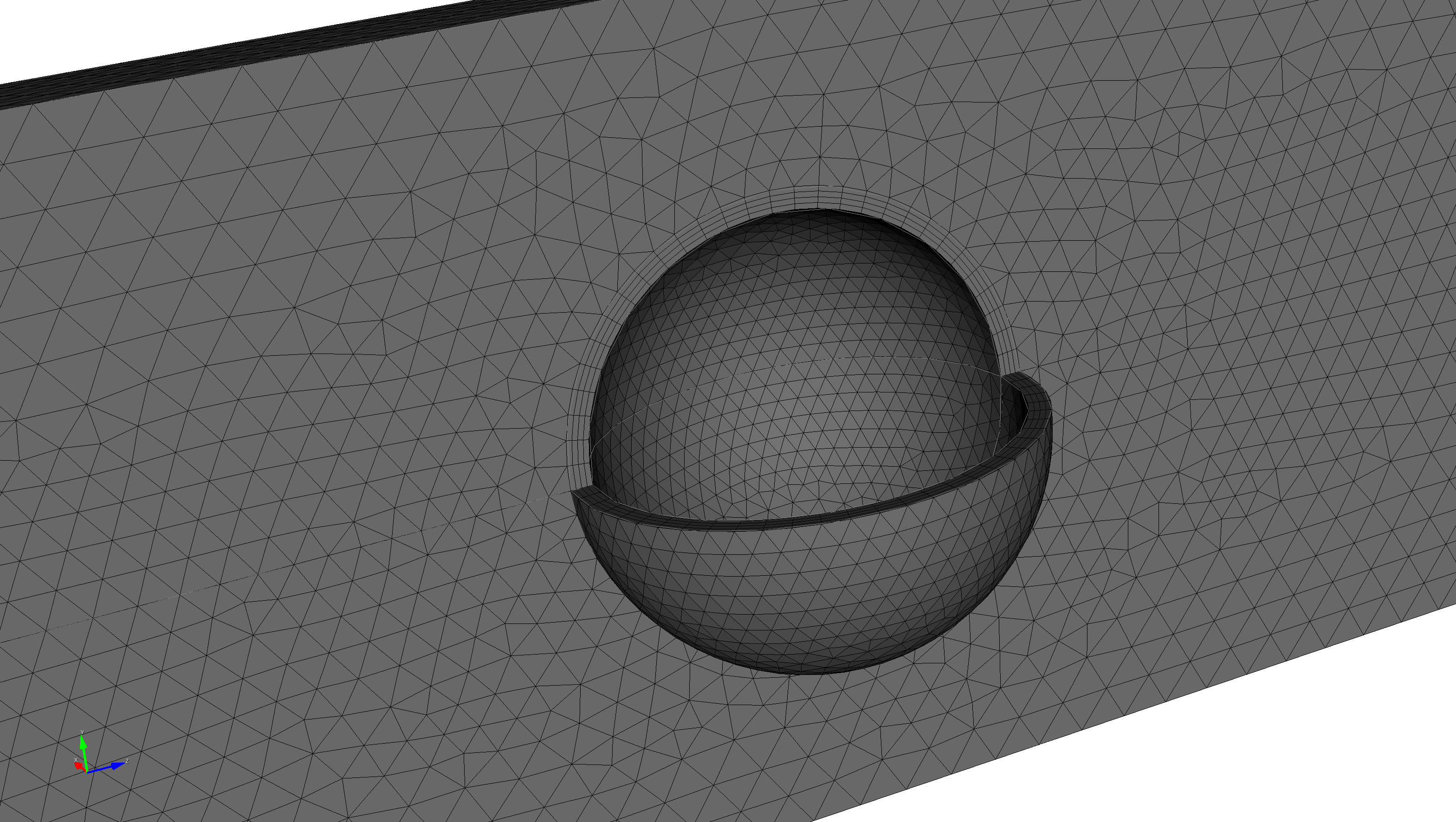}
    \caption{Details of the hybrid mesh surrounding the sphere.\label{fig:sphere_mesh}}
  \end{subfigure}
  \\[6pt]
  \begin{subfigure}[t]{.9\linewidth}
    \includegraphics[width=\linewidth]{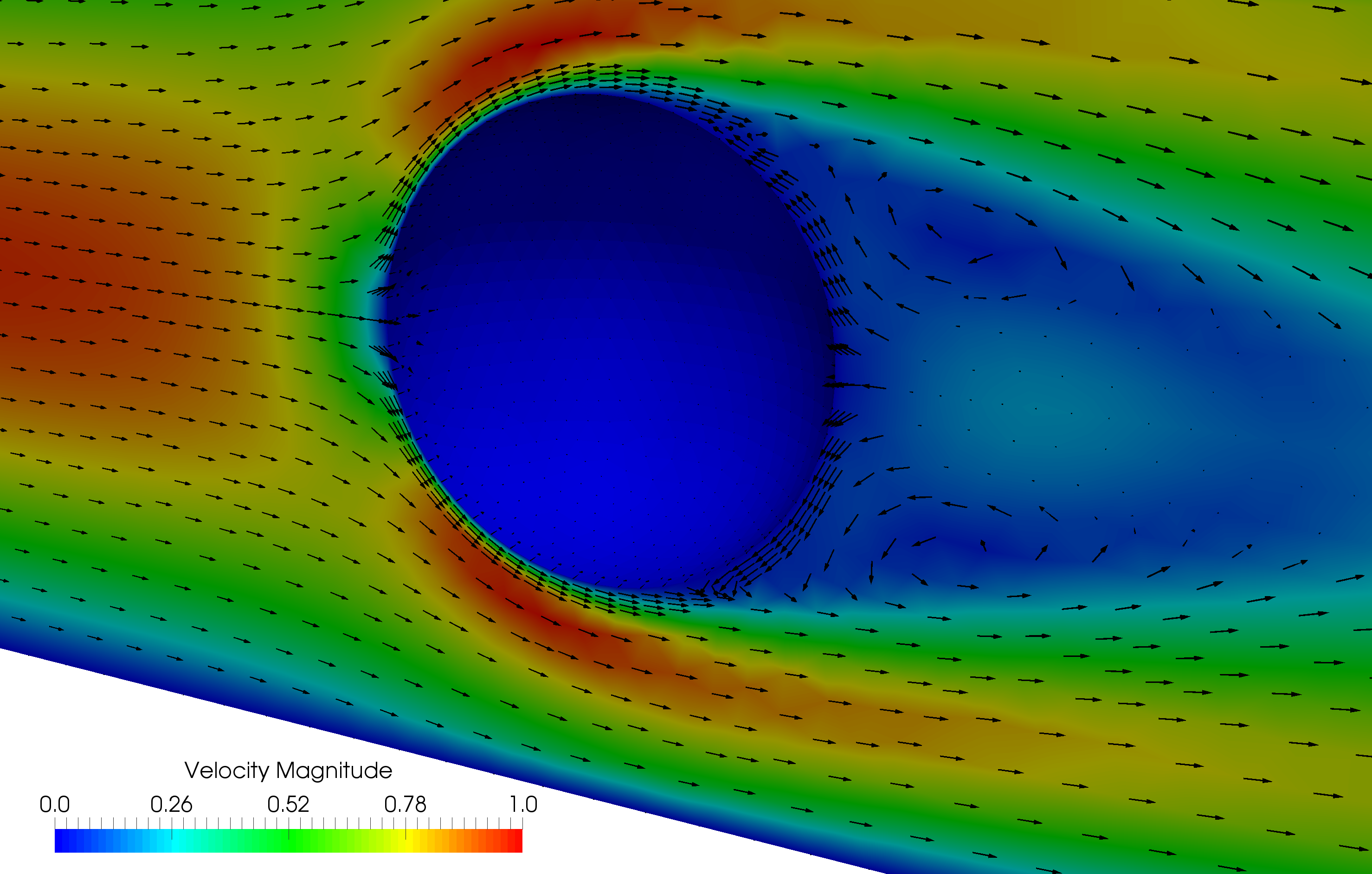}
    \caption{Contour plot of velocity magnitude and (unscaled) velocity vectors.\label{fig:sphere_unscaled}}
  \end{subfigure}
  \caption{Plots showing  (\subref{fig:sphere_mesh}) details of the hybrid prismatic/tetrahedral mesh near the sphere, and
    (\subref{fig:sphere_unscaled}) contour plot of the centerline velocity magnitude combined with unscaled
    velocity vectors to indicate the direction of the flow field. The approximate location of the boundary layer
    separation point on the downstream side of the sphere is visible, as are the extent of the recirculation
    zones in the wake.\label{fig:sphere_mesh_and_boundary_layer}}
\end{figure}

In this example, we simulate flow in the channel $\Omega=[-2,2]^2 \times [-5,5]$ around a sphere of diameter 2
centered at the origin. The inlet flow profile is given by:
\begin{align}
  u_1 &= u_2 = 0
  \\
  \label{eq:inlet_profile_sphere}
  u_3 &= \frac{(4-x^2)(4-y^2)}{16}
\end{align}
The Reynolds number based on the maximum inlet velocity, sphere
diameter, and viscosity $\mu = 4\times 10^{-3}$ is $5\times 10^2$.  (Based on the
area-averaged inlet velocity, the Reynolds number is $\approx
222.22$.) A 3D hybrid mesh, pictured in Fig.~\ref{fig:sphere_mesh},
with 104069 elements (93317 tetrahedra and 10752 prisms) comprising 23407
nodes, which is locally graded near the boundary of the sphere,
is used for the calculation.
A transient, \gls{SUPG}/\gls{PSPG}-stabilized Laplace formulation of the \gls{INS} equations
on equal-order elements (nominally \pp{1}{1}, although the prismatic elements
also contain bilinear terms in the sphere-normal direction) is employed, and
the equations are stepped to steady state using an implicit (backward Euler)
time integration routine with adaptive timestep selection based on the number
of nonlinear iterations.

Natural boundary conditions are once again imposed at the
outlet. Because of the proximity of the channel walls and the outlet
to the sphere, we do not necessarily expect to observe vortex shedding
or non-stationary steady states even at this relatively high Reynolds
number.  The problem was solved in parallel on a workstation with 24
processors using Newton's method combined with a parallel direct
solver (SuperLU\_DIST~\cite{Xiaoye_2003}) at each linear iteration, although
an inexact Newton scheme based on a parallel Additive Schwarz
preconditioner with local ILU sub-preconditioners was also found to
work well.

\begin{figure}[htpb]
  \centering
  \begin{subfigure}[t]{.9\linewidth}
    \includegraphics[width=\linewidth]{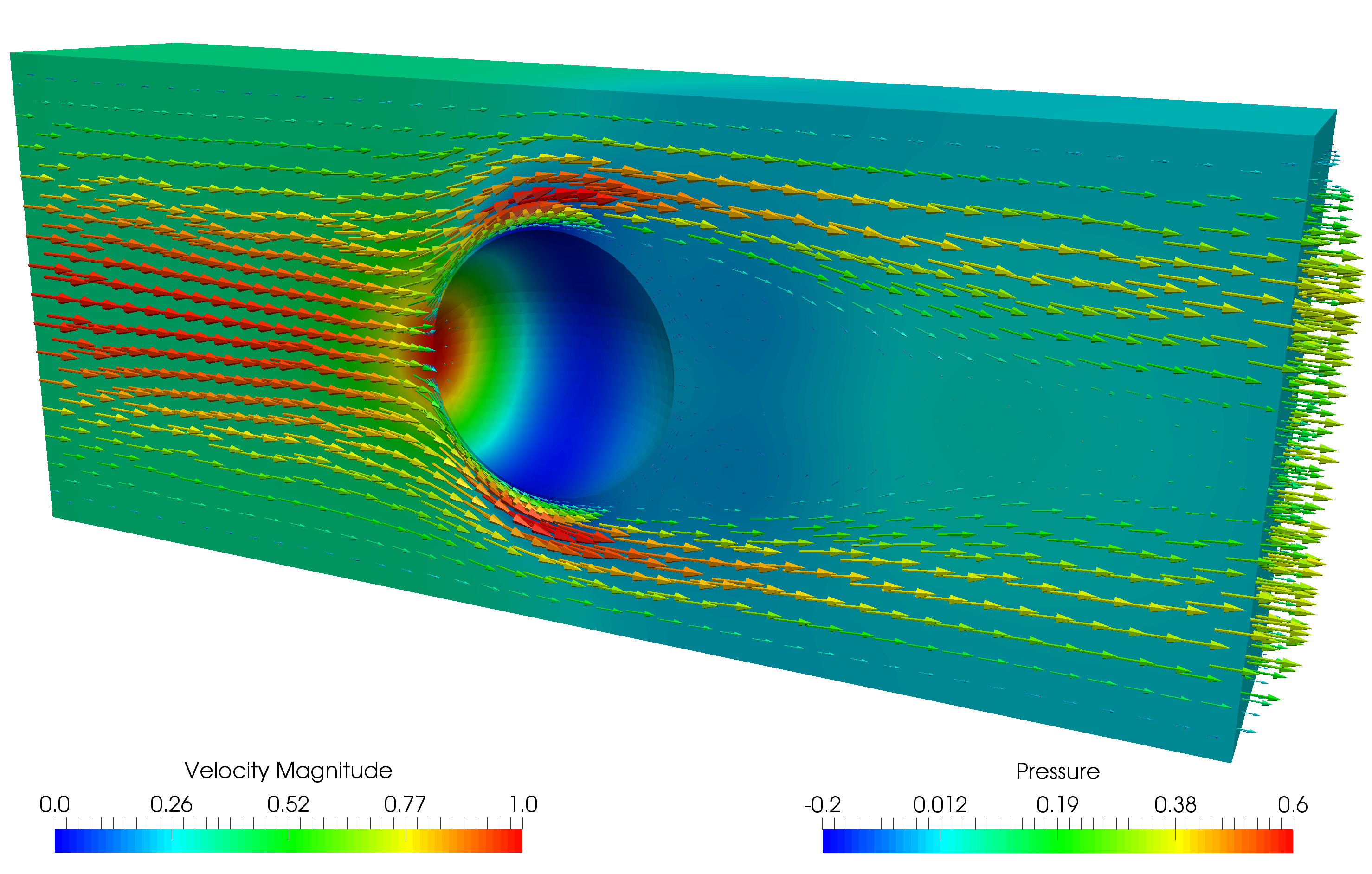}
    \caption{Centerline pressure contour plot and velocity vectors colored by magnitude.\label{fig:sphere_vectors}}
  \end{subfigure}
  \begin{subfigure}[t]{.9\linewidth}
    \includegraphics[width=\linewidth]{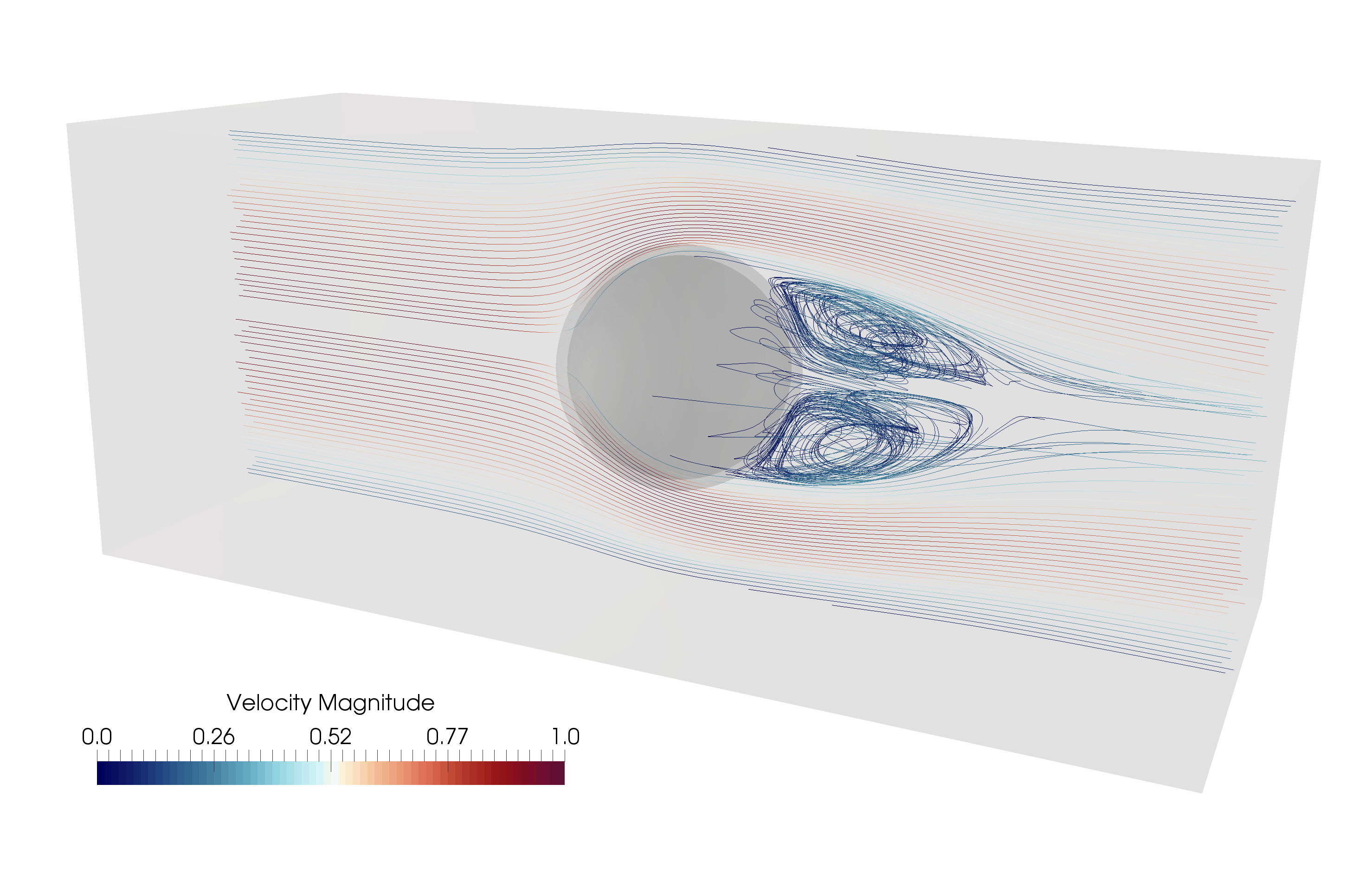}
    \caption{Streamtraces showing entrainment behind the sphere.\label{fig:sphere_streamtraces}}
  \end{subfigure}
  \caption{Plots showing (\subref{fig:sphere_vectors}) the steady state centerline pressure and velocity vectors colored by $|\vec{u}|$, and
    (\subref{fig:sphere_streamtraces}) streamtraces originating from just downstream of the sphere and colored
    by the velocity magnitude. In (\subref{fig:sphere_vectors}) the vector length is proportional to $|\vec{u}|$, so it is
    difficult to see the three-dimensional nature of the region of recirculating flow directly behind the sphere, but this region is clearly visible in
    (\subref{fig:sphere_streamtraces}).\label{fig:sphere}}
\end{figure}

A contour plot of the centerline velocity magnitude and the
corresponding (unscaled) velocity vector field at steady state is shown
in~Fig.~\ref{fig:sphere_unscaled}. In this view, the cone-shaped recirculation
region behind the sphere is visible, as is the rough location of the boundary
layer separation point on the back half of the sphere.
The steady state velocity and pressure fields along the
channel centerline are shown in Fig.~\ref{fig:sphere_vectors}. In this figure, the
velocity vectors are colored and scaled by the velocity magnitude, so the
region of slow-moving flow in the wake of the sphere is once again visible.
Finally, in Fig.~\ref{fig:sphere_streamtraces}, streamtraces originating from
a line source placed directly behind the sphere are shown. The streamtraces
help to visualize the three-dimensional nature of the wake region, and the
entrainment of particles from the surrounding flow.

\section{Conclusions and future work\label{sec:future}}
The stabilized finite element formulation of the \gls{INS} equations
in the \ns module of \gls{MOOSE} is currently adequate for solving
small to medium-sized problems over a wide range of
Reynolds numbers, on non-trivial 2D and 3D geometries. The software
implementation of the \gls{INS} equations has been verified for a
number of simple test problems, and is developed under a modern and
open continuous integration workflow to help guarantee high software
quality.  The \ns module is also easy to integrate (without writing
additional code) with existing and new \gls{MOOSE}-based application
codes since it is distributed directly with the \gls{MOOSE} framework
itself: all that should be required is a simple change to the
application Makefile.

There are still a number of challenges and opportunities for
collaboration to improve the quality and effectiveness of
the \gls{INS} implementation in the \ns module. These opportunities
range from the relatively straightforward, such as implementing
\gls{LSIC} stabilization, investigating and developing non-conforming
velocity, discontinuous pressure, and divergence-free finite element discretizations, and adding
\gls{MOOSE} \texttt{Actions} to simplify and streamline the input file
writing experience, to more involved tasks like implementing a
coupled, stabilized temperature convection-diffusion equation and
Boussinesq approximation term, implementing and testing turbulence
models, and adding support for simulating non-Newtonian fluids.  The
closely-related low Mach number equations~\cite{Rehm_1978} are also of
theoretical and practical interest for cases where temperature
variations in the fluid are expected to be large, and would be
appropriate to include in the \ns module.

The largest outstanding challenges primarily involve improving the
basic field-split preconditioning capabilities which are already
present in the \ns module, but are not currently widely used due to
their cost and the lack of a general and efficient procedure for
approximating the action of the Schur complement matrix. In this
context, we also note the important role of multigrid methods (both
geometric and algebraic variants) in the quest for developing scalable
algorithms. The field-split approach, combined with geometric
multigrid preconditioners, is probably the most promising avenue for
applying the \ns module to larger and more realistic applications.
These techniques require a tight coupling between the finite element
discretization and linear algebra components of the code, and a deep
understanding of the physics involved. We are confident that the
\ns module itself will provide a valuable and flexible test bed
for researching such methods in the future.

\section*{Acknowledgments}
The authors are grateful to the MOOSE users and members of the open
science community that have contributed to the \ns framework with
code, feature requests, capability clarifications, and insightful
mailing list questions, and for the helpful comments from the
anonymous reviewers which helped to improve this paper.

This manuscript has been authored by a contractor of the
U.S. Government under Contract DE-AC07-05ID14517. Accordingly, the
U.S. Government retains a non-exclusive, royalty-free license to
publish or reproduce the published form of this contribution, or allow
others to do so, for U.S. Government purposes.

\clearpage
\printglossary[type=\acronymtype]
\bibliographystyle{ieeetr}
\bibliography{ins}
\end{document}